\documentclass[oneside, a4paper,reqno]{amsart}
\usepackage{pdfsync}
\usepackage{stmaryrd}
\usepackage{mathrsfs}
\usepackage{multicol}
\usepackage{amsmath, amsthm, amscd, amssymb, latexsym, eucal}
\usepackage[all]{xy}
 \calclayout \makeatletter
 \calclayout \makeatletter
\def\serieslogo@{} \def\@setcopyright{} \makeatother

\usepackage{multienum}

\usepackage{hyperref}
\usepackage{color}
\usepackage{cite}

\hypersetup{colorlinks=true,
     breaklinks=true,
     linkcolor=blue,    
     bookmarks=true,
      pageanchor=true}

\makeatletter
\renewcommand*\env@matrix[1][c]{\hskip -\arraycolsep
  \let\@ifnextchar\new@ifnextchar
  \array{*\c@MaxMatrixCols #1}}
\makeatother

\numberwithin{equation}{section}
\newtheorem{thm}{Theorem}[section]
\newtheorem{cor}[thm]{Corollary}
\newtheorem{lem}[thm]{Lemma}
\newtheorem{prop}[thm]{Proposition}

\theoremstyle{definition}
\newtheorem{defn}[thm]{Definition}
\newtheorem{rem}[thm]{Remark}
\newtheorem{exam}[thm]{Example}

\newtheorem*{thmA}{Theorem~A}
\newtheorem*{thmB}{Theorem~B}

\newtheorem*{ackn}{Acknowledgment}

\newcommand{\lxr}{\longrightarrow}

\newcommand{\A}{\mathscr A}
\newcommand{\B}{\mathscr B}
\newcommand{\C}{\mathscr C}

\newcommand{\D}{\mathscr D}
\newcommand{\E}{\mathscr E}
\newcommand{\F}{\mathcal F}
\newcommand{\G}{\mathcal G}

\newcommand{\M}{\mathcal M}

\newcommand{\mQ}{\mathsf{Q}}

\newcommand{\mf}{\mathsf{F}}

\newcommand{\mt}{\mathsf{T}}
\newcommand{\mh}{\mathsf{H}}
\newcommand{\mU}{\mathsf{U}}
\newcommand{\mz}{\mathsf{Z}}

\newcommand{\mD}{\mathsf{D}}

 \DeclareMathOperator{\inc}{\mathsf{inc}}
  
\DeclareMathOperator*{\Ker}{\mathsf{Ker}}
 \DeclareMathOperator*{\Image}{\mathsf{Im}}
\DeclareMathOperator*{\Coker}{\mathsf{Coker}}

 \DeclareMathOperator{\pd}{\mathsf{pd}}
\DeclareMathOperator*{\id}{\mathsf{id}}

  \DeclareMathOperator*{\gld}{\mathsf{gl.dim}}

\DeclareMathOperator*{\Mod}{\mathsf{Mod}-\!}

\DeclareMathOperator*{\End}{\mathsf{End}}
 \DeclareMathOperator*{\smod}{\mathsf{mod}-\!}

\DeclareMathOperator*{\inj}{\mathsf{inj}}
\DeclareMathOperator*{\proj}{\mathsf{proj}}
\DeclareMathOperator*{\Inj}{\mathsf{Inj}}
\DeclareMathOperator*{\Proj}{\mathsf{Proj}}

\newcommand{\GProj}{\operatorname{\mathsf{GProj}}\nolimits}
 \newcommand{\Gproj}{\operatorname{\mathsf{Gproj}}\nolimits}

\DeclareMathOperator*{\add}{\mathsf{add}}

\DeclareMathOperator{\Hom}{\mathsf{Hom}}

\DeclareMathOperator*{\Ext}{\mathsf{Ext}}
\DeclareMathOperator*{\Tor}{\mathsf{Tor}}

  \DeclareMathOperator*{\op}{\mathsf{op}}
   
  \DeclareMathOperator*{\silp}{\mathsf{silp}}
  \DeclareMathOperator*{\spli}{\mathsf{spli}}
  
   \DeclareMathOperator*{\Ab}{\A\!\textit{b}}

   \DeclareMathOperator*{\du}{\mathsf{D}}

\DeclareMathOperator*{\Cok}{\mathsf{Cok}}

\DeclareMathOperator*{\Mono}{\mathsf{Mono}}
\DeclareMathOperator*{\mono}{\mathsf{mono}}

\DeclareMathOperator*{\Mor}{\mathsf{Mor}}
\DeclareMathOperator*{\DMor}{\mathsf{DMor}}
\DeclareMathOperator*{\Mon}{\mathsf{Mon}}
\DeclareMathOperator*{\mon}{\mathsf{mon}}

\newcommand{\iso}{\cong}

\newcommand{\iden}{\operatorname{Id}\nolimits}

\newsavebox{\proofbox}
\savebox{\proofbox}{\begin{picture}(7,7)%
  \put(0,0){\framebox(7,7){}}\end{picture}}

\begin{document}


\title[]{
Gorenstein homological aspects of \\ monomorphism categories via Morita rings
}

\author[C. Psaroudakis]{Nan Gao and Chrysostomos Psaroudakis}
\address{Nan Gao \\ Department of Mathematics, Shanghai University, Shanghai 200444, PR China}
\email{nangao@shu.edu.cn}

\address{Chrysostomos Psaroudakis\\
Department of Mathematical Sciences, Norwegian University of Science and Technology \\
 7491 Trondheim, Norway} 
\email{chrysostomos.psaroudakis@math.ntnu.no}

\date{\today}




\begin{abstract} 
For any ring $R$ the category of monomorphisms is a full subcategory of the morphsim category over $R$, where the latter is equivalent to the module category of the triangular matrix ring with entries the ring $R$. In this work, we consider the monomorphism category as a full subcategory of the module category over the Morita ring with all entries the ring $R$ and zero bimodule homomorphisms. This approach provides an interesting link between Morita rings and monomorphism categories. The aim of this paper is two-fold. First, we construct Gorenstein-projective modules over Morita rings with zero bimodule homomorphisms and we provide sufficient conditions for such rings to be Gorenstein Artin algebras. This is the first part of our work which is strongly connected with monomorphism categories. In the second part, we investigate monomorphisms where the domain has finite projective dimension. In particular, we show that the latter category is a Gorenstein subcategory of the monomorphism category over a Gorenstein algebra. Finally, we consider the category of coherent functors over the stable category of this Gorenstein subcategory and show that  it carries a structure of a Gorenstein abelian category.
\end{abstract}

\maketitle

\setcounter{tocdepth}{1} \tableofcontents

\section{Introduction and Main Results}
This article deals with Gorenstein homological aspects of Morita rings and monomorphism categories. We highlight the connection between monomorphism categories and Morita rings and we explain how this setting is the proper framework for the aims and objectives of this paper. In what follows, we first give some background, motivation and indicate how these concepts are connected from our perspective.

For an abelian category $\A$ we denote by ${\Mor{\A}}$ the category of morphisms over $\A$. The monomorhism category $\Mon{\A}$ is by definition the full subcategory of $\Mor{\A}$ consisting of all monomorphisms in $\A$. If $R$ is a ring and $\A$ is the category $\Mod{R}$ of $R$-modules, then the category of monomorphisms $\Mon{(\Mod{R})}$ can be considered as a full subcategory of the module category $\Mod{\mathsf{T}_2(R)}$, where $\mathsf{T}_2(R)=\bigl(\begin{smallmatrix}
R & R \\
0 & R
\end{smallmatrix}\bigr)$, since it is known that the morphism category $\Mor{(\Mod{R})}$ is equivalent to $\Mod{\mathsf{T}_2(R)}$. Note that $\Mon{\A}$ is an extension closed subcategory of the abelian category $\Mor{\A}$ and therefore is an exact category in the sense of Quillen. Monomorphism categories appear quite naturally in various settings and are omnipresent in representation theory. In fact, there are connections with classification problems (Ringel, Schmidmeier \cite{RS1, RS2, RS3}, Xiong, Zhang, Zhang \cite{XZZ}), with weighted projective lines (Kussin, Lenzing, Meltzer \cite{KLM}), and with aspects of Gorenstein homological algebra (Beligiannis \cite{Bel:virtually, Bel:finiteCMtype}, Zhang \cite{AZ}, Chen \cite{Chen}). In a series of papers \cite{XiongZhang, LiZhang, Zhang}, the authors determined the Gorenstein-projective modules over the triangular matrix algebra $\bigl(\begin{smallmatrix}
A & N \\
0 & B
\end{smallmatrix}\bigr)$ under some conditions on the bimodule $_AN_B$. In the special case where $A=N=B$ and $A$ is a Gorenstein Artin algebra, i.e. $\mathsf{T}_2(A)=\bigl(\begin{smallmatrix}
A & A \\
0 & A
\end{smallmatrix}\bigr)$, the authors in \cite{LiZhang} showed that a module over $\mathsf{T}_2(A)$, i.e. a triple $(X,Y,f)$, is Gorenstein-projective if and only if $(X,Y,f)$ belongs to the monomorphism category $\mon(A)$ and the modules $X$, $Y$ and $\Coker{f}$ are Gorenstein-projective. 

A natural extension of triangular matrix rings is the class of Morita rings with zero bimodule homomorphisms. First, recall that Morita rings are $2\times 2$ matrix rings associated to Morita contexts (Bass \cite{Bass_Morita}, Cohn \cite{Cohn}). A {\em Morita context} over two unital associative rings $R$ and $S$, is a tuple $\M = (R,N,M,S, \phi, \psi)$, where  ${_{S}}M_{R}$ is an $S$-$R$-bimodule,  ${_{R}}N_{S}$ is an $R$-$S$-bimodule, $\phi \colon   M\otimes_{R}N \lxr S$ is an $S$-$S$-bimodule homomorphism, and $\psi \colon N\otimes_{S}M \lxr R$ is an $R$-$R$-bimodule homomorphism, satisfying certain associativity conditions. 
If $\M$ is a Morita context, there is the 
associative ring $\Lambda_{(\phi,\psi)}=\bigl(\begin{smallmatrix}
A & _AN_B \\
_BM_A & B
\end{smallmatrix}\bigr)$ called the {\em Morita ring} of $\M$, where $\Lambda_{(\phi,\psi)} = A \oplus N \oplus M \oplus B$ as an abelian group and the matrix multiplication is given by the ordinary operation on matrices using the maps $\phi$ and $\psi$, see subsection~\ref{subsectionMoritarings}. In this paper the Morita rings that we deal with admit the structure of an Artin algebra. We state this as: Morita rings which are Artin algebras, in order to distinguish it from the notion of Morita algebras due to Kerner-Yamagata \cite{KernerYamagata}. We refer to \cite{McConnel, Rowen} for the terminology of Morita rings, and to \cite[Introduction]{GP}, among other papers, for a thorough discussion of Morita rings as well as examples and situations where Morita rings appear.

A particular case of interest for us is when $A=B=N=M=R$, where $R$ is an associative unital ring and $\phi=\psi=0$. This class of matrix rings provides on one hand a natural extension of the triangular matrix ring $\mathsf{T}_2(R)$ and on the other provides the link between monomorphsim categories and Morita rings. We now explain this.
For an abelian category $\A$ we define the {\em double morphism category} $\DMor(\A)$, see subsection~\ref{subMonocat}, whose objects are pairs of morphisms in $\A\colon$$\xymatrix@C=0.5cm{
X  \ar@<-.7ex>[rr]|-{f} && Y \ar@<-.7ex>[ll]|-{g }}$ such that $f\circ g=0$ and $g\circ f=0$ and morphisms are commutative squares defined in a natural way. In particular, if $R$ is a ring then the double morphism category $\DMor(\Mod{R})$ of $\Mod{R}$ is equivalent to the category of modules over the Morita ring $\Delta_{(0,0)}=\bigl(\begin{smallmatrix}
R & R \\
R & R
\end{smallmatrix}\bigr)$ with zero bimodule homomorphisms. Consider now the upper triangular matrix ring $\mathsf{T}_2(R)$. The monomorphism category $\Mon(R)$ is a full subcategory of $\Mod{\mathsf{T}_2(R)}$ defined by $\Mon(R)=\{(X,Y,f) \ | \ f\colon X\lxr Y \ \ \text{is a monomorphism}\}$. As mentioned already, a natural extension of the ring $\mathsf{T}_2(R)$ is the Morita ring $\Delta_{(0,0)}$, since there is a full embedding $\Mod{\mathsf{T}_2(R)}\lxr \Mod{\Delta_{(0,0)}}$ (see subsection~\ref{subringepi}). Then we define the monomorphism category $\Mono(R)$ as the full subcategory of $\Mod{\Delta_{(0,0)}}$ defined by all tuples $(X,Y,f,0)$ such that $f\colon X\lxr Y$ is a monomorphism. It follows easily that the exact categories $\Mon(R)$ and $\Mono(R)$ are equivalent. However, there are many properties which are not satisfied when we view the monomorphism category inside the double morphism category. 
From the discussion so far, one of the main problems considered in this paper is as follows$\colon$ 

\smallskip

\noindent {\bf \textsf{Problem}}$\colon$Construct Gorenstein-projective modules over Morita rings with zero bimodule homomorphisms. 

\smallskip

This problem provides the link between the monomorhism categories and Morita rings as explained above. In particular, we obtain in Section~\ref{GPmodulesoverMorita} families of Gorenstein-projective modules over $\Delta_{(0,0)}$ which belong to the monomorphism category $\Mono(R)$. The situation mentioned above and the important role that monomorphism categories as well as Morita rings play in various different contexts, provide a strong motivation for studying monomorphism categories in the context of Morita rings using homological and representation-theoretic tools.  Our aim in this paper is two-fold and can be summarized as follows$\colon$

\smallskip

(i) Construct Gorenstein-projective modules over Morita rings with zero bimodule homomorphisms and provide sufficient conditions for such rings to be Gorenstein Artin algebras.

(ii) Construct Gorenstein abelian categories from exact subcategories of the monomorphism category.

\smallskip

The organization and the main results of the paper are as follows. In Section~\ref{sectionMorringsMonocat} we collect preliminary notions and results on Morita rings and monomorphism categories that will be useful throughout the paper and we fix notation. Moreover we introduce the double morphism category of an abelian category and we define the monomorphism category in this general setting. The rest of the paper is divided into two parts that we describe briefly below and we present the main results.

The first part of our paper deals with item (i) and consists of Sections~\ref{GPmodulesoverMorita} and~\ref{HomolEmbGorAlgebras}. For an algebra $\Lambda$, we denote by $\Gproj{\Lambda}$ the full subcategory of $\smod{\Lambda}$ consisting of the finitely generated Gorenstein-projective $\Lambda$-modules. Let $\Lambda_{(0,0)}=\bigl(\begin{smallmatrix}
A & _AN_B \\
_BM_A & B
\end{smallmatrix}\bigr)$ be a Morita ring which is an Artin algebra and has zero bimodule homomorphisms. In order to construct Gorenstein-projective modules over $\Lambda_{(0,0)}$, we need to assume some natural conditions on the bimodules $_AN_B$ and $_BM_A$, similar to the conditions considered by Zhang \cite{Zhang} in the triangular matrix case. We refer to these assumptions as the {\em compatibility conditions} on $_BM_A$ and $_AN_B$, see Section~\ref{GPmodulesoverMorita}. These conditions have a nice interpretation via finiteness of the projective dimension of the bimodules $N$ and $M$, see Corollary~\ref{corcompatibilitycondpdfinite}.
Our first main result is Theorem~A (i), which provides a method to construct Gorenstein-projective modules over Morita rings with zero bimodule homomorphisms. We refer to Theorem~\ref{thmGorproj} for the proof as well as its dual version.
On the other hand, we give sufficient conditions for a Morita ring $\Lambda_{(0,0)}$ with zero bimodule homomorphisms to be a Gorenstein Artin algebra. This constitutes our second main result and is Theorem~A (ii), see Theorem~\ref{thmGorenstein}. Recall that an Artin algebra $\Lambda$ is Gorenstein if and only if $\spli{\Lambda}=\sup\{\pd{_\Lambda I} \ | \ I\in \inj{\Lambda}\}<\infty$ and $\silp{\Lambda}=\sup\{\id{_\Lambda P} \ | \ P\in \proj{\Lambda}\}<\infty$.
Our second main result is closely related to the property of the functors $\mathsf{Z}_B\colon \Mod{B} \lxr \Mod{\Lambda_{(0,0)}}$ and $\mathsf{Z}_A\colon \Mod{A} \lxr \Mod{\Lambda_{(0,0)}}$ being homological embeddings, i.e. $\Ext_{B}^n(Y,Y')\simeq \Ext_{\Lambda_{(0,0)}}(\mathsf{Z}_B(Y),\mathsf{Z}_B(Y'))$ for all $n\geq 0$ and $Y, Y'\in \Mod{B}$, and similarly for $\mathsf{Z}_A$. In this connection, we characterize for which Morita rings the above functors are homological embeddings, see Proposition~\ref{propstratifyingideals}.
Our main results in this part are summarized in the following theorem. For simplicity we state the result in part (ii) only for $\silp$.

\begin{thmA}
\label{thmAintroduction}
Let $\Lambda_{(0,0)}$ be a Morita ring which is an Artin algebra and has zero bimodule homomorphisms.
\begin{enumerate}
\item ({\bf \textsf{Gorenstein-projectives}}) Assume that the bimodules $_BM_A$ and $_AN_B$ satisfy the above compatibility conditions. Let $Z$ be a Gorenstein-projective $B$-module with a monomorphism $s\colon N\otimes_BZ\lxr X$, for some $A$-module $X$, such that $\Coker{s}$ lies in $\Gproj{A}$ and there is a monomorphism $t\colon M\otimes_{A}\Coker{s}\lxr Y$ with $\Coker{t}=Z$ and $Y$ an $B$-module. We set $\pi_X$, resp. $\pi_Y$, for the map $M\otimes_AX\lxr \Coker{s}$, resp. $N\otimes_BY\lxr \Coker{t}$. Then the tuple$\colon$
\[
\big(X,Y, (\iden_{M}\otimes \pi_X)\circ t, (\iden_N\otimes \pi_Y)\circ s \big) \ \in \ \Gproj\Lambda_{(0,0)}
\]

\item ({\bf \textsf{Gorenstein algebras}}) Assume that the following conditions hold$\colon$
\begin{enumerate}
\item $M_A$ is projective and $\pd{_BM}<\infty$.

\item $N_B$ is projective and $\pd{_AN}<\infty$.

\item The functors $\mz_{A}$ and $\mz_{B}$ are homological embeddings.
\end{enumerate}
If $\silp{A}<\infty$ and $\silp{B}<\infty$, then $\silp{\Lambda_{(0,0)}}<\infty$.
\end{enumerate}
\end{thmA}

As an application of Theorem~A (i), we construct examples of Gorenstein-projective modules over the double morphism category $\DMor(\smod\Lambda)$. 
We remark that we obtain Gorenstein-projective modules that lie in $\mono(\Lambda)$, see Corollary~\ref{corGorprojDelta}. Also, from Theorem~A (ii) we get examples of Morita rings which are Gorenstein algebras (Corollary~\ref{corGorensteinMorita}). This constitutes the first part of this work which provides the link with monomorphism categories via Morita rings that we study further in the second part. 

In Section~\ref{SectionGorsubcatcoherentfun} we study the subcategory $\C$ of $\mono(\Lambda)$, where $\Lambda$ is an Artin algebra, consisting of all monomorphisms $f\colon X\lxr Y$ such that the projective dimension of $X$ is finite. Our third main result is Theorem~B (i) where assuming that $\Lambda$ is Gorenstein we show that $\C$ is a Gorenstein subcategory of $\mono(\Lambda)$. We refer to Theorem~\ref{Gro} for its proof and to Definition~\ref{defnGorsubcat} for the precise notion of Gorenstein subcategories in the setting of exact categories. Moreover, inspired by recent work of Matsui and Takahashi \cite{MatsuiTakahashi} we continue our study on the exact subcategory $\C$ of $\mono(\Lambda)$ by considering the category of coherent functors $\smod{\underline{\C}}$ over the stable category $\underline{\C}$ of $\C$. Also, we define the subcategory $\Omega^{n}(\C)$ of $\C$ consisting of all nth syzygies of objects in $\C$ (subsection~\ref{subsectioncoherentgor}). In this context, the fourth main result of this paper is Theorem~B (ii), see Corollary~\ref{corcoherentfunctorgoren}, which shows that the category of coherent functors over $\underline{\C}$ is a Gorenstein abelian category in the sense of \cite{BR}. Finally, using a result of Beligiannis \cite{Bel:ABcontexts} we realize the singularity category \cite{MatsuiTakahashi} of $\smod{\underline{\C}}$ as the stable category of Cohen-Macaulay objects over $\smod{\underline{\C}}$.

\begin{thmB}
\label{thmBintr}
({\bf \textsf{Gorenstein categories}})
Let $\Lambda$ be an $n$-Gorenstein Artin algebra. 
\begin{enumerate}
\item $\C=\{(X, Y, f,0)\in \mono(\Lambda) \ | \ \pd{_{\Lambda}X}<\infty \}$ is an $n$-Gorenstein subcategory of $\mono(\Lambda)$.

\item For the category of coherent functors over $\underline{\C}$ and $\underline{\Omega^n(\C)}$ the following hold$\colon$
 \begin{enumerate}

\item $\smod\underline{\C}$ is a $3n$-Gorenstein abelian category.

\item $\smod{\underline{\Omega^n(\C)}}$ is a Frobenius abelian category.
\end{enumerate}
\noindent Moreover, there are the following triangle equivalences$\colon$
\[
\xymatrix{
\mD_{\mathsf{sg}}(\smod\underline{\C}) \ar[r]^{\simeq \ \ } & \underline{\Gproj}(\smod{\underline{\C}}) } \ \ \text{and} \ \ \xymatrix{
\mD_{\mathsf{sg}}(\smod{\underline{\Omega^n(\C)}}) \ar[r]^{ \ \ \simeq } & \underline{\smod} \, \underline{\Omega^n(\C)} }
\]
\end{enumerate}
\end{thmB}

Statement (ii) above is a consequence of Theorem~\ref{mainthmexactcat} which provides sufficient conditions on a subcategory $\B$ of an exact category $\A$ with enough projectives such that $\smod{\underline{\B}}$ is a Gorenstein abelian category. It should be noted that this result generalizes, and is inspired by, a result of Matsui and Takahashi \cite{MatsuiTakahashi}. 

\smallskip

\noindent{\bf Conventions and Notation.} We compose morphisms in a given category in a diagrammatic order. Our subcategories are assumed to be closed under isomorphisms and direct summands. For a ring $R$ we usually work with left $R$-modules and the corresponding category is denoted by $\Mod{R}$. By a module over an Artin algebra $\Lambda$, we mean a finitely generated left $\Lambda$-module and we denote by $\smod{\Lambda}$ the category of finitely generated left $\Lambda$-modules. For all unexplained notions and results concerning the representation theory of Artin algebras we refer to \cite{ARS}.

\section{Morita Rings and Monomorphism Categories}
\label{sectionMorringsMonocat}

In this section we fix notation and 
we collect several preliminary results on Morita rings and monomorphism categories that  will be used throughout the paper.

\subsection{Morita Rings}
\label{subsectionMoritarings}
Let $A$ and $B$ be two rings, $_AN_B$ an $A$-$B$-bimodule, $_BM_A$ a $B$-$A$-bimodule, and $\phi \colon   M\otimes_{A}N \lxr B$ a $B$-$B$-bimodule homomorphism, and $\psi \colon N\otimes_{B}M \lxr A$ an $A$-$A$-bimodule homomorphism. Then from the Morita context $\M = (A,N,M,B, \phi, \psi)$ we define the \textsf{Morita ring}$\colon$
\[
\Lambda_{(\phi,\psi)}=
         \begin{pmatrix}
           A & _AN_B \\
           _BM_A & B \\
         \end{pmatrix}
\]
where the addition of elements of $\Lambda_{(\phi,\psi)}$ is componentwise and multiplication is given by
\[
         \begin{pmatrix}
           a & n \\
           m & b \\
         \end{pmatrix}
       \cdot
         \begin{pmatrix}
           a' & n' \\
           m' & b' \\
         \end{pmatrix}=
         \begin{pmatrix}
           aa'+\psi(n\otimes m') & an'+nb' \\
           ma'+bm' & bb'+\phi(m\otimes n') \\
         \end{pmatrix}
\]
We assume that $\phi(m\otimes n)m'=m\psi(n\otimes m')$ and $n\phi(m\otimes n')=\psi(n\otimes m)n'$ for all $m,m'\in M$ and $n,n'\in N$. This condition ensures that $\Lambda_{(\phi,\psi)}$ is an associative ring. 

The description of the modules over a Morita ring $\Lambda_{(\phi,\psi)}$ is well known, see for instance \cite{Green}, but for completeness and due to our needs we also include it here. We introduce the following category.

Let $\M(\Lambda)$ be the category whose objects are tuples $(X,Y,f,g)$
where $X\in \Mod{A}$, $Y\in \Mod{B}$, $f\in\Hom_B(M\otimes_AX,Y)$ and $g\in \Hom_A(N\otimes_BY,X)$ such that the following diagrams are commutative$\colon$
\begin{equation}
\label{diagramsfortuples}
\xymatrix{
  N\otimes_B M\otimes_A X \ar[d]_{\psi\otimes \iden_X} \ar[r]^{ \ \ \ \ \ \iden_N\otimes f} &  N\otimes_BY \ar[d]^{g}     \\
  A\otimes_AX    \ar[r]^{\simeq} & X                  } \ \ \ \ \ \ \  \ \ \  \xymatrix{
  M\otimes_A N\otimes_B Y \ar[d]_{\phi\otimes \iden_Y} \ar[r]^{ \ \ \ \ \iden_M\otimes g} &  M\otimes_AX \ar[d]^{f}     \\
  B\otimes_BY    \ar[r]^{\simeq} & Y                  }
\end{equation}
We denote by $\Psi_X$ and $\Phi_Y$ the following compositions$\colon$
\[
\xymatrix@C=0.5cm{
  N\otimes_BM\otimes_AX \ar[rr]^{ \ \ \ \ \psi\otimes \iden_X}  \ar @/^1.5pc/[rrrr]^{{\Psi_X}} && A\otimes_AX  \ar[rr]^{\simeq} && X  } \ \ \ \ \xymatrix@C=0.5cm{
  M\otimes_AN\otimes_BY \ar[rr]^{ \ \ \ \ \phi\otimes \iden_Y}  \ar @/^1.5pc/[rrrr]^{{\Phi_Y}} && B\otimes_BY  \ar[rr]^{\simeq} && Y  }
\]
Let $(X,Y,f,g)$ and $(X',Y',f',g')$ be objects of $\M(\Lambda)$. Then a morphism $(X,Y,f,g)\lxr (X',Y',f',g')$ in $\M(\Lambda)$ is a pair of homomorphisms $(a,b)$, where $a\colon X\lxr X'$ is an $A$-morphism and $b\colon Y\lxr Y'$ is a $B$-morphism, such that the following diagrams are commutative$\colon$
\[
\xymatrix{
  M\otimes_A X \ar[d]_{\iden_{M}\otimes a} \ar[r]^{\ \ \ f} & Y \ar[d]^{b}     \\
  M\otimes_AX'    \ar[r]^{\ \ \ f'} & Y'                  } \ \ \ \ \ \ \  \ \ \  \xymatrix{
  N\otimes_B Y \ar[d]_{\iden_{N}\otimes b} \ar[r]^{\ \ \  g} &  X \ar[d]^{a}     \\
  N\otimes_BY'    \ar[r]^{\ \ \ g'} & X'                  }
\]
The relationship between $\Mod{\Lambda_{(\phi,\psi)}}$ and $\M(\Lambda)$ is given via the functor $\mf\colon \M(\Lambda)\lxr \Mod{\Lambda_{(\phi,\psi)}}$ which is defined on objects $(X,Y,f,g)$ of $\M(\Lambda)$ as follows$\colon$$\mf(X,Y,f,g)=X\oplus Y$ as abelian groups, with a $\Lambda_{(\phi,\psi)}$-module structure given by
\[
         \begin{pmatrix}
           a & n \\
           m & b \\
         \end{pmatrix}(x,y)=(ax+g(n\otimes y), by+f(m\otimes x))
\]
for all $a\in A, b\in B, n\in N, m\in M, x\in X$ and $y\in Y$. If $(a,b)\colon (X,Y,f,g)\lxr (X',Y',f',g')$ is a morphism in
$\M(\Lambda)$ then $\mf(a,b)=\bigl(\begin{smallmatrix}
a & 0 \\
0 & b
\end{smallmatrix}\bigr)\colon X\oplus Y\lxr X'\oplus Y'$. Then the functor $\mf$ turns out to be an equivalence of categories, see \cite[Theorem~$1.5$]{Green}. We refer to \cite{GP} for an extensive discussion of Morita rings mainly in the context of Artin algebras. We also refer to \cite[Chapter $3$]{Psaroud:PhD} for a thorough discussion on the abelian structure of Morita rings in a more general setting. For our purpose and for completeness, we summarize in the next remark several properties of Morita rings that we need throughout the paper.

\begin{rem}
\label{remMoritarings}
Let $\Lambda_{(\phi,\psi)}=\bigl(\begin{smallmatrix}
A & _AN_B \\
_BM_A & B
\end{smallmatrix}\bigr)$ be a Morita ring.
\begin{enumerate}

\item Throughout the paper we deal mainly with Morita rings which are Artin algebras. Then it is easy to observe that, see also \cite[Proposition $2.2$]{GP}, a Morita ring  $\Lambda_{(\phi,\psi)}$ is an Artin algebra if and only if there is a commutative artin ring $R$ such that $A$ and $B$ are Artin $R$-algebras and $M$ and $N$ are finitely generated over $R$ which acts centrally both on $M$ and $N$.

\item From now on we identify the modules over $\Lambda_{(\phi,\psi)}$ with the objects of $\M(\Lambda)$. 

\item A sequence of tuples $0 \lxr (X_1,Y_1,f_1,g_1) \lxr (X_2,Y_2,f_2,g_2) \lxr (X_3,Y_3,f_3,g_3) \lxr 0$ is exact in $\Mod{\Lambda_{(\phi,\psi)}}$ if and only if the sequences $0\lxr X_1\lxr X_2\lxr X_3\lxr 0$ and $0\lxr Y_1\lxr Y_2\lxr Y_3\lxr 0$ are exact in $\Mod{A}$ and $\Mod{B}$ respectively.

\item Let $(a,b)\colon (X,Y,f,g)\lxr (X',Y',f',g')$ be a morphism in $\Mod{\Lambda_{(\phi,\psi)}}$ and consider the maps $c\colon \Ker{a}\lxr X$ and $d\colon \Ker{b}\lxr Y$. Then the kernel of $(a,b)$ is the object $(\Ker{a}, \Ker{b}, h, j)$ where the maps $h$ and $j$ are induced from the following commutative diagrams$\colon$
\begin{equation}
\label{diagramskernels}
\ \  \ \ \ \ \ \ \ \ \ \ \ \xymatrix{
    M\otimes_A\Ker{a} \ar@{-->}[d]_{h} \ar[r]^{ \ \ \iden_{M}\otimes c} & M\otimes_{A}X \ar[d]_{f}
  \ar[r]^{ \iden_{M}\otimes a} & M\otimes_A X' \ar[d]^{f'} \\
   \Ker{b} \ \ar@{>->}[r]^{d} & Y \ar[r]^{b} & Y' } \ \ \ \ \  \xymatrix{
    N\otimes_B\Ker{b} \ar@{-->}[d]_{j} \ar[r]^{ \ \ \iden_N\otimes d} & N\otimes_BY \ar[d]_{g}
  \ar[r]^{\iden_{N}\otimes b} & N\otimes_BY' \ar[d]^{g'} \\
   \Ker{a} \ \ar@{>->}[r]^{c} & X \ar[r]^{a} & X'   }
\end{equation}
Similarly, we derive a description for the cokernel of the morphism $(a,b)$.
\end{enumerate}
\end{rem} 

As in \cite{GP} we define the following functors$\colon$

\begin{enumerate}
\item The functor $\mt_{A}\colon \Mod{A}\lxr \Mod{\Lambda_{(\phi,\psi)}}$ is defined by
$\mt_{A}(X)=(X,M\otimes_AX,\iden_{M\otimes X},\Psi_X)$ on the objects $X\in \Mod{A}$ and given an $A$-morphism $a\colon X\lxr X'$ then
$\mt_{A}(a)=(a,\iden_{M}\otimes a)$.

\item The functor $\mU_{A}\colon \Mod{\Lambda_{(\phi,\psi)}}\lxr \Mod{A}$ is defined by
$\mU_{A}(X,Y,f,g)=X$ on the objects $(X,Y,f,g)\in \Mod{\Lambda_{(\phi,\psi)}}$ and
given a morphism $(a,b)\colon (X,Y,f,g)\lxr (X',Y',f',g')$ in
$\Mod{\Lambda_{(\phi,\psi)}}$ then $\mU_{A}(a,b)=a$.

\item The functor $\mt_{B}\colon \Mod{B}\lxr \Mod{\Lambda_{(\phi,\psi)}}$ is defined by
$\mt_{B}(Y)=(N\otimes_BY,Y,\Phi_Y,\iden_{N\otimes Y})$ on the objects $Y\in \Mod{B}$ and given a $B$-morphism $b\colon Y\lxr Y'$ then $\mt_{B}(b)=(\iden_{N}\otimes b,b)$.

\item The functor $\mU_{B}\colon \Mod{\Lambda_{(\phi,\psi)}}\lxr \Mod{B}$ is defined by
$\mU_{B}(X,Y,f,g)=Y$ on the $\Lambda_{(\phi,\psi)}$-modules $(X,Y,f,g)$ and given a $\Lambda_{(\phi,\psi)}$-morphism $(a,b)\colon (X,Y,f,g)\lxr
(X',Y',f',g')$ then $\mU_{B}(a,b)=b$.

\item The functor $\mh_{A}\colon \Mod{A}\lxr \Mod{\Lambda_{(\phi,\psi)}}$ is defined by
$\mh_{A}(X)=(X,\Hom_A(N,X),\delta'_{M\otimes X}\circ \Hom_A(N,\Psi_X),\epsilon'_X)$ on the objects $X\in \Mod{A}$ and given an $A$-morphism $a\colon X\lxr X'$ then
$\mh_{A}(a)=(a,\Hom_A(N,a))$.

\item The functor $\mh_{B}\colon \Mod{B}\lxr \Mod{\Lambda_{(\phi,\psi)}}$ is defined by
$\mh_{B}(Y)=(\Hom_B(M,Y),Y,\epsilon_Y,\delta_{N\otimes Y}\circ \Hom_B(M,\Phi_Y))$ on the objects $Y\in \Mod{B}$ and given a $B$-morphism $b\colon Y\lxr Y'$ then
$\mh_{B}(b)=(\Hom_B(M,b),b)$.

\item Suppose that $\phi=0=\psi$. Then we define the functor $\mz_{A}\colon \Mod{A}\lxr \Mod{\Lambda_{(0,0)}}$ by $\mz_{A}(X)=(X,0,0,0)$ on the objects $X\in \Mod{A}$ and if
$a\colon X\lxr X'$ is an $A$-morphism then $\mz_{A}(a)=(a,0)$. Dually we define the functor $\mz_{B}\colon \Mod{B}\lxr \Mod{\Lambda_{(0,0)}}$.
\end{enumerate}

When a Morita ring is an Artin algebra we have the following description of the indecomposable projective and injective modules.

\begin{prop}\textnormal{\cite[Propositions $3.1$ and $3.2$]{GP}}
\label{prop:projmod}
Let $\Lambda_{(\phi,\psi)}$ be a Morita ring which is an Artin algebra. Then the following hold.
\begin{enumerate}
\item The indecomposable projective $\Lambda_{(\phi,\psi)}$-modules
are objects of the form$\colon$
\[
\left\{
  \begin{array}{lll}
   \mt_A(P)=(P, {M\otimes_{A}P}, \iden_{M\otimes_{A}P}, \Psi_P)  & \hbox{} \\
           & \hbox{} \\
   \mt_B(Q)=({N\otimes_BQ}, Q, \Phi_{Q}, \iden_{N\otimes_{B}Q})  & \hbox{}
  \end{array}
\right.
\]
where $P$ is an indecomposable projective $A$-module and
$Q$ is an indecomposable projective $B$-module.

\item The indecomposable injective $\Lambda_{(\phi,\psi)}$-modules
are objects of the form$\colon$
\[
\left\{
  \begin{array}{ll}
   \mh_A(I)=(I,\Hom_{A}(N,I),\delta'_{M\otimes I}\circ \Hom_A(N,\Psi_I),\epsilon'_I)  & \hbox{} \\
     & \hbox{} \\
   \mh_B(J)=(\Hom_{B}(M,J),J,\epsilon_{J},\delta_{N\otimes J}\circ \Hom_B(M,\Phi_J)) & \hbox{}
  \end{array}
\right.
\]
where $I$ is an indecomposable injective $A$-module and $J$ is an indecomposable injective $B$-module.
\end{enumerate}
\end{prop}
 
We continue now with examples of Morita rings which will be used in the sequel.

\begin{exam}
\label{examplesMoritarings}
\begin{enumerate}
\item Let $R$ be a ring with an idempotent element $e$. Then, from the Pierce decomposition of $R$ with respect to the idempotents $e$ and $f=1_{R}-e$, it follows that $R$ is the Morita ring with $A=eRe$, $B=fRf$, $N=eRf$, $M=fRe$ and the bimodule homomorphisms $\phi$, $\psi$ are induced by the multiplication in $R$.

\item Any pair $(A,P_{A})$, where $A$ is a ring and $P_{A}$ is a right $A$-module, induces a Morita ring as follows$\colon$
\[
\Lambda_{(\phi,\psi)}=\begin{pmatrix}
A & {\Hom_{A}(P,A)} \\
P & \End_{A}(P)
\end{pmatrix}
\]
with bimodule homomorphisms $\phi\colon P\otimes_A\Hom_{A}(P,A)\lxr \End_{A}(P)$, $p\otimes f\mapsto \phi(p\otimes f)(p')=pf(p')$ and $\psi\colon \Hom_{A}(P,A)\otimes_{\End_{A}(P)}P\lxr A$, $f\otimes p\mapsto \psi(f\otimes p)=f(p)$. It is well known that if the $A$-module $P$ is a progenerator, then the rings $A$ and $\End_{A}(P)$ are Morita equivalent.

\item Let $\Lambda$ be an Artin algebra and $X, Y$ two $\Lambda$-modules. Then the endomorphism algebra $\End_{\Lambda}(X\oplus Y)$ is a Morita ring$\colon$
\[
{\End}_{\Lambda}(X\oplus Y)\iso
         \begin{pmatrix}
           \End_{\Lambda}(X) & \Hom_{\Lambda}(X,Y) \\
           \Hom_{\Lambda}(Y,X) & \End_{\Lambda}(Y) \\
         \end{pmatrix}_{(\phi,\psi)}
\]
The bimodule homomorphisms $\phi\colon \Hom_{\Lambda}(Y,X)\otimes \Hom_{\Lambda}(X,Y)\lxr \End_{\Lambda}(Y) $ and $\psi\colon \Hom_{\Lambda}(X,Y)\otimes \Hom_{\Lambda}(Y,X)\lxr \End_{\Lambda}(X)$ are given by composition.

\item Let $\Lambda_{(0,0)}=\bigl(\begin{smallmatrix}
A & _AN_B \\
_BM_A & B
\end{smallmatrix}\bigr)$ be a Morita ring where the bimodule morphisms $\phi$ and $\psi$ are zero. Then we have an isomorphism of rings$\colon$
\[
\xymatrix@C=0.5cm{
  \Lambda_{(0,0)} \ar[r]^{\iso \ \ \ \ \ \ \ \ \ \ } & (A\times B)\ltimes M\oplus N  }
\]
where $(A\times B)\ltimes M\oplus N$ is the trivial
extension ring of $A\times B$ by the $(A\times B)$-$(A\times
B)$-bimodule $M\oplus N$. For the notion of trivial extension of rings and for the above isomorphism we refer to \cite{FGR}, see also \cite[Proposition $2.5$]{GP}. 

\item  Suppose that we have the following Morita ring$\colon$
\[
\Lambda_{(\phi,\psi)} = \begin{pmatrix} A & A \\ A & A \end{pmatrix}
\]
where every entry is a ring $A$. Then, it follows from the associativity of the multiplication that $\phi=\psi$, see \cite[Corollary $2.13$]{GP} for more details. A special case is when $\phi=0$, that is $\Delta_{(0,0)}:=\bigl(\begin{smallmatrix}
\Lambda & \Lambda \\
\Lambda & \Lambda
\end{smallmatrix}\bigr)\iso (\Lambda\times \Lambda)\ltimes \Lambda\oplus \Lambda$. In the next subsection we analyze the module category of $\Delta_{(0,0)}$ via recollements of abelian categories. 
\end{enumerate}
\end{exam}

We close this subsection with the next result, which shows that always a Morita ring gives rise to a recollement situation. This provides a way to relate the module category of a Morita ring with the module categories of its underlying rings. For the proof see \cite[Proposition $2.4$]{GP} and for more details on recollements of abelian categories we refer to \cite{Pira, Psaroud:homolrecol}.

\begin{prop}
\label{prop:propertiesrecollement}
Let $\Lambda_{(\phi,\psi)}$ be a Morita ring. Then the following diagrams$\colon$
\[
\xymatrix@C=0.5cm{
\Mod{{B/\Image{\phi}}} \ar[rrr]^{\inc} &&& \Mod{\Lambda_{(\phi,\psi)}} \ar[rrr]^{\mU_A } \ar @/_1.5pc/[lll]_{\mathsf{Q}_B}  \ar
 @/^1.5pc/[lll]_{\mathsf{P}_B} &&& \Mod{A}
\ar @/_1.5pc/[lll]_{\mt_A} \ar
 @/^1.5pc/[lll]_{\mh_A}
 }
\]
and
\[
\xymatrix@C=0.5cm{
\Mod{{A/\Image{\psi}}} \ar[rrr]^{\inc} &&& \Mod{\Lambda_{(\phi,\psi)}} \ar[rrr]^{\mU_B } \ar @/_1.5pc/[lll]_{\mathsf{Q}_A}  \ar
 @/^1.5pc/[lll]_{\mathsf{P}_A} &&& \Mod{B}
\ar @/_1.5pc/[lll]_{\mt_B} \ar
 @/^1.5pc/[lll]_{\mh_B}
 }
\]

\smallskip

\noindent are recollements of abelian categories, that is$\colon$

\smallskip

\begin{minipage}{0.49\textwidth}
\begin{enumerate}
\item $(\mt_A, \mU_A,\mh_A)$ is an adjoint triple.

\item The functors $\mt_A$ and $\mh_A$ are fully faithful.

\item $\Ker{\mU_A}=\Mod{{B/\Image{\phi}}}$.
\end{enumerate}
\end{minipage}
\begin{minipage}{0.49\textwidth}
\begin{enumerate}
\item $(\mt_B, \mU_B,\mh_B)$ is an adjoint triple.

\item The functors $\mt_B$ and $\mh_B$ are fully faithful.

\item $\Ker{\mU_B}=\Mod{{A/\Image{\psi}}}$.
\end{enumerate}
\end{minipage}

\smallskip

\noindent In particular, if $\phi=0=\psi$ then we have the following recollements of module categories$\colon$
\[
\xymatrix@C=0.4cm{
\Mod{B} \ar[rrr]^{\mz_B \ \ } &&& \Mod{\Lambda_{(0,0)}} \ar[rrr]^{\mU_A} \ar @/_1.5pc/[lll]_{\mathsf{Q}_B \ }  \ar
 @/^1.5pc/[lll]_{\mathsf{P}_B \ } &&& \Mod{A}
\ar @/_1.5pc/[lll]_{\mt_A} \ar
 @/^1.5pc/[lll]_{\mh_A}
 } \ \ \ \ \xymatrix@C=0.4cm{
\Mod{A} \ar[rrr]^{\mz_A \ \ } &&& \Mod{\Lambda_{(0,0)}} \ar[rrr]^{\mU_B} \ar @/_1.5pc/[lll]_{\mathsf{Q}_{A} \ }  \ar
 @/^1.5pc/[lll]_{\mathsf{P}_A \ } &&& \Mod{B}
\ar @/_1.5pc/[lll]_{\mt_B} \ar
 @/^1.5pc/[lll]_{\mh_B}
 }
\]
\end{prop}

Our aim next is to analyze the recollement of the module category of the Morita ring $\Delta_{(0,0)}=\bigl(\begin{smallmatrix}
R & R \\
R & R
\end{smallmatrix}\bigr)$, where $R$ is a unital associative ring. For this reason, we introduce in the next subsection the double morphism category of an abelian category. This construction can be considered as an abstract model for the category of modules over $\Delta_{(0,0)}$.

\subsection{The Double Morphism Category}
\label{subsectiondoulemorphism}
Let $\A$ be an abelian category. The {\bf \textsf{double morphism category}} of $\A$, denote by $\DMor(\A)$, has as objects diagrams of the form$\colon$
\[
\xymatrix@C=0.5cm{
X \ \ar@<-.7ex>[rr]_{f} && \ Y \ \ar@<-.7ex>[ll]_-{g } }
\]
where $f$ and $g$ are morphisms in $\A$ such that $f\circ g=0$ and $g\circ f=0$. We simply denote the objects as tuples $(X,Y,f,g)$. A morphism $(X,Y,f,g)\lxr (X',Y',f',g')$ in $\DMor(\A)$ is a pair $(a,b)$ of morphisms in $\A$, where $a\colon X\lxr X'$ and $b\colon Y\lxr Y'$, such that the following diagram commutes$\colon$
\[
\xymatrix@C=0.5cm{
X \ \ar@<-.7ex>[rr]_{f} \ar[d]^{a} && \ Y \ \ar@<-.7ex>[ll]_-{g } \ar[d]^{b}  \\
X' \ \ar@<-.7ex>[rr]_{f'} && \ Y' \ \ar@<-.7ex>[ll]_-{g' }}
\]
that is, $g\circ a=b\circ g'$ and $f\circ b=a\circ f'$. We show that the double morphism category $\DMor(\A)$ is an abelian category and that there is a recollement which relates $\DMor(\A)$ and $\A$. In order to give an abelian structure on $\DMor(\A)$, we provide another description of $\DMor(\A)$. In particular, we show that there is an equivalence of categories between $\DMor(\A)$ and $(\A\times \A)\ltimes H$, where $(\A\times \A)\ltimes H$ is the trivial extension of $\A\times \A$ by an endofunctor $H$, see Fossum-Griffith-Reiten \cite{FGR}.

We define the functor $H\colon \A\times \A \lxr \A\times \A$, $H(X,Y)=(Y,X)$, and given a morphism $(a,b)\colon (X,Y)\lxr (X',Y')$ then
$H(a,b)=(b,a)$. Then we can define the trivial extension
$(\A\times \A)\ltimes H$, where the objects are morphisms
$\alpha\colon H(X,Y)\lxr (X,Y)$ such that $H(\alpha)\circ \alpha=0$, and if
$\alpha\colon H(X,Y)\lxr (X,Y)$ and $\beta\colon H(X',Y')\lxr
(X',Y')$ are two objects in $(\A\times \A)\ltimes H$, then a
morphism between the objects $\alpha$ and $\beta$ is a morphism
$\gamma\colon (X,Y)\lxr (X',Y')$ such that the diagram
\[
\xymatrix{
  H(X,Y) \ar[d]_{\alpha} \ar[r]^{H(\gamma) \ \ } & H(X',Y') \ar[d]^{\beta} \\
  (X,Y) \ar[r]^{\gamma \ \ } & (X',Y')   }
\]
is commutative, where $\alpha=(a_1,a_2)$, $\beta=(b_1,b_2)$ and $\gamma=(c_1,c_2)$. Since the endofunctor $H$ is (right) exact, it follows from \cite{FGR} that the trivial extension $(\A\times \B)\ltimes H$ is an abelian category.

\begin{prop}
\label{propdoublemorrecol}
Let $\A$ be an abelian category.
\begin{enumerate}
\item There is an equivalence of categories$\colon$
\[
\xymatrix{
  \DMor(\A) \ \ar[r]^{\simeq \ \ \ \ } & \ (\A\times \A)\ltimes H }
\]
In particular, the double morphism category $\DMor(\A)$ is abelian.

\item There is a recollement of abelian categories$\colon$
\begin{equation}
\label{recdoublemorphismcat}
\xymatrix@C=0.5cm{
\A \ar[rrr]^{\inc} &&& \DMor(\A) \ar[rrr]^{{\mU_{\A}}} \ar @/_1.5pc/[lll]_{\mathsf{Q}_{\A}}  \ar @/^1.5pc/[lll]_{\mathsf{P}_{\A}} &&& \A
\ar @/_1.5pc/[lll]_{{\mt_{\A}}} \ar
 @/^1.5pc/[lll]_{\mh_{\A}}
 }
\end{equation}
\end{enumerate}
\begin{proof}
(i) Let $(X,Y,f,g)$ be an object of $\DMor(\A)$. We define the
functor
\[
\xymatrix{
  \F\colon \DMor(\A) \ \ar[r]^{ } & \ (\A\times \A)\ltimes H}, \ \F(X,Y,f,g)=\xymatrix{
    H(X,Y) \ar[r]^{ (g,f)} & (X,Y) }
\]
and given a morphism $(a,b)\colon (X,Y,f,g)\lxr
(X',Y',f',g')$ in $\DMor(\A)$ then $\F(a,b)=H(a,b)$. The functor $\F$
is well defined since the following composition
\[
\xymatrix{
  H^2(X,Y) \ar[r]^{(f,g) } & H(X,Y) \ar[r]^{(g,f) } & (X,Y)  }
\]
is zero, i.e. the object $\F(X,Y,f,g)$ lies in $(\A\times
\A)\ltimes H$. It is clear that the functor $\F$ is faithful. Let
\[
\xymatrix{
  \big[(Y,X) \ar[r]^{ (g,f) \ } & (X,Y)\big] \ar[r]^{(a,b) \ \  } & \big[(Y',X') \ar[r]^{ (g',f') \ \ } & (X',Y')\big] }
\]
be a morphism in $(\A\times \A)\ltimes H$. Then the following commutative diagram
\[
\xymatrix{
  (Y,X) \ar[d]_{(g,f)} \ar[rr]^{(b,a) \ } && (Y',X') \ar[d]^{(g',f')} \\
  (X,Y) \ar[rr]^{(a,b)} && (X',Y')   }
\]
implies that $(a,b)\colon (X,Y,f,g)\lxr (X',Y',f',g')$ is a morphism in $\DMor(\A)$ and $\F(a,b)=H(a,b)$. Thus the functor $\F$ is full. Finally, if $(a_1,a_2)\colon H(X,Y)\lxr (X,Y)$ is an object of $(\A\times \A)\ltimes H$, then since $H(a_1,a_2)\circ
(a_1,a_2)=0$ we infer that $(X,Y,a_2,a_1)\in \DMor(\A)$ such that $\F(X,Y,a_2,a_1)=(a_1,a_2)$. This shows that the functor
$\F$ is essentially surjective. Hence, the categories $\DMor(\A)$ and $(\A\times
\A)\ltimes H$ are equivalent and therefore the double morphism category $\DMor(\A)$ is abelian.

(ii) The functors appearing in diagram $(\ref{recdoublemorphismcat})$ were defined in subsection~\ref{subsectionMoritarings} for the module category of a Morita ring. In this case, if $A$ is an object in $\A$ then $\mt_{\A}(A)=(A,A,\iden_{A},0)$, $\mh_{\A}(A)=(A,A,0,\iden_{A})$ and for a tuple $(X,Y,f,g)$ in $\DMor(\A)$ we have $\mU_{\A}(X,Y,f,g)=X$. Similarly with subsection~\ref{subsectionMoritarings}, we get a description of these functors on morphisms. Then, it is easy to check that $(\mt_{\A},\mU_{\A},\mh_{\A})$ is an adjoint triple with $\mt_{\A}$ (equivalently, $\mh_{\A}$) fully faithful and the kernel of $\mU_{\A}$ is equivalent with $\A$, see also Proposition~\ref{prop:propertiesrecollement}. We infer that $(\A,\DMor(\A),\A)$ is a recollement of abelian categories.
\end{proof}
\end{prop}

In the following remark we collect some interesting properties of the recollement diagram of a double morphism category.

\begin{rem}
\label{remdoublemorphismcat}
Let $\DMor(\A)$ be the double morphism category of an abelian category $\A$ and consider the recollement situation $(\ref{recdoublemorphismcat})$.
\begin{enumerate}
\item The functors $\mt_{\A}\colon \A\lxr \DMor(\A)$ and $\mh_{\A}\colon \A\lxr \DMor(\A)$ are exact. Thus, the recollement $(\ref{recdoublemorphismcat})$ of $\DMor(\A)$ has the property that the left and right adjoint of the quotient functor $\mU_{\A}\colon \DMor(\A)\lxr \A$ are exact. In general, this property doesn't hold in a recollement situation.

\item Let $(X,Y,f,g)$ be an object in $\DMor(\A)$. Then the tuple $(Y,X,g,f)$ is also an object in $\DMor(\A)$ since the composition of morphisms is still zero. This gives a functor $\F\colon \DMor(\A)\lxr \DMor(\A)$, $\F(X,Y,f,g)=(Y,X,g,f)$, which turns out to be an auto-equivalence.

\item For an object $(X,Y,f,g)$ in $\DMor(\A)$ we define the exact functor $\mU_{\A}'\colon \DMor(\A)\lxr \A$ given by $\mU'_{\A}(X,Y,f,g)=Y$, and $\mU'_{\A}(a,b)=b$ for a morphism $(a,b)$ in $\DMor(\A)$. It is easy to check that $\mU'_{\A}$ is the middle functor of the adjoint triple $(\mh_{\A},\mU'_{\A},\mt_{\A})$ and therefore we obtain a recollement of abelian categories $(\A,\DMor(\A),\A)$. Then the following commutative diagram$\colon$
\[
\xymatrix{
    \DMor(\A) \ar[d]_{\F}^{\simeq} \ar[r]^{ \ \ \ \mU_{\A}} & \A \ar[d]^{\iden_{\A}} \\
   \DMor(\A) \ar[r]^{ \ \ \ \mU_{\A}'} & \A }
\]
shows that there is a natural equivalence of functors between $\mU'_{\A}\F$ and $\mU_{\A}$. Thus, from \cite[Definition $4.1$, Lemma $4.2$]{PsaroudVitoria:1} we infer that the two recollements of $\DMor(\A)$, i.e. the recollement $(\ref{recdoublemorphismcat})$ and the recollement $(\A,\DMor(\A),\A)$ given by $\mU'_{\A}$, are equivalent. From now on, we fix the recollement diagram $(\ref{recdoublemorphismcat})$ for the double morphism category $\DMor(\A)$.
\end{enumerate}
\end{rem}

We close this subsection with the next example, which was the starting point for introducing the double morphism category of an abelian category.

\begin{exam}
\label{examrecollementofDelta}
Let $R$ be a ring and consider the Morita ring $\Delta_{(0,0)}=\bigl(\begin{smallmatrix}
R & R \\
R & R
\end{smallmatrix}\bigr)$, see Example~\ref{examplesMoritarings} (v). From Remark~\ref{remMoritarings} (ii), the category $\Mod{\Delta_{(0,0)}}$ is equivalent to the double morphism category $\DMor{(\Mod{R})}$ of $\Mod{R}$. Then, from Proposition~\ref{propdoublemorrecol} (ii) we have the following recollement$\colon$
\begin{equation}
\label{recollementofMoritaringDelta}
\xymatrix@C=0.5cm{
\Mod{R} \ar[rrr]^{\mz_2} &&& \Mod{\Delta_{(0,0)}} \ar[rrr]^{\mU_1} \ar @/_1.5pc/[lll]_{\Cok}  \ar @/^1.5pc/[lll]_{\Ker} &&& \Mod{R}
\ar @/_1.5pc/[lll]_{\mt_1} \ar
 @/^1.5pc/[lll]_{\mh_1}
 }
\end{equation}
For later use, we fix the above notation for the functors of the recollement of $\Mod{\Delta_{(0,0)}}$. In particular, and relative to subsection~\ref{subsectionMoritarings} and Proposition~\ref{propdoublemorrecol}, we have$\colon$
\begin{enumerate}
\item The functor $\mt_1\colon \Mod{R}\lxr \Mod{\Delta_{(0,0)}}$ is given by $\mt_{1}(X)=(X, X, \iden_{X}, 0)$ on the objects $X\in \Mod{R}$ and for an $R$-morphism $a\colon X\lxr X'$ then $\mt_{1}(a)=(a, a)$. Moreover, the functor $\mt_1$ is exact. Similarly, the functor $\mt_2\colon \Mod{R}\lxr \Mod{\Delta_{(0,0)}}$ is given by $\mt_{2}(X)=(X, X, 0, \iden_{X})$ on the objects $X\in \Mod{R}$ and for an $R$-morphism $a\colon X\lxr X'$ then $\mt_{2}(a)=(a, a)$. Note that in this case $\mt_2$ is precisely the functor $\mh_1$ appearing in the recollement $(\ref{recollementofMoritaringDelta})$.

\item The functor $\mU_{\A}$ of Proposition~\ref{propdoublemorrecol} is now denoted by  $\mU_1\colon \Mod{\Delta_{(0,0)}}\lxr \Mod{R}$.

\item The functor $\mz_{2}\colon \Mod{R}\lxr \Mod{\Delta_{(0,0)}}$, given by $\mz_{2}(X)=(0,X,0,0)$ for $X\in \Mod{R}$, is the functor $\mz_B$ defined in subsection~\ref{subsectionMoritarings}.

\item The cokernel functor $\Cok\colon \Mod{\Delta_{(0,0)}}\lxr \Mod{R}$ is given by $\Cok(X,Y,f,g)=\Coker{f}$ on the objects $(X,Y,f,g)\in \Mod{\Delta_{(0,0)}}$ and for a $\Delta_{(0,0)}$-morphism $(a, b)\colon (X, Y, f, g)\lxr (X', Y', f', g')$ we have $\Cok(a, b)=c$, where $c\colon \Coker{f}\lxr \Coker{f'}$ is the induced morphism such that $b\circ \pi'=\pi\circ c$, where $\pi\colon Y\lxr \Coker{f}$ and $\pi'\colon Y'\lxr \Coker{f'}$. This is the functor $\mathsf{Q}_{\A}$ in Proposition~\ref{propdoublemorrecol}. 

\item The kernel functor $\Ker\colon \Mod{\Delta_{(0,0)}}\lxr \Mod{R}$ is given by $\Ker(X,Y,f,g)=\Ker{g}$ on the objects $(X,Y,f,g)$ of $\Mod{\Delta_{(0,0)}}$ and for a $\Delta_{(0,0)}$-morphism $(a, b)\colon (X, Y, f, g)\lxr (X',Y',f',g')$ we have $\Ker(a, b)=c$, where $c$ is the restriction map of $b$ to $\Ker{g}$. This is the functor $\mathsf{P}_{\A}$ in Proposition~\ref{propdoublemorrecol}. 
\end{enumerate}
\end{exam}

\subsection{Monomorphism Categories}
\label{subMonocat}
Let $\A$ be an abelian category. We denote by ${\Mor}{\A}$ the category of morphisms of $\A$. Recall that the objects of ${\Mor}{\A}$ are triples $(X,Y,f)$, where $f\colon X\lxr Y$ is a morphism in $\A$, and given two objects $(X,Y,f)$ and $(X',Y',f')$ then a morphism is a pair $(a,b)$ of maps in $\A$ such that the following diagram commutes$\colon$
\[
\xymatrix{
    X \ar[d]_{a} \ar[r]^{f} & Y \ar[d]^{b} \\
   X' \ar[r]^{f'} & Y' }
\]
Since the morphism category ${\Mor}{\A}$ is a special case of a trivial extension of abelian categories, see \cite{FGR}, it follows that ${\Mor}{\A}$ is an abelian category. Our particular interest is the monomorphism category ${\Mon}{\A}$ of $\A$, which is the full subcategory of ${\Mor}{\A}$ consisting of monomorphisms in $\A$.
The monomorphism category ${\Mon}{\A}$ is an extension closed additive subcategory of ${\Mor}{\A}$, and this implies that ${\Mon}{\A}$ is an exact category in the sense of Quillen \cite{Quillen}. Due to our needs in
the sequel (see Section~\ref{SectionGorsubcatcoherentfun}) we recall the notion of exact categories. Our approach follows the appendix of Keller \cite[Appendix A]{Keller}, see also \cite{Buhler}.

Let $\A$ be an additive category. A pair of composable morphisms $\xymatrix{
 X \ar[r]|-{f} & Y \ar[r]|-{g} & Z }$
is called \textsf{exact}, if $f$ is the kernel of $g$ and $g$ is the cokernel of $f$. Two exact pairs $(f,g)$ and $(f',g')$ are isomorphic if there are isomorphisms $a\colon X\lxr X'$, $b\colon Y\lxr Y'$ and $c\colon Z\lxr Z'$ such that $f\circ b=a\circ f'$ and $g\circ c=b\circ g'$. Let $\E$ be a class of exact pairs which is closed under isomorphisms. A pair $(f,g)$ in $\E$ is called a \textsf{conflation}, while the map $f$ is called an \textsf{inflation} and the map $g$ is called a \textsf{deflation}. Then the class $\E$ is an \textsf{exact structure} of $\A$ and $(\A,\E)$ is called an \textsf{exact category}, if the following axioms hold$\colon$

\smallskip

\noindent $(\textsf{Ex} \, 0)\colon$ The identity morphism of the zero object is a deflation.

\noindent $(\textsf{Ex} \, 1)\colon$ The composition of two deflations is a deflation.

\noindent $(\textsf{Ex} \, 2)\colon$ If $g\colon Y\lxr Z$ is a deflation and $c\colon Z'\lxr Z$ is a morphism, then there exists a pullback diagram$\colon$
\[
\xymatrix{
  Y' \ar@{-->}[d]_{c'} \ar@{-->}[r]^{g'} &  Z' \ar[d]^{c}   \\
  Y \ar[r]^{g} & Z                  }
\]
such that $g'$ is a deflation.

\noindent $(\textsf{Ex} \, 2)^{\op}\colon$ If $f\colon X\lxr Y$ is an inflation and $a\colon X\lxr X'$ is a morphism, then there exists a pushout diagram$\colon$
\[
\xymatrix{
  X \ar[d]_{a} \ar[r]^{f} & Y \ar@{-->}[d]^{a'}   \\
  X' \ar@{-->}[r]^{f'} & Y'                  }
\]
such that $f'$ is an inflation.

Let $(\A,\E)$ be an exact category. Recall that an object $P$ in $\A$ is {\bf \textsf{projective}} if the functor $\Hom_{\A}(P,-)$ sends conflations to short exact sequences, and $\A$ has {\bf \textsf{enough projectives}} if for any object $A$ in $\A$ there exists a deflation $g\colon P\lxr A$ with $P\in \Proj{\A}$. Dually we have the notions of injective objects and enough injectives. 

Of particular interest in representation theory is the morphism category over the category of finitely generated modules over an Artin algebra $\Lambda$. In particular, there is an equivalence of categories between the morphism category $\Mor(\smod{\Lambda})$ and the module category $\smod{\mathsf{T}_2(\Lambda)}$, see \cite{ARS}. In this case, the monomorphism category ${\Mon}{(\smod{\Lambda})}$, simply denoted by ${\mon}{(\Lambda)}$, is the full subcategory of $\smod{\mathsf{T}_2(\Lambda)}$ consisting of all monomorphisms of $\Lambda$-modules. Note that this holds for any ring $R$, i.e. $\Mod{\mathsf{T}_2(R)}\simeq \Mor(\Mod{R})$.

We now return to the double morphism category. For an abelian category $\A$, define the monomorphism categories of $\DMor(\A)$
as follows$\colon$$\Mono_1(\A)=\{(X,Y,f,0) \ | \ f\colon X
\lxr Y \ \text{monomorphism in} \ \A \}$ and $\Mono_2(\Lambda)=\{(X,Y,0,g) \ | \ g\colon Y\lxr
X \ \text{monomorphism in} \ \A \}$. In the next result we show that the above monomorphism categories are exact and that they are equivalent.

\begin{lem}
\label{lemequivmono}
Let $\A$ be an abelian category. Then the monomorphism categories $\Mon(\A)$, $\Mono_1(\A)$ and $\Mono_2(\A)$ are equivalent as exact categories.
\begin{proof}
There is a functor $\F\colon \Mor(\A)\lxr \DMor(\A)$ defined by
$\F(X,Y,f)=(X,Y,f,0)$ on the objects $(X,Y,f)\in \Mor(\A)$ and given a morphism $(a,b)\colon (X,Y,f)\lxr (X',Y',f')$ in $\Mor(\A)$ then
$\F(a,b)=(a,b)$. Note that the functor $\F$ is an exact full embedding and using $\F$ we view $\Mon(\A)$ as a full subcategory of the double morphism category $\DMor(\A)$. From the Snake Lemma it follows easily that the monomorphism categories $\Mon(\A)$, $\Mono_1(\A)$ and $\Mono_2(\A)$ are extension closed additive subcategories of the abelian category $\DMor(\A)$. This implies that they are exact categories, where the conflations are short exact sequences in $\DMor(\A)$ with terms in the corresponding monomorphism categories. Then, it follows easily that the functor $\F$ provides an equivalence between $\Mon(\A)$ and $\Mono_1(\A)$. Also, if $(X,Y,f,0)$ lies in $\Mono_1(\A)$ then the object $(Y,X,0,f)$ lies in $\Mono_2(\A)$. Then, using this correspondence, we infer that the exact categories $\Mono_1(\A)$ and $\Mono_2(\A)$ are equivalent.
\end{proof}
\end{lem}

From now on the {\bf \textsf{monomorphism category}} of an abelian category $\A$, denoted by $\Mono(\A)$, is the category $\Mono_1(\A)$, i.e. the image of $\Mon(\A)$ via the functor $\F$. Note that when $\A$ is exact, the monomorphism category of $\A$ is the {\em inflation category} of $\A$, that is, the morphisms of $\A$ which are inflations. We continue to call this category the monomorphism category of $\A$ and we denote it by $\Mono(\A)$ as well. The next result provides a description of the projective and injective objects in $\Mono(\A)$ and will be used in Section~\ref{SectionGorsubcatcoherentfun}. Note that this is proved in \cite[Lemma $2.1$]{Chen} for the monomorphism subcategory $\Mon(\A)$ of $\Mor(\A)$, but for completeness we provide a proof in our setting. 
 
\begin{lem}
\label{lemprojinjMono}
Let $\A$ be an exact $($abelian$)$ category with enough projective and injective objects. Then the monomorphism category $\Mono(\A)$ has enough projective and injective objects, in particular$\colon$
\begin{enumerate}
\item $\Proj(\Mono(\A))=\add\{ \mt_1(P)\oplus \mz_2(Q) \ | \ P, Q\in \Proj{\A}\}$, and

\item $\Inj(\Mono(\A))=\add\{ \mt_1(I)\oplus \mz_2(J) \ | \ I, J\in \Inj{\A}\}$.

\end{enumerate}
\begin{proof}
We show (i), statement (ii) follows similarly. We first claim that if $P$ is a projective object of $\A$, then the objects $\mz_{2}(P)=(0,P,0,0)$ and $\mt_{1}(P)=(P,P,\iden_{P},0)$ are projectives in $\Mono(\A)$. Indeed, let $(X,Y,f,0)$ be an object in $\Mono(\A)$ and let $(f,g)$ be a conflation of $\A$. It is easy to check that we have the isomorphisms  $\Hom_{\Mono(\A)}(\mz_{2}(P),(X,Y,f,0)) \iso \Hom_{\A}(P,Y)$ and $\Hom_{\Mono(\A)}(\mt_{1}(P),(X,Y,f,0)) \iso \Hom_{\A}(P,X)$. Let $(X_1,Y_1,f_1,0) \lxr (X_2,Y_2,f_2,0) \lxr (X_3,Y_3,f_3,0)$ be a conflation in $\Mono(\A)$. Then, from \cite[Corollary $8.13$]{Buhler}, we have the following commutative diagram$\colon$
\[
\xymatrix{
 X_1 \ar[d]^{f_1} \ar[r] & X_2 \ar[d]^{f_2} \ar[r]^{} & X_3 \ar[d]^{f_3} \\
Y_1 \ar[d]^{p_1} \ar[r] & Y_2 \ar[d]^{p_2} \ar[r]^{} & Y_3 \ar[d]^{p_3} \\
Z_1 \ar[r] & Z_2  \ar[r]^{} & Z_3 }
\]
where all rows and columns are conflations in $\A$. From the above isomorphisms and since the functor $\Hom_{\A}(P,-)$ sends conflations to short exact sequences, we get that the objects $\mz_{2}(P)$ and $\mt_{1}(P)$ lie in $\Proj(\Mono(\A))$. Now we show that the exact category $\Mono(\A)$ has enough projective objects. Let $(X,Y,f,0)$ be an object in $\Mono(\A)$. Since $\A$ has enough projectives, there are deflations $a\colon P\lxr X$ and $b\colon Q\lxr Z$ with $P, Q\in \Proj{\A}$. Then, from \cite[Corollary $8.13$]{Buhler}, we have the commutative diagram$\colon$
\begin{equation}
\label{comdiagramenoughproj}
\xymatrix{ K' \ar[r]^{i} \ar[d]^{a'} & K \ar[r]^{j } \ar[d]^{c'} & K'' \ar[d]^{b'} \\
 P \ar[d]^{a} \ar[r]^{(\begin{smallmatrix}
1 & 0 \\
\end{smallmatrix}) \ \ \ \ } & P\oplus Q \ar[r]^{ \ \bigl(\begin{smallmatrix}
0  \\
1
\end{smallmatrix}\bigr) } \ar[d]^{c} & Q \ar[d]^{b}  \\
 X \ar[r]^{f} & Y \ar[r]^{p} & Z }
\end{equation}
where all rows and columns are conflations of $\A$. Note that $c=\bigl(\begin{smallmatrix}
a\circ f  \\
d
\end{smallmatrix}\bigr)$, where $d\colon Q\lxr Y$ such that $d\circ p=b$. Then we have the following conflation in $\Mono(\A)\colon$
\[
\xymatrix{
(K',K,i,0) \ar[r]^{(a',c') \ \ \ \ \ \ \  \ } & (P,P\oplus Q,(\iden_{P} \ 0),0) \ar[r]^{ \ \ \ \ \ \ (a,c) } & (X,Y,f,0) }
\]
where $(P,P\oplus Q,(\iden_{P} \ 0),0)=\mt_1(P)\oplus \mz_2(Q)\in \Proj{(\Mono(\A))}$ and $\Ker(a,c)=(\Ker{a},\Ker{c},i,0)\in \Mono(\A)$.
Hence, the exact category $\Mono(\A)$ has enough projective objects. \end{proof}
\end{lem}

Let $\Lambda$ be an Artin algebra and consider the Morita ring $\Delta_{(0,0)}$ which is an Artin algebra. In this case, the monomorphism category of $\Lambda$ is the following full subcategory of $\smod{\Delta_{(0,0)}}$, i.e. of $\DMor(\smod\Lambda)\colon$  
\begin{equation}
\label{eqmono}
\mono(\Lambda) = \{ (X,Y,f,0) \ | \ f\colon X\lxr Y \ \ \text{is a monomorphism}\}
\end{equation}

We close this section with the next result where we collect some useful properties of $\mono(\Lambda)$ that we need in the sequel.

\begin{lem}
\label{adjointtriples}
Let $\Lambda$ be an Artin algebra. Then the following statements hold.
\begin{enumerate}
\item The monomorphism category $\mono(\Lambda)$ is an exact category which is closed under kernels.

\item We have the adjoint triples $(\Cok,\mz_2,\mU_2)$ and $(\mU_2,\mt_1,\mU_1)\colon$
\[
\xymatrix@C=0.5cm{
\smod{\Lambda} \ar[rrr]^{\mz_2   } &&& {\mono(\Lambda)}
\ar @/_1.5pc/[lll]_{\Cok }  \ar
 @/^1.5pc/[lll]_{ \mU_2 } \ar @/^1.5pc/[rrr]^{\mU_2} \ar @/_1.5pc/[rrr]^{\mU_1 }   &&& \smod{\Lambda} \ar[lll]_{ \mt_1 }
 }
\]
The above functors are exact and preserve projective objects, and $\mz_2$ and $\mt_1$ are fully faithful.
\end{enumerate}
\begin{proof}
(i) From Lemma~\ref{lemequivmono} the monomorphism category $\mono(\Lambda)$ is exact, since it is an extension closed subcategory of $\smod{\Delta_{(0,0)}}$. Let $(a,b)\colon (X,Y,f,0)\lxr (X',Y',f',0)$ be a morphism in $\smod{\Delta_{(0,0)}}$ with $(X,Y,f,0)$ and $(X',Y',f',0)$ in $\mono(\Lambda)$. Consider the following exact commutative diagram$\colon$
\[
\xymatrix{
0\ar[r] & \Ker{a} \ar[r]^{i } \ar[d]^{h} & X \ar[r]^{a} \ar[d]^{f} & X' \ar[d]^{f'}    \\
0 \ar[r] & \Ker{b} \ar[r]^{j} & Y \ar[r]^{b} & Y'  }
\]
Since the composition $f\circ i$ is a monomorphism, it follows that the map $h$ is a monomorphism. Then $\Ker(a,b)=(\Ker{a},\Ker{b},h,0)$ lies in $\mono(\Lambda)$. We infer that $\mono(\Lambda)$ is closed under kernels.

(ii) It is easy to check that the above functors form adjoint pairs, see Proposition~\ref{prop:propertiesrecollement} and Example~\ref{examrecollementofDelta}. Since $\Cok$, $\mz_2$, $\mU_2$ and $\mt_1$ are left adjoint functors of exact functors it follows that they preserve projective objects. The functor $\mU_1$ preserves projectives by the description of $\proj(\mono(\Lambda))$ given in Lemma~\ref{lemprojinjMono}, see also Example~\ref{examrecollementofDelta}. Moreover, it follows easily from the definition that the functors $\mz_2$, $\mU_2$, $\mt_1$ and $\mU_1$ are exact, and moreover that $\mz_2$ and $\mt_1$ are fully faithful, see again Example~\ref{examrecollementofDelta}. It remains to show that the cokernel functor $\Cok\colon \mono(\Lambda)\lxr \smod\Lambda$ is exact. Let $(X_1,Y_1,f_1,0)\lxr (X_2,Y_2,f_2,0)\lxr (X_3,Y_3,f_3,0)$ be a conflation in $\mono(\Lambda)$. Then we have the following exact commutative diagram$\colon$
\[
\xymatrix{
0 \ar[r] & X_1 \ar[r]^{} \ar[d]^{f_1} & X_2 \ar[d]^{f_2} \ar[r]^{} & X_3 \ar[r] \ar[d]^{f_3} & 0   \\
0 \ar[r] & Y_1 \ar[r] & Y_2 \ar[r] & Y_3 \ar[r]^{} & 0 }
\]
where the maps $f_1, f_2$ and $f_3$ are monomorphisms. From the Snake Lemma in the above diagram, it follows that the sequence $0\lxr \Coker{f_1}\lxr \Coker{f_2}\lxr \Coker{f_3}\lxr 0$ is exact and therefore the functor $\Cok_1\colon \mono(\Lambda)\lxr \smod{\Lambda}$ is exact.
\end{proof}
\end{lem}

\section{Gorenstein-Projective Modules over Morita Rings}
\label{GPmodulesoverMorita}
Our aim in this section is to study Gorenstein-projective modules over Morita rings. In particular we provide a method for constructing Gorenstein-projective modules over Morita rings, which are Artin algebras and satisfy certain conditions, from Gorenstein-projective modules of the underlying algebras. This section is divided into three subsections and the main result is stated in the second subsection. We start by recalling the notion of Gorenstein-projective modules and we also fix notation.

Let $\Lambda$ be an Artin algebra. An acyclic complex of projective $\Lambda$-modules $\textsf{P}^{\bullet}\colon \cdots \lxr P^{i-1}\lxr P^i\lxr P^{i+1}\lxr \cdots$ is called {\bf \textsf{totally acyclic}}, if the complex $\Hom_{\Lambda}(\textsf{P}^{\bullet},\Lambda)$ is acyclic. Then, a $\Lambda$-module $X$ is called {\bf \textsf{Gorenstein-projective}}, if it is of the form $X=\Coker{(P^{-1}\lxr P^0)}$ for some totally acyclic complex $\textsf{P}^{\bullet}$ of projective $\Lambda$-modules. We denote by $\Gproj{\Lambda}$ the full subcategory of $\smod{\Lambda}$ consisting of the finitely generated Gorenstein-projective $\Lambda$-modules. Moreover, we denote $\{X\in \smod{\Lambda} \ | \ \Ext_{\Lambda}^1(\Gproj{\Lambda},X)=0\}$ by $(\Gproj{\Lambda})^{\bot}$. Recall also from \cite{Bel:virtually, Bel:finiteCMtype}, that an Artin algebra $\Lambda$ is said to be of {\bf \textsf{finite Cohen-Macaulay type}}, if the category $\Gproj{\Lambda}$ is of finite representation type, i.e. the set of isomorphism classes of indecomposable finitely generated Gorenstein-projective modules is finite. Finally, for a $\Lambda$-module $X$ we denote by $\add{X}$ the full subcategory of $\smod{\Lambda}$ consisting of all direct summands of finite direct sums of $X$.

\subsection{Lifting Gorenstein-Projective Modules}\label{subLifting}
From Proposition~\ref{prop:propertiesrecollement} it follows that the functors $\mt_{A}\colon$ $\smod{A}\lxr \smod{\Lambda_{(\phi,\psi)}}$ and $\mt_{B}\colon \smod{B}\lxr \smod{\Lambda_{(\phi,\psi)}}$ preserve projective modules. In this subsection we investigate when the functors $\mt_{A}$ and $\mt_{B}$ preserve Gorenstein-projective modules. The first step towards this problem, is to examine when the above functors preserve totally acyclic complexes. Under some conditions, this is achieved in the next result.

\begin{prop}\label{proptotallyacyclic}
Let $\Lambda_{(\phi,\psi)}=\bigl(\begin{smallmatrix}
A & _AN_B \\
_BM_A & B
\end{smallmatrix}\bigr)$ be a Morita ring which is an Artin algebra.
\begin{enumerate}
\item Assume that the functor $M\otimes_{A}-\colon \smod{A}\lxr \smod{B}$ sends acyclic complexes of projective $A$-modules to acyclic complexes of $B$-modules and $\add{_AN}\subseteq (\Gproj{A})^{\bot}$. Then a complex $\mathsf{P}^{\bullet}$ in $\smod{A}$
is totally acyclic if and only if the complex $\mt_A(\mathsf{P}^{\bullet})$ is totally acyclic in $\smod{\Lambda_{(\phi,\psi)}}$.

\item Assume that the functor $N\otimes_{B}-\colon \smod{B}\lxr \smod{A}$ sends acyclic complexes of projective $B$-modules to acyclic complexes of $A$-modules and $\add{_BM}\subseteq (\Gproj{B})^{\bot}$. Then a complex $\mathsf{P}^{\bullet}$ in $\smod{B}$ is totally acyclic if and only if the complex $\mt_B(\mathsf{P}^{\bullet})$ is totally acyclic in $\smod{\Lambda_{(\phi,\psi)}}$.
\end{enumerate}
\begin{proof}
We prove only (i), the statement (ii) is dual. Assume that
\[
\mathsf{P}^{\bullet}\colon \
\xymatrix{
  \cdots \ar[r]_{}^{} & P^{-1} \ar[r]^{d^{-1}} & P^0 \ar[r]^{d^0} & P^{1} \ar[r] & \cdots }
\]
is a totally acyclic complex of projectives in $\smod{A}$. Then, by the assumption on the functor $M\otimes_{A}-$ and Remark~\ref{remMoritarings} (iv), we obtain that the following complex$\colon$
\[
\mt_{A}(\mathsf{P}^{\bullet})\colon \
\xymatrix{
  \cdots \ar[r]_{}^{} & \mt_{A}(P^{-1}) \ar[rr]^{\mt_{A}(d^{-1})} && \mt_{A}(P^0) \ar[rr]^{\mt_{A}(d^0)} && \mt_{A}(P^{1}) \ar[r] & \cdots }
\]
is exact, where each $\mt_A(P^i)$ lies in $\proj{\Lambda_{(\phi,\psi)}}$ by Proposition~\ref{prop:propertiesrecollement}. We show now that the complex $\Hom_{\Lambda_{(\phi,\psi)}}(\mt_{A}(\mathsf{P}^{\bullet}), (X,Y,f,g))$ is acyclic for all $(X,Y,f,g)$ in $\proj{\Lambda_{(\phi,\psi)}}$. In fact, from Proposition~\ref{prop:projmod} (i) it is enough to consider only the complexes $\Hom_{\Lambda_{(\phi,\psi)}}(\mt_{A}(\mathsf{P}^{\bullet}), \mt_{A}(P))$ and $\Hom_{\Lambda_{(\phi,\psi)}}(\mt_{A}(\mathsf{P}^{\bullet}), \mt_{B}(Q))$, where $P$ lies in $\proj{A}$ and $Q$ lies in $\proj{B}$. In the first case, the complex $\Hom_{\Lambda_{(\phi,\psi)}}(\mt_{A}(\mathsf{P}^{\bullet}), \mt_{A}(P))$ is acyclic since the complex $\Hom_{A}(\mathsf{P}^{\bullet},P)$ is acyclic and from Proposition~\ref{prop:propertiesrecollement} the functor $\mt_{A}$ is fully faithful. Let $Q$ be a projective $B$-module. Then, by using the adjoint pair $(\mt_{A},\mU_{A})$, we have the following commutative diagram$\colon$
\[
\xymatrix@C=0.4cm{
\cdots \ar[r] &  \Hom_{\Lambda_{(\phi,\psi)}}(\mt_A(P^1),\mt_{B}(Q)) \ar[d]_{\iso} \ar[r]^{} & \Hom_{\Lambda_{(\phi,\psi)}}(\mt_A(P^0),\mt_{B}(Q)) \ar[d]_{\iso}
  \ar[r]^{} & \Hom_{\Lambda_{(\phi,\psi)}}(\mt_A(P^{-1}),\mt_{B}(Q)) \ar[d]_{\iso} \ar[r]^{} & \cdots  \\
 \cdots \ar[r] & \Hom_{A}(P^1,N\otimes_BQ) \ar[r]^{} & \Hom_{A}(P^0,N\otimes_BQ) \ar[r]^{} & \Hom_{A}(P^{-1},N\otimes_BQ) \ar[r]^{} & \cdots   }
\]
Since $N\otimes_BQ$ is a direct sum of summands of $N$ and $\add{_AN}\subseteq (\Gproj{A})^{\bot}$, it follows that the complex $\Hom_{A}(\mathsf{P}^{\bullet},N\otimes_BQ)$ is acyclic and therefore the complex $\Hom_{\Lambda_{(\phi,\psi)}}(\mt_{A}(\mathsf{P}^{\bullet}), \mt_{B}(Q))$ is also acyclic. We infer that the complex $\mt_A(\mathsf{P}^{\bullet})$ is totally acyclic.

Conversely, assume that $\mathsf{P}^{\bullet}$ is a complex of $A$-modules such that $\mt_{A}(\mathsf{P}^{\bullet})$ is totally acyclic. If we apply the functor $\mU_{A}$ to the complex $\mt_{A}(\mathsf{P}^{\bullet})$, we get that the complex $\mathsf{P}^{\bullet}\colon \cdots \lxr P^{-1}\lxr P^0\lxr P^1\lxr \cdots $ is acyclic. Note that since the functor $\mt_{A}$ is right exact and fully faithful it follows that each $P^i$ lies in $\proj{A}$, see Proposition~\ref{prop:propertiesrecollement}. Then, for every projective $A$-module $P$, we derive as above that the complex $\Hom_{A}(\mathsf{P}^{\bullet},P)$ is acyclic.
We remark that in this direction we did not make use of our assumptions.
\end{proof}
\end{prop}

We refer to the above conditions as the {\bf \textsf{compatibility conditions}} on the bimodules $_AN_B$ and $_BM_A$.

\begin{exam}
\label{examconditionsformainthm}
Let $\Lambda_{(\phi,\psi)}=\bigl(\begin{smallmatrix}
A & _AN_B \\
_BM_A & B
\end{smallmatrix}\bigr)$ be a Morita ring which is an Artin algebra. 
\begin{enumerate}
\item Assume that $M_A$ is projective as a right $A$-module and $_AN$ is projective as a left $A$-module. Then the functor $M\otimes_{A}-\colon \smod{A}\lxr \smod{B}$ is exact and the subcategory $\add{_AN}$ lies in $(\Gproj{A})^{\bot}$. Hence, from Proposition~\ref{proptotallyacyclic} (i) it follows that a complex $\mathsf{P}^{\bullet}$ is totally acyclic in $\smod{A}$ if and only if the complex $\mt_A(\mathsf{P}^{\bullet})$ is totally acyclic in $\smod{\Lambda_{(\phi,\psi)}}$. Similarly, if $N_B$ is projective as a right $B$-module and $_BM$ is projective as a left $B$-module, then the statement of Proposition~\ref{proptotallyacyclic} (ii) holds. In particular, consider the case of the Morita ring 
 $\Delta_{(\phi,\phi)}=\bigl(\begin{smallmatrix}
\Lambda & \Lambda \\
\Lambda & \Lambda
\end{smallmatrix}\bigr)$, see Example~\ref{examrecollementofDelta}. Then it follows that a complex $\mathsf{P}^{\bullet}$ in $\smod{\Lambda}$ is totally acyclic if and only if the complex $\mt_1(\mathsf{P}^{\bullet})$ is totally acyclic in $\smod{\Delta_{(\phi,\phi)}}$ if and only if the complex $\mt_2(\mathsf{P}^{\bullet})$ is totally acyclic in $\smod{\Delta_{(\phi,\phi)}}$.

\item Assume that $\pd{M_A}<\infty$ and $\pd{_AN}<\infty$ (or $\id{_AN}<\infty$). Then from \cite[Proposition $1.3$]{Zhang}, it follows that the functor $M\otimes_{A}-\colon \smod{A}\lxr \smod{B}$ sends acyclic complexes of projective $A$-modules to acyclic complexes of $B$-modules and $\add{_AN}\subseteq (\Gproj{A})^{\bot}$. Hence, if $\Lambda_{(\phi,\psi)}$ is a Morita ring which is an Artin algebra such that $\pd{M_A}<\infty$ and $\pd{_AN}<\infty$ (or $\id{_AN}<\infty$), then from Proposition~\ref{proptotallyacyclic} (i) we get that the functor $\mt_{A}$ preserves totally acyclic complexes. Dually, if we assume that $\pd{N_B}<\infty$ and $\pd{_BM}<\infty$ (or $\id{_BM}<\infty$), then the conditions of Proposition~\ref{proptotallyacyclic} (ii) are satisfied and therefore the functor $\mt_B$ preserves totally acyclic complexes. Note that this example generalizes the situation mentioned in (i).
\end{enumerate}
\end{exam}

As a consequence of Proposition~\ref{proptotallyacyclic} we have the following result, which provides sufficient conditions such that the functors $\mt_{A}$ and $\mt_{B}$ lift Gorenstein-projective modules. In particular, we derive that Cohen-Macaulay finiteness of the Morita ring is inherited to the underlying algebras as well.

\begin{cor}\label{corliftinggorproj}
Let $\Lambda_{(\phi,\psi)}$ be a Morita ring which is an Artin algebra.
\begin{enumerate}
\item Assume that the functor $M\otimes_{A}-\colon \smod{A}\lxr \smod{B}$ sends acyclic complexes of projective $A$-modules to acyclic complexes of $B$-modules and $\add{_AN}\subseteq (\Gproj{A})^{\bot}$.
\begin{enumerate}
\item If $X\in \Gproj{A}$ then $\mt_{A}(X)\in \Gproj{\Lambda_{(\phi,\psi)}}$.

\item If $\Lambda_{(\phi,\psi)}$ is of finite Cohen-Macaulay type, then $A$ is also of finite Cohen-Macaulay type.
\end{enumerate}

\item Assume that the functor $N\otimes_{B}-\colon \smod{B}\lxr \smod{A}$ sends acyclic complexes of projective $B$-modules to acyclic complexes of $A$-modules and $\add{_BM}\subseteq (\Gproj{B})^{\bot}$.
\begin{enumerate}
\item If $Y\in \Gproj{B}$ then $\mt_{B}(Y)\in \Gproj{\Lambda_{(\phi,\psi)}}$.

\item If $\Lambda_{(\phi,\psi)}$ is of finite Cohen-Macaulay type, then $B$ is also of finite Cohen-Macaulay type.
\end{enumerate}
\end{enumerate}
\end{cor}

Now we turn our attention to the algebra $\Delta_{(\phi,\phi)}=\bigl(\begin{smallmatrix}
\Lambda & \Lambda \\
\Lambda & \Lambda
\end{smallmatrix}\bigr)$. We recall the following.

\begin{prop}
\label{proppropertiesDelta}
Let $\Lambda$ be an Artin algebra and let $n\geq 0$ be a natural number.
\begin{enumerate}
\item \textnormal{\cite[Corollary $6.4$]{GP}} $\Lambda$ is $n$-Gorenstein if and only if the Morita ring $\Delta_{(\phi,\phi)}$ is $n$-Gorenstein algebra.

\item \textnormal{\cite[Corollary $6.6$]{GP}} Assume that $\Lambda$ is Gorenstein. Then a $\Delta_{(\phi,\phi)}$-module $(X,Y,f,g)$ is Gorenstein-projective if and only if $X$ and $Y$ are Gorenstein-projective $\Lambda$-modules.
\end{enumerate}
\end{prop}

In the next result we show the one direction of Proposition~\ref{proppropertiesDelta} (ii) without assuming $\Lambda$ to be Gorenstein.

\begin{lem}
\label{lemgorprojinvere}
Let $\Lambda$ be an Artin algebra and let $\Delta_{(\phi,\phi)}=\bigl(\begin{smallmatrix}
\Lambda & \Lambda \\
\Lambda & \Lambda
\end{smallmatrix}\bigr)$. If $(X,Y,f,g)$ is an object in $\Gproj{\Delta_{(\phi,\phi)}}$ then the $\Lambda$-modules $X$ and $Y$ lie in $\Gproj{\Lambda}$.
\begin{proof}
Let $(X,Y,f,g)$ be a Gorenstein-projective $\Delta_{(\phi,\phi)}$-module. Then from Proposition~\ref{prop:projmod}, there exists a totally acyclic complex of the following form$\colon$
\[
\xymatrix@C=0.5cm{
 \mathsf{T}^{\bullet}\colon  & \cdots  \ar[rr]^{} && \mt_{1}(P^{-1})\oplus \mt_{2}(Q^{-1}) \ar@{-->}[rr]^{} \ar@{->>}[dr]_{} & & \mt_{1}(P^0)\oplus \mt_{2}(Q^0) \ar[rr]^{} && \cdots  \\
 &   & & & (X,Y,f,g) \ \ar@{>->}[ur]_{ } & &&  }
\]
where $P^i$ and $Q^i$ are projective $\Lambda$-modules. Then, if we apply the exact functor $\mU_{1}\colon \smod{\Delta_{(\phi,\phi)}}\lxr \smod{\Lambda}$, we get the exact sequence of projective $\Lambda$-modules$\colon$
\[
\xymatrix@C=0.5cm{
\mathsf{P}^{\bullet}\colon  & \cdots  \ar[rr]^{} && P^{-1}\oplus (\Lambda\otimes_{\Lambda }Q^{-1}) \ar@{-->}[rr]^{} \ar@{->>}[dr]_{} & & P^0\oplus (\Lambda\otimes_{\Lambda }Q^0) \ar[rr]^{} && \cdots  \\
 &   & & & X \ \ar@{>->}[ur]_{ } & &&  }
\]
We claim that the above complex is totally acyclic. Let $P$ be a projective $\Lambda$-module. Then from Example~\ref{examrecollementofDelta}, we have the following isomorphisms$\colon$
\[
\Hom_{\Lambda}(P^i\oplus (\Lambda\otimes_{\Lambda} Q^i),P) \iso \Hom_{\Delta_{(\phi,\phi)}}\big(\mt_{1}(P^{i})\oplus \mt_{2}(Q^{i}), \mh_{1}(P)\big) \iso \Hom_{\Delta_{(\phi,\phi)}}\big(\mt_{1}(P^{i})\oplus \mt_{2}(Q^{i}), \mt_{2}(P)\big)
\]
Since the complex $\Hom_{\Delta_{(\phi,\phi)}}(\mathsf{T}^{\bullet},\mt_{2}(P))$ is acyclic, it follows from the above isomorphisms that the complex $\Hom_{\Lambda}(\mathsf{P}^{\bullet},P)$ is also acyclic. We infer that the complex $\mathsf{P}^{\bullet}$ is totally acyclic and therefore the $\Lambda$-module $X$ is Gorenstein-projective. Similarly we show that $Y$ is a Gorenstein-projective $\Lambda$-module.
\end{proof}
\end{lem}

As a consequence of Corollary~\ref{corliftinggorproj} and Lemma~\ref{lemgorprojinvere} we obtain the following. Note that if $\Lambda$ is Gorenstein, Proposition~\ref{proppropertiesDelta} (ii) gives a direct proof of the next result.

\begin{cor}\label{corGorprojDelta}
Let $\Lambda$ be an Artin algebra and let $\Delta_{(\phi,\phi)}=\bigl(\begin{smallmatrix}
\Lambda & \Lambda \\
\Lambda & \Lambda
\end{smallmatrix}\bigr)$. Then for a $\Lambda$-module $X$ the following statements are equivalent$\colon$
\begin{enumerate}
\item $X\in \Gproj{A}$.

\item $\mt_{1}(X)\in \Gproj{\Delta_{(\phi,\phi)}}$.

\item $\mt_{2}(X)\in \Gproj{\Delta_{(\phi,\phi)}}$.
\end{enumerate}
\end{cor}

In the special case where $\phi=0$, the Gorenstein-projective modules $\mt_1(X)$ lie in the monomorphism category $\mono(\Lambda)$ as defined in subsection~\ref{subMonocat}. We close this subsection with the next example. 

\begin{exam}
Let $\mathbb{K}$ be a field and  $R=\mathbb{K}[[X_{1}, X_{2}]]/(X_{1}X_{2})$.
Consider the Morita ring $\Delta_{(0,0)}=\bigl(\begin{smallmatrix}
R & R \\
R & R
\end{smallmatrix}\bigr)$. By~\cite[Example 4.1.5]{Ch} the $R$-modules $\overline{X}_{1}$ and $\overline{X}_{2}$ are Gorenstein-projective, where
$\overline{X}_{i}$ is the residue class in $R$ of $X_{i}$ for $i=1, 2$. Thus, from Corollary~\ref{corGorprojDelta} and for $i=1,2$ it follows that the objects $\mt_{1}(\overline{X}_{i})=(\overline{X}_{i},\overline{X}_{i},\iden_{\overline{X}_{i}},0)$ and $\mt_{2}(\overline{X}_{i})=(\overline{X}_{i},\overline{X}_{i},0,\iden_{\overline{X}_{i}})$ are Gorenstein-projective $\Delta_{(0,0)}$-modules.
\end{exam}

\subsection{Constructing Gorenstein-Projective Modules}

In this subsection, we construct Gorenstein-projective modules over Morita rings $\Lambda_{(0,0)}$ which are Artin algebras and satisfy
the compatibility conditions on the bimodules $_AN_B$ and $_BM_A$, as discussed in subsection~\ref{subLifting}.

Before we proceed to the main result of this subsection (see Theorem~\ref{thmGorproj}), we need some preparations. 

\begin{lem}\label{lem:exactseq}
Let $\Lambda_{(0,0)}=\bigl(\begin{smallmatrix}
A & _AN_B \\
_BM_A & B
\end{smallmatrix}\bigr)$ be a Morita ring. Then for every $A$-module $X$ and $B$-module $Y$ we have the following exact sequences in $\Mod{\Lambda_{(0,0)}}\colon$
\[
\left\{
  \begin{array}{ll}
   \xymatrix{
  0 \ar[r]_{}^{} & \mz_{B}(M\otimes_{A}X) \ar[r]^{} & \mt_{A}(X) \ar[r]^{} & \mz_{A}(X) \ar[r] & 0 }  & \hbox{} \\
     & \hbox{} \\
  \xymatrix{
  0 \ar[r]_{}^{} & \mz_{A}(N\otimes_{B}Y) \ar[r]^{} & \mt_{B}(Y) \ar[r]^{} & \mz_{B}(Y) \ar[r] & 0 }  & \hbox{}
  \end{array}
\right.
\]
and
\[
\left\{
  \begin{array}{ll}
   \xymatrix{
  0 \ar[r]_{}^{} & \mz_{A}(X) \ar[r]^{} & \mh_{A}(X) \ar[r]^{} & \mz_{B}(\Hom_A(N,X)) \ar[r] & 0 }  & \hbox{} \\
     & \hbox{} \\
  \xymatrix{
  0 \ar[r]_{}^{} & \mz_{B}(Y) \ar[r]^{} & \mh_{B}(Y) \ar[r]^{} & \mz_{A}(\Hom_B(M,Y)) \ar[r] & 0 }  & \hbox{}
  \end{array}
\right.
\]
\begin{proof}
Let $X$ be an $A$-module. Then the map $(\iden_{X},0)\colon \mt_{A}(X) \lxr \mz_{A}(X)$ is an epimorphism in the category $\Mod{\Lambda_{(0,0)}}$, where $\mt_{A}(X)=(X,M\otimes_{A}X,\iden_{M\otimes_{A}X},0)$, $\mz_{A}(X)=(X,0,0,0)$, and the kernel of the morphism $(\iden_{X},0)$ is the object $\mz_{B}(M\otimes_{A}X)=(0,M\otimes_{A}X,0,0)$. We infer that the sequence $ 0 \lxr \mz_{B}(M\otimes_{A}X) \lxr \mt_{A}(X) \lxr \mz_{A}(X) \lxr 0$ is exact. In the same way we derive that the rest sequences are exact, the details are left to the reader.
\end{proof}
\end{lem}

\begin{lem}\label{lem:isomadjoint}
Let $\Lambda_{(0,0)}$ be a Morita ring. Then for every $X, X'\in \Mod{A}$ and $Y, Y'\in \Mod{B}$, we have the following isomorphisms$\colon$
\[
\Hom_{\Lambda_{(0,0)}}\big(\mt_{A}(X)\oplus \mt_{B}(Y),\mz_{A}(X')\big) \ \iso \ \Hom_{A}(X,X')
\]
and
\[
\Hom_{\Lambda_{(0,0)}}\big(\mt_{A}(X)\oplus \mt_{B}(Y),\mz_{B}(Y')\big) \ \iso \ \Hom_{B}(Y,Y^{\prime})
\]
\begin{proof}
We show the first isomorphism. From Proposition~\ref{prop:propertiesrecollement}, we have the adjoint pair $(\mQ_{A},\mz_{A})$ and from the recollement $(\Mod{A},\Mod{\Lambda_{(0,0)}},\Mod{B})$ it follows that $\mQ_{A}\mt_B=0$. Then, we have the isomorphism
\[
\Hom_{\Lambda_{(0,0)}}\big(\mt_{A}(X)\oplus \mt_{B}(Y),\mz_{A}(X')\big) \iso \Hom_{A}\big(\mQ_{A}\mt_{A}(X), X' \big)
\]
and it remains to compute the object $\mQ_{A}\mt_{A}(X)$. From the counit of the adjoint pair $(\mt_B,\mU_B)$ we have the following exact sequence in $\Mod{\Lambda_{(0,0)}}\colon$
\[
\xymatrix{
  \mt_{B}\mU_B\big(\mt_A(X) \big) \ar[rr]^{ \ \ \ \ (0,\iden_{M\otimes X})} && \mt_{A}(X) \ar[r]^{} & \mz_A\mQ_A(\mt_A(X)) \ar[r]^{} & 0 }
\]
see Proposition~\ref{prop:propertiesrecollement}, where $\mt_{B}\mU_B(\mt_A(X))=(N\otimes_{B}M\otimes_{A}X, M\otimes_AX,0,\iden_{N\otimes M\otimes X})$ and $\mz_A\mQ_A(\mt_A(X)) \iso \Coker{(0,\iden_{M\otimes X})} \iso \mz_{A}(X)$. This implies that $\mQ_{A}\mt_{A}(X)\iso X$ and therefore we have the isomorphism $\Hom_{\Lambda_{(0,0)}}(\mt_{A}(X)\oplus \mt_{B}(Y),\mz_{A}(X')) \iso \Hom_{A}(X,X')$.
The second isomorphism follows similarly by using the adjoint pair $(\mQ_{B},\mz_{B})$.
\end{proof}
\end{lem}

We are ready to prove the main result of this section which constructs Gorenstein-projective modules over Morita rings $\Lambda_{(0,0)}$. This result constitutes the first part of Theorem A presented in the Introduction.

\begin{thm}
\label{thmGorproj}
Let $\Lambda_{(0,0)}$ be a Morita ring which is an Artin algebra such that the bimodules $_AN_B$ and $_BM_A$ satisfy the compatibility conditions, that is, the following conditions hold$\colon$
\begin{enumerate}
\item The functor $M\otimes_{A}-\colon \smod{A}\lxr \smod{B}$ sends acyclic complexes of projective $A$-modules to acyclic complexes of $B$-modules.
\item $\add{_AN}\subseteq (\Gproj{A})^{\bot}$.

\item The functor $N\otimes_{B}-\colon \smod{B}\lxr \smod{A}$ sends acyclic complexes of projective $B$-modules to acyclic complexes of $A$-modules.

\item $\add{_BM}\subseteq (\Gproj{B})^{\bot}$.
\end{enumerate}
$(\alpha)$ Assume that there exists a Gorenstein-projective $B$-module $Z$
with a monomorphism $s\colon N\otimes_BZ\lxr X$, for some $A$-module $X$, such that $\Coker{s}$ lies in $\Gproj{A}$ and there is a monomorphism $t\colon M\otimes_{A}\Coker{s}\lxr Y$ with $\Coker{t}=Z$ and for some $B$-module $Y$. Then the tuple
\begin{equation}\label{gorprojone}
\big(X,Y, (\iden_{M}\otimes \pi_X)\circ t, (\iden_N\otimes \pi_Y)\circ s \big)
\end{equation}
is a Gorenstein-projective $\Lambda_{(0,0)}$-module, where $\pi_{X}\colon X\lxr \Coker{s}$ and $\pi_{Y}\colon Y\lxr \Coker{t}$.

\noindent $(\beta)$ Assume that there exists a Gorenstein-projective $A$-module $Z$ 
with a monomorphism $t\colon M\otimes_AZ\lxr Y$, for some $B$-module $Y$, such that $\Coker{t}$ lies in $\Gproj{B}$ and there is a monomorphism $s\colon N\otimes_{B}\Coker{t}\lxr X$ with $\Coker{s}=Z$ and for some $A$-module $X$. Then the tuple
\begin{equation}\label{gorprojtwo}
\big(X,Y, (\iden_{M}\otimes \pi_X)\circ t, (\iden_N\otimes \pi_Y)\circ s \big)
\end{equation}
is a Gorenstein-projective $\Lambda_{(0,0)}$-module, where $\pi_{X}\colon X\lxr \Coker{s}$ and $\pi_{Y}\colon Y\lxr \Coker{t}$.
\begin{proof}
We prove $(\alpha)$, the statement $(\beta)$ is dual.
The proof for $(\alpha)$ is divided into four steps. In the first two steps we construct (co)resolutions of $X$ and $Y$ by objects coming from the totally acyclic complexes of $\Coker{s}$ and $Z$. Then in the third step we lift this (co)resolutions to $\smod{\Lambda_{(0,0)}}$ and in the final step we show that this construction is indeed a totally acyclic complex of the object $(X,Y, (\iden_{M}\otimes \pi_X)\circ t, (\iden_N\otimes \pi_Y)\circ s )$.

\textsf{Step} $1\colon$
Since the $A$-module $\Coker{s}$ is Gorenstein-projective, there exists a totally acyclic complex of projective $A$-modules$\colon$
\[
\mathsf{P}^{\bullet}\colon \ \
\xymatrix{
  \cdots \ar[r]_{}^{} & P^{-1} \ar[r]^{d^{-1}_P} & P^0 \ar[r]^{d^0_P} & P^{1} \ar[r] & \cdots }
\]
such that $\Ker{d^0_P}=\Coker{s}$ and let $d^{-1}_P=\lambda^{-1}_P\circ \kappa^{-1}_P$ be the canonical factorization through $\Coker{s}$. Also, since the $B$-module $Z$ is Gorenstein-projective there exists a totally acyclic complex of projective $B$-modules$\colon$
\[
\mathsf{Q}^{\bullet}\colon \ \
\xymatrix{
  \cdots \ar[r]_{}^{} & Q^{-1} \ar[r]^{d^{-1}_Q} & Q^0 \ar[r]^{d^0_Q} & Q^{1} \ar[r] & \cdots }
\]
such that $\Ker{d^0_Q}=Z$ and let $d^{-1}_Q=\lambda^{-1}_Q\circ \kappa^{-1}_Q$ be the canonical factorization through $Z$. Then from the assumption (iii), it follows that the complex of $A$-modules $N\otimes_B \mathsf{Q}^{\bullet}$ is acyclic. Applying to the exact sequence$\colon$
\[
\xymatrix{
  0 \ar[r]^{} & N\otimes_BZ \ar[r]^{ \ \ \ s} & X \ar[r]^{\pi_X \ \ \ \ } & \Coker{s} \ar[r] & 0 }
\]
the functor $\Hom_{A}(-,N\otimes_BQ^0)$, we get the exact sequence$\colon$
\[
\xymatrix{
  0 \ar[r]^{} & \Hom_{A}(\Coker{s},N\otimes_BQ^0)  \ar[r]^{} & \Hom_{A}(X,N\otimes_BQ^0) \ar[r]^{} & \Hom_{A}(N\otimes_BZ,N\otimes_BQ^0) \ar[r] & 0 }
\]
since $\Coker{s}\in \Gproj{A}$ and $N\otimes_BQ^0\in (\Gproj{A})^{\bot}$ from the assumption (ii). This implies that there is a map $\gamma^0\colon X\lxr N\otimes_BQ^0$ such that $s\circ \gamma^0=\iden_N\otimes \kappa^{-1}_Q$ and therefore we obtain the map
\[
\alpha^0:=\bigl(\begin{smallmatrix}
 \pi_X\circ \kappa_P^{-1} \ & \gamma^0 \\
\end{smallmatrix}\bigr)\colon X\lxr P^0\oplus (N\otimes_BQ^0)
\]
Then from the Horseshoe Lemma, see also \cite[Lemma $1.6$]{Zhang}, we obtain the exact commutative diagram$\colon$
\[
\xymatrix{
   & 0 \ar[d] & 0 \ar[d] & 0 \ar[d] & \\
  0 \ar[r] &  N\otimes_BZ \ar[d]^{s} \ar[r]^{\iden_{N}\otimes \kappa_Q^{-1}} & N\otimes_{B}Q^0 \ar[d]^{\bigl(\begin{smallmatrix}
 0 & 1 \\
\end{smallmatrix}\bigr)}
  \ar[r]^{\iden_{N}\otimes d^0_Q} & N\otimes_B Q^1 \ar[d]^{\bigl(\begin{smallmatrix}
 0 & 1 \\
\end{smallmatrix}\bigr)} \ar[r] & \cdots  \\
 0 \ar@{-->}[r] &  X \ar@{-->}[ur]^{\gamma^0} \ar@{-->}[r]^{\alpha^0 \ \ \ \ \ \ } \ar[d]^{\pi_{X}} & P^0\oplus(N\otimes_BQ^0) \ar@{-->}[r]^{\alpha^1} \ar[d]^{\bigl(\begin{smallmatrix}
 1 \\
0
\end{smallmatrix}\bigr)} & P^1\oplus(N\otimes_BQ^1) \ar@{-->}[r] \ar[d]^{\bigl(\begin{smallmatrix}
 1 \\
0
\end{smallmatrix}\bigr)} & \cdots   \\
 0 \ar[r] &  \Coker{s} \ar[d] \ar[r]^{\kappa^{-1}_P} & P^0 \ar[d] \ar[r]^{d^0_P} & P^1 \ar[d] \ar[r] & \cdots  \\
  & 0 & 0 & 0 &  }
\]
where for all $i\geq 1$ we have $\alpha^i = \bigl(\begin{smallmatrix}
 d^{i-1}_P & 0 \\
 \gamma^i & \iden_N\otimes d^{i-1}_Q
\end{smallmatrix}\bigr)\colon P^{i-1}\oplus (N\otimes_BQ^{i-1})\lxr P^{i}\oplus (N\otimes_BQ^{i})$ and $\gamma^i\colon P^{i-1}\lxr N\otimes_BQ^{i}$. Note that the existence of the maps $\gamma^i$ follow by using the assumption (ii). In the same way, we can construct a resolution of $X$ by objects of the form $P^i\oplus (N\otimes_BQ^i)$ but now we use that the modules $P^{-i}, i\geq 1$, are projective.
In particular, we get the map
\[
\alpha^{-1} = \bigl(\begin{smallmatrix}
 \gamma^{-1} \\
 (\iden_N\otimes \lambda^{-1}_Q)\circ s
\end{smallmatrix}\bigr)\colon P^{-1}\oplus (N\otimes_BQ^{-1})\lxr X
\]
where $\gamma^{-1}\colon P^{-1}\lxr X$ such that $\gamma^{-1}\circ \pi_X=\lambda^{-1}_P$, and for every $i\geq 2$ we have the maps$\colon$
\[
\alpha^{-i} = \bigl(\begin{smallmatrix}
 d^{-i}_P & 0 \\
 \gamma^{-i} & \iden_N\otimes d^{-i}_Q
\end{smallmatrix}\bigr)\colon P^{-i}\oplus (N\otimes_BQ^{-i})\lxr P^{-i+1}\oplus (N\otimes_BQ^{-i+1})
\]
similarly as described above. Thus, summarizing the construction so far, 
we have constructed the following exact sequence$\colon$
\[
\xymatrix@C=0.4cm{
  \cdots \ar[r]^{} & P^{-2}\oplus(N\otimes_BQ^{-2})  \ar[r]^{\alpha^{-2}} & P^{-1}\oplus(N\otimes_BQ^{-1}) \ar@{-->}[rr]^{} \ar@{->>}[dr]_{\alpha^{-1}} & & P^0\oplus(N\otimes_BQ^0) \ar[r]^{\alpha^1} & P^1\oplus(N\otimes_BQ^1) \ar[r] & \cdots  \\
   & & & X \ \ar@{>->}[ur]_{\alpha^0} & & & (*)  }
\]

\textsf{Step} $2\colon$ We construct an exact sequence similar to $(*)$ for the $B$-module $Y$.  Since $\add{_BM}\subseteq (\Gproj{B})^{\bot}$ (assumption (iv)) we have as in \textsf{Step} $1$ the following exact commutative diagram$\colon$
\[
\xymatrix{
   & 0 \ar[d] & 0 \ar[d] & 0 \ar[d] & \\
  0 \ar[r] &  M\otimes_A\Coker{s} \ar[d]^{t} \ar[r]^{ \ \ \iden_{M}\otimes \kappa_P^{-1}} & M\otimes_{A}P^0 \ar[d]^{\bigl(\begin{smallmatrix}
 1 & 0 \\
\end{smallmatrix}\bigr)}
  \ar[r]^{\iden_{M}\otimes d^0_P} & M\otimes_A P^1 \ar[d]^{\bigl(\begin{smallmatrix}
 1 & 0 \\
\end{smallmatrix}\bigr)} \ar[r] & \cdots  \\
 0 \ar@{-->}[r] &  Y \ar@{-->}[ur]^{\delta^0} \ar@{-->}[r]^{\beta^0 \ \ \ \ \ \ } \ar[d]^{\pi_{Y}} & (M\otimes_AP^0)\oplus Q^0 \ar@{-->}[r]^{\beta^1} \ar[d]^{\bigl(\begin{smallmatrix}
 0 \\
1
\end{smallmatrix}\bigr)} & (M\otimes_AP^1)\oplus Q^1 \ar@{-->}[r] \ar[d]^{\bigl(\begin{smallmatrix}
 0 \\
1
\end{smallmatrix}\bigr)} & \cdots & (**) \\
 0 \ar[r] &  Z \ar[d] \ar[r]^{\kappa^{-1}_Q} & Q^0 \ar[d] \ar[r]^{d^0_Q} & Q^1 \ar[d] \ar[r] & \cdots  \\
  & 0 & 0 & 0 &  }
\]
where $\beta^0:=\bigl(\begin{smallmatrix}
 \delta^0 \ & \pi_Y\circ \kappa_{Q}^{-1} \\
\end{smallmatrix}\bigr)\colon Y\lxr (M\otimes_AP^0)\oplus Q^0$
and for all $i\geq 1$ we have$\colon$
\[
\beta^i = \bigl(\begin{smallmatrix}
 \iden_{M}\otimes d^{i-1}_P & \delta^i \\
 0 & d^{i-1}_Q
\end{smallmatrix}\bigr)\colon (M\otimes_A P^{i-1})\oplus Q^{i-1}\lxr (M\otimes_A P^{i})\oplus Q^{i}
\]
for some $\delta^i\colon Q^{i-1}\lxr M\otimes_A P^{i}$. Then, as in \textsf{Step} $1$ we  construct a resolution of $Y$ by objects of the form $(M\otimes_AP^i)\oplus Q^i$ and putting together these, we obtain the following exact sequence$\colon$
\[
\xymatrix@C=0.4cm{
  \cdots \ar[r]^{} & (M\otimes_AP^{-2})\oplus Q^{-2}  \ar[r]^{\beta^{-2}} & (M\otimes_AP^{-1})\oplus Q^{-1} \ar@{-->}[rr]^{} \ar@{->>}[dr]_{\beta^{-1}} & & (M\otimes_AP^0)\oplus Q^0 \ar[r]^{\beta^1} & (M\otimes_AP^1)\oplus Q^1 \ar[r] & \cdots  \\
   & & & Y \ \ar@{>->}[ur]_{\beta^0} & & & (**)  }
\]

\textsf{Step} $3\colon$ We glue together the exact sequences $(*)$ and $(**)$ and we derive the following sequence$\colon$
\[
\mathsf{T}^{\bullet}\colon \ \
\xymatrix@C=0.5cm{
  \cdots  \ar[rr]^{(\alpha^{-2}, \beta^{-2}) \ \ \ \ \ \ \ \ \ \ \ \ \ } && \mt_{A}(P^{-1})\oplus \mt_{B}(Q^{-1}) \ar@{-->}[rr]^{} \ar@{->>}[dr]_{(\alpha^{-1}, \beta^{-1})} & & \mt_{A}(P^0)\oplus \mt_{B}(Q^0) \ar[rr]^{ \ \ \ \ \ \ \ \ (\alpha^1,\beta^1)} && \cdots  \\
    & & & (X,Y,f,g) \ \ar@{>->}[ur]_{(\alpha^0, \beta^0)} & &&   }
\]
We claim that the sequence $\mathsf{T}^{\bullet}$ is exact in $\smod{\Lambda_{(0,0)}}$. First, since the following diagrams are commutative
\[
\xymatrix{
  (M\otimes_AP^i)\oplus (M\otimes_AN\otimes_BQ^i) \ar[d]_{\bigl(\begin{smallmatrix}
 \iden_{M}\otimes d^i_P & 0 \\
 \iden_{M}\otimes \gamma^{i+1} & \iden_{M\otimes N}\otimes d^i_Q
\end{smallmatrix}\bigr)} \ar[rr]^{ \ \ \ \ \ \bigl(\begin{smallmatrix}
 \iden_{M\otimes P^i} & 0 \\
 0 & 0
\end{smallmatrix}\bigr)} && (M\otimes_AP^i)\oplus Q^i \ar[d]^{\bigl(\begin{smallmatrix}
 \iden_{M}\otimes d^i_P & \delta^{i+1} \\
 0 & d^i_Q
\end{smallmatrix}\bigr)} \\
 (M\otimes_AP^{i+1})\oplus (M\otimes_AN\otimes_BQ^{i+1}) \ar[rr]^{ \ \ \ \ \ \ \ \ \ \ \bigl(\begin{smallmatrix}
 \iden_{{M}\otimes P^{i+1}} & 0 \\
 0 & 0
\end{smallmatrix}\bigr)} && (M\otimes_AP^{i+1})\oplus Q^{i+1}  }
\]
and
\[
\xymatrix{
  (N\otimes_BM\otimes_AP^i)\oplus (N\otimes_BQ^i) \ar[d]_{\bigl(\begin{smallmatrix}
 \iden_{N\otimes M}\otimes d^i_P & \iden_{N}\otimes \delta^{i+1} \\
 0 & \iden_{N}\otimes d^i_Q
\end{smallmatrix}\bigr)} \ar[rr]^{ \ \ \ \ \ \ \bigl(\begin{smallmatrix}
 0 & 0 \\
 0 & \iden_{{N}\otimes Q^i}
\end{smallmatrix}\bigr)} && P^{i}\oplus (N\otimes_BQ^{i}) \ar[d]^{\bigl(\begin{smallmatrix}
 d^i_P & 0 \\
 \gamma^{i+1} & \iden_{N}\otimes d^i_Q
\end{smallmatrix}\bigr)} \\
  (N\otimes_BM\otimes_AP^{i+1})\oplus (N\otimes_BQ^{i+1}) \ar[rr]^{ \ \ \ \ \ \ \ \ \ \ \bigl(\begin{smallmatrix}
 0 & 0 \\
 0 & \iden_{N\otimes Q^{i+1}}
\end{smallmatrix}\bigr)} && P^{i+1}\oplus (N\otimes_BQ^{i+1})  }
\]
it follows that the maps
\[
(\alpha^i, \beta^i)\colon \mt_{A}(P^i)\oplus \mt_{B}(Q^i) \lxr \mt_{A}(P^{i+1})\oplus \mt_{B}(Q^{i+1})
\]
are morphisms in $\smod{\Lambda_{(0,0)}}$. Since the complexes $(*)$ and $(**)$ are acyclic, it follows from Remark~\ref{remMoritarings} (iii) that $\mathsf{T}^{\bullet}$ is an exact sequence in $\smod{\Lambda_{(0,0)}}$. Moreover, the object $(X,Y,f,g)$ arises as the kernel of the morphism $(\alpha^1,\beta^1)$ and
by Remark~\ref{remMoritarings} (iv) we observe that $f=(\iden_{M}\otimes \pi_X)\circ t$ and $g=(\iden_N\otimes \pi_Y)\circ s$.

\textsf{Step} $4\colon$ The final step of the proof is devoted to show that the acyclic complex $\mathsf{T}^{\bullet}$ is totally acyclic. From Proposition~\ref{prop:projmod}, it is enough to show that the complexes $\Hom_{\Lambda_{(0,0)}}(\mathsf{T}^{\bullet},\mt_{A}(P))$ and $\Hom_{\Lambda_{(0,0)}}(\mathsf{T}^{\bullet},\mt_{B}(Q))$ are acyclic,  where $P$ is a projective $A$-module and $Q$ is a projective $B$-module. 
From Lemma~\ref{lem:exactseq} we have the exact sequence $0 \lxr \mz_{B}(M\otimes_{A}P) \lxr \mt_{A}(P) \lxr \mz_{A}(P) \lxr 0$
and since each term of the complex $\mathsf{T}^{\bullet}$ is a projective $\Lambda_{(0,0)}$-module, it follows that the following sequence$\colon$
\begin{equation}\label{exactseqcomplexes}
\xymatrix@C=0.5cm{
  0 \ar[r]_{}^{} & \Hom_{\Lambda_{(0,0)}}(\mathsf{T}^{\bullet}, \mz_{B}(M\otimes_{A}P)) \ar[r]^{} &  \Hom_{\Lambda_{(0,0)}}(\mathsf{T}^{\bullet}, \mt_{A}(P)) \ar[r]^{} &  \Hom_{\Lambda_{(0,0)}}(\mathsf{T}^{\bullet}, \mz_{A}(P)) \ar[r] & 0 }
\end{equation}
is an exact sequence of complexes. Then, from Lemma~\ref{lem:isomadjoint} we have the isomorphism $\Hom_{\Lambda_{(0,0)}}(\mathsf{T}^{\bullet}, \mz_{A}(P))\iso \Hom_{A}(\mathsf{P}^{\bullet},P)$ and since $\mathsf{P}^{\bullet}$ is totally acyclic we infer that $\Hom_{\Lambda_{(0,0)}}(\mathsf{T}^{\bullet}, \mz_{A}(P))$ is acyclic. Also, from Lemma~\ref{lem:isomadjoint} we have $\Hom_{\Lambda_{(0,0)}}(\mathsf{T}^{\bullet}, \mz_{B}(M\otimes_{A}P))\iso \Hom_{B}(\mathsf{Q}^{\bullet},M\otimes_{A}P)$ and since $M\otimes_{A}P$ lies in $(\Gproj{B})^{\bot}$ by assumption (iv), it follows that the complex $\Hom_{B}(\mathsf{Q}^{\bullet},M\otimes_{A}P)$ is acyclic. Then, the complex $\Hom_{\Lambda_{(0,0)}}(\mathsf{T}^{\bullet}, \mz_{B}(M\otimes_{A}P))$ is acyclic and therefore from the exact sequence $(\ref{exactseqcomplexes})$, we infer that the  complex $\Hom_{\Lambda_{(0,0)}}(\mathsf{T}^{\bullet}, \mt_{A}(P))$ is acyclic. Similarly, using the exact sequence $0 \lxr \mz_{A}(N\otimes_{B}Q) \lxr \mt_{B}(Q) \lxr \mz_{B}(Q) \lxr 0$ we derive that the complex $\Hom_{\Lambda_{(0,0)}}(\mathsf{T}^{\bullet},\mt_{B}(Q))$ is acyclic.

In conclusion, the $\Lambda_{(0,0)}$-module $(X,Y, (\iden_{M}\otimes \pi_X)\circ t, (\iden_N\otimes \pi_Y)\circ s)$ is Gorenstein-projective.
\end{proof}
\end{thm}

\begin{cor}
\label{corcompatibilitycondpdfinite}
Let $\Lambda_{(0,0)}$ be a Morita ring which is an Artin algebra such that the conditions $(1)$ or $(3)$, and $(2)$ or $(4)$ hold$\colon$

\smallskip

\begin{minipage}[h]{7.2cm}
\begin{enumerate}
\item[$(1)$] $\pd{M_A}<\infty$ \ and \ $\pd{_AN}<\infty$.
\item[$(2)$] $\pd{N_B}<\infty$ \ and \ $\pd{_BM}<\infty$.
\end{enumerate}
\end{minipage}
\begin{minipage}[h]{7.2cm}
\begin{enumerate}
\item[$(3)$] $\pd{M_A}<\infty$ \ and \ $\id{_AN}<\infty$.
\item[$(4)$] $\pd{N_B}<\infty$ \ and \ $\id{_BM}<\infty$.
\end{enumerate}
\end{minipage}

\smallskip

\noindent $(\alpha)$ Assume that for an $A$-module $X$ there exists an exact sequence 
\[
\xymatrix{
  0 \ar[r]_{}^{} & N\otimes_{B}Z \ar[r]^{ \ \ \ s} & X \ar[r]^{\pi_X \ \ \ \ } & \Coker{s} \ar[r] & 0 }
\]
with $Z\in \Gproj{B}$ and $\Coker{s}\in \Gproj{A}$, such that there is an exact sequence 
\[
\xymatrix{
  0 \ar[r]_{}^{} & M\otimes_{A}\Coker{s} \ar[r]^{ \ \ \ \ \ \ t} & Y \ar[r]^{\pi_Y  } & Z \ar[r] & 0 }
\]
for some $B$-module $Y$. Then the objects$\colon$$(X,Y, (\iden_{M}\otimes \pi_X)\circ t, (\iden_N\otimes \pi_Y)\circ s)$, $\mt_A(\Coker{s})$, $\mt_B(Z)$ are Gorenstein-projective $\Lambda_{(0,0)}$-modules.

\noindent $(\beta)$ Assume that for a $B$-module $Y$ there exists an exact sequence 
\[
\xymatrix{
  0 \ar[r]_{}^{} & M\otimes_{A}Z \ar[r]^{ \ \ \ t} & Y \ar[r]^{\pi_Y \ \ \ \ } & \Coker{t} \ar[r] & 0 }
\]
with $Z\in \Gproj{A}$ and $\Coker{t}\in \Gproj{B}$, such that there is an exact sequence 
\[
\xymatrix{
  0 \ar[r]_{}^{} & N\otimes_{B}\Coker{t} \ar[r]^{ \ \ \ \ \ \ s} & X \ar[r]^{\pi_X  } & Z \ar[r] & 0 }
\]
for some $A$-module $X$. Then the objects$\colon$$(X,Y, (\iden_{M}\otimes \pi_X)\circ t, (\iden_N\otimes \pi_Y)\circ s)$, $\mt_A(Z)$, $\mt_B(\Coker{t})$ are Gorenstein-projective $\Lambda_{(0,0)}$-modules.
\begin{proof}
From Example~\ref{examconditionsformainthm} the conditions (i) -- (iv) of Theorem~\ref{thmGorproj} are satisfied. Then the result follows from Corollary~\ref{corliftinggorproj} and Theorem~\ref{thmGorproj}. 
\end{proof}
\end{cor}

The next result is a consequence of Theorem~\ref{thmGorproj} for the Morita ring $\Delta_{(0,0)}$. Recall that $\smod{\Delta_{(0,0)}}$ is the double morphism category $\DMor(\smod\Lambda)$ that we studied in subsection~\ref{subsectiondoulemorphism}.

\begin{cor}
\label{corexamgorprojexactseq}
Let $\Lambda$ be an Artin algebra and consider the algebra $\Delta_{(0,0)}=\bigl(\begin{smallmatrix}
\Lambda & \Lambda \\
\Lambda & \Lambda
\end{smallmatrix}\bigr)$. Let $(X,Y,f,g)$ be a $\Delta_{(0,0)}$-module such that there exist exact sequences
\[
\xymatrix{
  0 \ar[r]_{}^{} & Z \ar[r]^{s} & X \ar[r]^{\pi_{X} \ } & W \ar[r] & 0 }
\]
and
\[
\xymatrix{
  0 \ar[r]_{}^{} & W \ar[r]^{t} & Y \ar[r]^{\pi_{Y}} & Z \ar[r] & 0 }
\]
with $Z, W\in \Gproj{\Lambda}$ and set $f:=\pi_{X}\circ t$, $g:=\pi_{Y}\circ s$. Then the objects $(X,Y,f,g)$ and $(Y,X,g,f)$ are Gorenstein-projective $\Delta_{(0,0)}$-modules.
\end{cor}

\begin{rem}
Let $(X,Y,f,g)$ be a $\Lambda_{(0,0)}$-module. Assume that for $X$ and $Y$ the conditions of Theorem~\ref{thmGorproj} $(\alpha)$ are satisfied. Then, we cannot infer in general from Theorem~\ref{thmGorproj} that $(X,Y,f,g)$ lies in $\Gproj{\Delta_{(0,0)}}$. In other words, Theorem~\ref{thmGorproj} does not provide us with sufficient conditions for a tuple $(X,Y,f,g)$ to be Gorenstein-projective. We explain now where is the problem. Following the construction of Theorem~\ref{thmGorproj}, we conclude that the object $(X,Y, (\iden_{M}\otimes \pi_X)\circ t, (\iden_N\otimes \pi_Y)\circ s)$ is Gorenstein-projective. From Remark~\ref{remMoritarings} (v) we know that the maps $(\iden_{M}\otimes \pi_X)\circ t$ and $(\iden_N\otimes \pi_Y)\circ s$ are uniquely determined and satisfy the corresponding commutative diagrams $(\ref{diagramskernels})$. But since $f$ and $g$ are arbitrary maps, we don't know in general if they satisfy the diagrams $(\ref{diagramskernels})$. If $f$ and $g$ satisfy these diagrams, then from uniqueness it follows that $f=(\iden_{M}\otimes \pi_X)\circ t$, $g=(\iden_N\otimes \pi_Y)\circ s$ and therefore $(X,Y,f,g)$ is Gorenstein-projective.
Hence, we cannot conclude from Theorem~\ref{thmGorproj} that $(X,Y,f,g)$ is Gorenstein-projective.
\end{rem}

The next example shows how we can apply Corollary~\ref{corexamgorprojexactseq} to construct Gorenstein-projective modules over $\Delta_{(0,0)}$ from Gorenstein-projective modules over the underlying triangular matrix algebras. 

\begin{exam}
\label{examconstructgorproj}
Let $\Lambda$ be an Artin algebra and consider the Morita ring $\Delta_{(0,0)}=\bigl(\begin{smallmatrix}
\Lambda & \Lambda \\
\Lambda & \Lambda
\end{smallmatrix}\bigr)$.
\begin{enumerate}
\item Consider the lower triangular matrix algebra $\mathsf{T}_2(\Lambda)=\bigl(\begin{smallmatrix}
\Lambda & 0 \\
\Lambda & \Lambda
\end{smallmatrix}\bigr)$. From \cite[Theorem $1.4$]{Zhang}, a triple $(X,Y,f)$ is a Gorenstein-projective $\Gamma$-module if and only if there is an exact sequence
\begin{equation}
\label{equationexample}
\xymatrix{
  0 \ar[r]_{}^{} & X \ar[r]^{f} & Y \ar[r]^{\pi \ \ \ \ } & \Coker{f} \ar[r] & 0 }
\end{equation}
such that the $\Lambda$-modules $X$ and $\Coker{f}$ are Gorenstein-projective. Let $(X,Y,f)$ be a Gorenstein-projective $\Gamma$-module. Thus, we have the sequence $(\ref{equationexample})$ and we also form the split exact sequence$\colon$
\[
\xymatrix{
  0 \ar[r]_{}^{} & \Coker{f} \ar[r]^{\bigl(\begin{smallmatrix}
1 & 0
\end{smallmatrix}\bigr) \ \ \ \ } & \Coker{f}\oplus X \ar[r]^{ \ \ \ \ \ \bigl(\begin{smallmatrix}
 0 \\
 1
\end{smallmatrix}\bigr)} & X \ar[r] & 0 }
\]
Then, Corollary~\ref{corexamgorprojexactseq} yields that the objects
\[
\big(Y, \Coker{f}\oplus X, \pi\circ \bigl(\begin{smallmatrix}
1 & 0
\end{smallmatrix}\bigr), \bigl(\begin{smallmatrix}
 0 \\
 1
\end{smallmatrix}\bigr)\circ f\big) \ \ \ \ \text{and} \ \ \ \ \big(\Coker{f}\oplus X, Y, \bigl(\begin{smallmatrix}
 0 \\
 1
\end{smallmatrix}\bigr)\circ f, \pi\circ \bigl(\begin{smallmatrix}
1 & 0
\end{smallmatrix}\bigr) \big)
\]
are Gorenstein-projective $\Delta_{(0,0)}$-modules. Consider now the upper triangular matrix algebra $\Sigma=\bigl(\begin{smallmatrix}
\Lambda & \Lambda \\
0 & \Lambda
\end{smallmatrix}\bigr)$ and let $(Z,W,g)\in \Gproj{\Sigma}$. Then, from \cite[Theorem $1.4$]{Zhang} there is an exact sequence$\colon$
\[
\xymatrix{
  0 \ar[r]_{}^{} & W \ar[r]^{g} & Z \ar[r]^{\rho \ \ \ \ } & \Coker{g} \ar[r] & 0 }
\]
such that the $\Lambda$-modules $W$ and $\Coker{g}$ lie in $\Gproj{\Lambda}$, and we also have the split exact sequence$\colon$
\[
\xymatrix{
  0 \ar[r]_{}^{} & \Coker{g} \ar[r]^{\bigl(\begin{smallmatrix}
1 & 0
\end{smallmatrix}\bigr) \ \ \ \ } & \Coker{g}\oplus W \ar[r]^{ \ \ \ \ \ \bigl(\begin{smallmatrix}
 0 \\
 1
\end{smallmatrix}\bigr)} & W \ar[r] & 0 }
\]
Hence, by Corollary~\ref{corexamgorprojexactseq} it follows that the following objects$\colon$
\[
\big(Z, \Coker{g}\oplus W, \rho\circ \bigl(\begin{smallmatrix}
1 & 0
\end{smallmatrix}\bigr), \bigl(\begin{smallmatrix}
 0 \\
 1
\end{smallmatrix}\bigr)\circ g\big) \ \ \ \ \text{and} \ \ \ \ \big(\Coker{g}\oplus W, Z, \bigl(\begin{smallmatrix}
 0 \\
 1
\end{smallmatrix}\bigr)\circ g, \rho \circ \bigl(\begin{smallmatrix}
1 & 0
\end{smallmatrix}\bigr) \big)
\]
are Gorenstein-projective $\Delta_{(0,0)}$-modules.

\item Let $X$ be a Gorenstein-projective $\Lambda$-module.  From (i) the objects $(X,X,0,\iden_X)$ and  $(X,X,\iden_X,0)$ are
Gorenstein-projective $\Delta_{(0,0)}$-modules. Note that this was also observed in Corollary~\ref{corGorprojDelta}.
\end{enumerate}
\end{exam}

The above example shows that using Theorem~\ref{thmGorproj}, we obtain non-trivial examples of Gorenstein-projective modules over the Morita ring $\Delta_{(0,0)}$ from Gorenstein-projective modules of the triangular matrix algebras $\Gamma$ and $\Sigma$. It should be noted that we don't know if all Gorenstein-projective modules over $\Delta_{(0,0)}$ arises in this way, as well as
how many objects from $\Gproj{\Delta_{(0,0)}}$ we finally obtain. We recall that in the case where $\Lambda$ is a Gorenstein Artin algebra, a tuple $(X,Y,f,g)$ lies in $\Gproj{\Delta_{(0,0)}}$ if and only if $X$ and $Y$ lie in $\Gproj{\Lambda}$, see Proposition~\ref{proppropertiesDelta} (ii). But also in this case it seems to be a difficult problem to determine the Gorenstein-projectives $\Delta_{(0,0)}$-modules. On the other hand, the above example shows the connections of the module categories of $\Gamma$ and $\Sigma$ with the module category of $\Delta_{(0,0)}$. This fact will become clear in the next subsection, see Lemma~\ref{lemadjointtriangularmorita}.

We close this subsection with the following consequence of Corollary~\ref{corexamgorprojexactseq} and an example. We mention that Example~\ref{examstronglygorproj} provides an interesting connection between our main result (Theorem~\ref{thmGorproj}) and the class of strongly Gorenstein-projective modules.

\begin{cor}
\label{corollarychargorproj}
Let $\Lambda$ be an Artin algebra and consider the algebra $\Delta_{(0,0)}=\bigl(\begin{smallmatrix}
\Lambda & \Lambda \\
\Lambda & \Lambda
\end{smallmatrix}\bigr)$. Let $(X,Y,f,g)$ be a $\Delta_{(0,0)}$-module such that $\Image{f}=\Ker{g}$, $\Image{g}=\Ker{f}$ and assume that $\Image{f}$ lies in $\Gproj{\Lambda}$. Then $(X,Y,f,g)\in \Gproj{\Delta_{(0,0)}}$ if and only if $X, Y\in \Gproj{\Lambda}$ if and only if $(Y,X,g,f)\in \Gproj{\Delta_{(0,0)}}$.
\begin{proof}
Suppose first that $X$ and $Y$ are Gorenstein-projective $\Lambda$-modules. Then, from our assumptions the following complex$\colon$
\[
\xymatrix{
  \cdots \ar[r]_{}^{} & X \ar[r]^{f} & Y \ar[r]^{g} & X \ar[r]^{f} & Y \ar[r] & \cdots }
\]
is acyclic. Thus, we have the following short exact sequences$\colon$
\[
\xymatrix{
  0 \ar[r]_{}^{} & \Image{g} \ar[r]^{} & X \ar[r]^{} & \Image{f} \ar[r] & 0 }
\]
and
\[
\xymatrix{
  0 \ar[r]_{}^{} & \Image{f} \ar[r]^{} & Y \ar[r]^{} & \Image{g} \ar[r] & 0 }
\]
Since $\Gproj{\Lambda}$ is closed under kernels of epimorphisms, it follows that $\Image{f}\in \Gproj{\Lambda}$ if and only if $\Image{g}\in \Gproj{\Lambda}$. Then, for $Z=\Image{g}$ in Corollary~\ref{corexamgorprojexactseq}, we get that the module $(X,Y,f,g)$ is Gorenstein-projective. Note that, in this case, the maps of the tuple that we obtain from Corollary~\ref{corexamgorprojexactseq} are precisely $f$ and $g$. Similarly, if $Z=\Image{f}$ then the tuple $(Y,X,g,f)$ is Gorenstein-projective. The converse directions follow from Lemma~\ref{lemgorprojinvere}.
\end{proof}
\end{cor}

\begin{exam} 
\label{examstronglygorproj}
Let $\Lambda$ be an Artin algebra and consider the matrix algebra $\Delta_{(0,0)}=\bigl(\begin{smallmatrix}
\Lambda & \Lambda \\
\Lambda & \Lambda
\end{smallmatrix}\bigr)$. Let 
\[
\xymatrix{
  \cdots \ar[r]_{}^{} & P \ar[r]^{f} & P \ar[r]^{f} & P \ar[r]^{f} & P \ar[r] & \cdots }
\]
be a totally acyclic complex of projective $\Lambda$-modules. Then, from Corollary~\ref{corollarychargorproj} it follows that $(P,P,f,f)$ is a Gorenstein-projective $\Delta_{(0,0)}$-module. In this case, the $\Lambda$-module $\Image{f}$ is called strongly Gorenstein-projective. We refer to \cite{BM} for more details on this class of modules. 

As a particular example, let $\mathbb{K}$ be a field, $\Lambda=\mathbb{K}[X]/(X^{2})$ and consider the matrix algebra $\Delta_{(0,0)}$. Denote by $\overline{X}$ the residue class of $X$ in $\Lambda$. Then by \cite[Example 2.5]{BM} the following sequence
\[
\xymatrix{
  \cdots \ar[r]_{}^{} & \Lambda \ar[r]^{x} & \Lambda \ar[r]^{x} & \Lambda \ar[r]^{x} & \Lambda \ar[r] & \cdots }
\]
 is a totally acyclic complex of projective $\Lambda$-modules and $\overline{X}=\Image{x}=\Ker x$ is a strongly Gorenstein-projective $\Lambda$-module. We infer that $(\Lambda,\Lambda,x,x)$ is a Gorenstein-projective $\Delta_{(0,0)}$-module. 
\end{exam} 

\begin{rem} 
By Corollary~\ref{corollarychargorproj} we can instantly derive Example~\ref{examconstructgorproj} (i). Indeed, let $f\colon X\lxr Y$ be a monomorphism with $\Coker{f}$ in $\Gproj{\Lambda}$. Consider the maps $\bigl(\begin{smallmatrix}
1 & 0
\end{smallmatrix}\bigr)\colon \Coker{f}\lxr \Coker{f}\oplus X$ and $\bigl(\begin{smallmatrix}
 0 \\
 1
\end{smallmatrix}\bigr)\colon \Coker{f}\oplus X \lxr X$. Then by Corollary~\ref{corollarychargorproj} we get that $(Y, \Coker{f}\oplus X, \pi\circ \bigl(\begin{smallmatrix}
1 & 0
\end{smallmatrix}\bigr),  \bigl(\begin{smallmatrix}
 0 \\
 1
\end{smallmatrix}\bigr)\circ f)$ is a Gorenstein projective $\Delta_{(0,0)}$-module if and only if $Y$ and $\Coker{f}\oplus X$ are Gorenstein-projectives if and only if $Y$ and $X$ are Gorenstein-projectives. 
\end{rem}

\subsection{Ring Epimorphisms and Gorenstein-Projective Modules}
\label{subringepi}
In this subsection we address an interesting connection between Morita rings with zero bimodule maps and triangular matrix rings.

Let $\Lambda_{(0,0)}=\bigl(\begin{smallmatrix}
A & _AN_B \\
_BM_A & B
\end{smallmatrix}\bigr)$ be a Morita ring. Then we have the triangular matrix rings $\bigl(\begin{smallmatrix}
 A & 0 \\
 _BM_A & B
\end{smallmatrix}\bigr)$ and $\bigl(\begin{smallmatrix}
 A & _AN_B \\
 0 & B
\end{smallmatrix}\bigr)$. Before we define functors between the module categories of the above rings, we recall the adjoint isomorphisms$\colon$$\pi:\Hom_{B}(M\otimes_{A}X,Y) \stackrel{\simeq}{\lxr} \Hom_{A}(X,\Hom_{B}(M,Y))$ and 
$\rho:\Hom_{A}(N\otimes_{B}Y,X) \stackrel{\simeq}{\lxr} \Hom_{B}(Y,\Hom_{A}(N,X))$.

We define the following functors$\colon$
\begin{enumerate}
\item The functor $\F\colon \Mod{\bigl(\begin{smallmatrix}
 A & 0 \\
 _BM_A & B
\end{smallmatrix}\bigr)}\lxr \Mod{\Lambda_{(0,0)}}$ is defined by
$\F(X,Y,f)=(X,Y,f,0)$ on the objects $(X,Y,f)\in \Mod{\bigl(\begin{smallmatrix}
 A & 0 \\
 _BM_A & B
\end{smallmatrix}\bigr)}$ and given a homomorphism $(a,b)\colon (X,Y,f)\lxr (X',Y',f')$ in $\Mod{\bigl(\begin{smallmatrix}
 A & 0 \\
 _BM_A & B
\end{smallmatrix}\bigr)}$, then
$\F(a,b)=(a,b)$.

\item The functor $\G\colon \Mod{\Lambda_{(0,0)}}\lxr \Mod{\bigl(\begin{smallmatrix}
 A & 0 \\
 _BM_A & B
\end{smallmatrix}\bigr)}$ is defined by
$\G(X,Y,f,g)=(\Coker{g},Y,h)$ on the objects $(X,Y,f,g)\in \Mod{\Lambda_{(0,0)}}$, where the morphism $h\colon M\otimes_A\Coker{g}\lxr Y$ is obtained from the following commutative diagram$\colon$
\[
\xymatrix{
  M\otimes_AN\otimes_B Y \ar@{-->}[dr]_{0} \ar[r]^{\ \ \ \  \iden_{M}\otimes g} &  M\otimes_AX \ar[d]^{f}  \ar[r]^{\iden_{M}\otimes \pi \ \ \  } & M\otimes_A\Coker{g} \ar[r] \ar[dl]^{h} & 0    \\
   & Y  & &      }
\]
If $(a,b)\colon (X,Y,f,g)\lxr (X',Y',f',g')$ is a morphism in
$\Mod{\Lambda_{(0,0)}}$, then $\G(a,b)=(\xi,b)$ where $\xi\colon \Coker{g}\lxr \Coker{g'}$ is the unique morphism such that $\pi\circ \xi=a\circ \pi'$, where $\pi\colon X\lxr \Coker{g}$ and $\pi'\colon X'\lxr \Coker{g'}$. 

\item The functor $\mathcal{H}\colon \Mod{\Lambda_{(0,0)}}\lxr \Mod{\bigl(\begin{smallmatrix}
 A & 0 \\
 _BM_A & B
\end{smallmatrix}\bigr)}$ is defined by
$\mathcal{H}(X,Y,f,g)=(X,\Ker{\rho(g)},\pi^{-1}(j))$ on the objects $(X,Y,f,g)\in \Mod{\Lambda_{(0,0)}}$, where the morphism $\pi^{-1}(j)\colon M\otimes_AX\lxr \Ker{\rho(g)}$ is obtained from the following commutative diagram$\colon$
\[
\ \ \ \ \ \ \ \ \xymatrix{
0\ar[r] & \Hom_B(M,\Ker{\rho(g)}) \ar[rr]^{ \ \ \Hom(M,i) } && \Hom_B(M,Y)  \ar[rr]^{\Hom(M,\rho(g)) \ \ \ \ \ \ \ \ } && \Hom_B(M,\Hom_A(N,X))    \\
   & && X \ar[u]_{\pi(f)} \ar[llu]^{j} \ar@{-->}[rru]_{0} &&      }
\]
If $(a,b)\colon (X,Y,f,g)\lxr (X',Y',f',g')$ is a morphism in
$\Mod{\Lambda_{(0,0)}}$, then $\mathcal{H}(a,b)=(a,\zeta)$ where $\zeta\colon \Ker{\rho(g)}\lxr \Ker{\rho(g')}$ is the unique morphism which makes the next diagram commutative$\colon$
\[
\xymatrix{
  0\ar[r] & \Ker{\rho(g)} \ar[d]_{\zeta} \ar[r]^{ \ \ i} &  Y \ar[d]^{b}  \ar[r]^{\rho(g) \ \ \ \ \ \ \ \ } & \Hom_A(N,X) \ar[d]^{\Hom(N,a)}    \\
0\ar[r] & \Ker{\rho(g')}   \ar[r]^{ \ \ i'} & Y'  \ar[r]^{\rho(g') \ \ \ \ \ \ \ \ } & \Hom_A(N,X') }
\]
\end{enumerate}
Similarly we have the functor $\F'\colon \Mod{\bigl(\begin{smallmatrix}
 A & _AN_B \\
 0 & B
\end{smallmatrix}\bigr)}\lxr \Mod{\Lambda_{(0,0)}}$ given by $\F'(X,Y,g)=(X,Y,0,g)$ where $g\colon N\otimes_BY\lxr X$, and dually with $\G$, $\mathcal{H}$ we define the functors $\G', \mathcal{H}'\colon \Mod{\Lambda_{(0,0)}}\lxr \Mod{\bigl(\begin{smallmatrix}
 A & _AN_B \\
 0 & B
\end{smallmatrix}\bigr)}$ by $\G'(X,Y,f,g)=(X,\Coker{f},h)$ and $\mathcal{H}'=(\Ker{\pi(f)},Y,\rho^{-1}(j))$. It is straightforward that $\F$, $\G$, $\mathcal{H}$ and $\F'$, $\G'$, $\mathcal{H}'$ are additive functors.

\begin{lem}
\label{lemadjointtriangularmorita}
Let $\Lambda_{(0,0)}$ be a Morita ring. Then the triples of functors $(\G,\F,\mathcal{H})$ and $(\G',\F',\mathcal{H}')\colon$
\[
\xymatrix@C=0.5cm{
\Mod{\bigl(\begin{smallmatrix}
 A & 0 \\
 _BM_A & B
\end{smallmatrix}\bigr)} \ar[rrr]^{\F} &&& \Mod{\Lambda_{(0,0)}} \ar @/_1.5pc/[lll]_{\G}  \ar
 @/^1.5pc/[lll]^{\mathcal{H} } }
\ \ \ \text{and} \ \ \ \xymatrix@C=0.5cm{
\Mod{\bigl(\begin{smallmatrix}
 A & _AN_B \\
 0 & B
\end{smallmatrix}\bigr)} \ar[rrr]^{\F'} &&& \Mod{\Lambda_{(0,0)}} \ar @/_1.5pc/[lll]_{\G'}  \ar
 @/^1.5pc/[lll]^{\mathcal{H}' } }
\]
are adjoint triples.
\begin{proof}
We leave to the reader to check that the above functors provide adjoint pairs. However, we give another direct proof following \cite[Theorem $4.6$]{GP}. Since the bimodule homomorphisms $\phi$ and $\psi$ are zero, there are surjective ring homomorphisms $\Lambda_{(0,0)}\lxr \bigl(\begin{smallmatrix}
 A & 0 \\
 _BM_A & B
\end{smallmatrix}\bigr)$ and $\Lambda_{(0,0)}\lxr \bigl(\begin{smallmatrix}
 A & _AN_B \\
 0 & B
\end{smallmatrix}\bigr)$, defined in an obvious way. Then it is well known that the restriction functors $\F\colon \Mod{\bigl(\begin{smallmatrix}
 A & 0 \\
 _BM_A & B
\end{smallmatrix}\bigr)}\lxr \Mod{\Lambda_{(0,0)}}$ and $\F'\colon \Mod{\bigl(\begin{smallmatrix}
 A & _AN_B \\
 0 & B
\end{smallmatrix}\bigr)}\lxr \Mod{\Lambda_{(0,0)}}$ are fully faithful and have left and right adjoints.
\end{proof}
\end{lem}

We close this subsection with the next result, which shows that the left adjoints of the full embeddings $\F$ and $\F'$ preserve the Gorenstein-projective modules constructed in Theorem~\ref{thmGorproj}. In the following result we keep the assumptions of Theorem~\ref{thmGorproj}.

\begin{prop}
The functors $\G\colon \smod{\Lambda_{(0,0)}}\lxr \smod{\bigl(\begin{smallmatrix}
 A & 0 \\
 _BM_A & B
\end{smallmatrix}\bigr)}$ and $\G'\colon \smod{\Lambda_{(0,0)}}\lxr \smod{\bigl(\begin{smallmatrix}
 A & _AN_B \\
 0 & B
\end{smallmatrix}\bigr)}$ preserve the Gorenstein-projective objects $(\ref{gorprojone})$ and $(\ref{gorprojtwo})$.
\begin{proof}
We show that the functor $\G$ preserves the Gorenstein-projective objects of Theorem~\ref{thmGorproj} $(\alpha)$. The other cases are treated with the same way. Let $(X,Y, f, g)$ be a Gorenstein-projective $\Lambda_{(0,0)}$-module of the form $(\ref{gorprojone})$. This means that we have the following exact sequences$\colon$
\[
\xymatrix{
  0 \ar[r]^{} & N\otimes_BZ \ar[r]^{ \ \ \ s} & X \ar[r]^{\pi_X \ \ \ \ } & \Coker{s} \ar[r] & 0 } \ \ \ \ \xymatrix{
  0 \ar[r]^{} & M\otimes_A\Coker{s} \ar[r]^{ \ \ \ \ \ \ t} & Y \ar[r]^{\pi_Y \ } & Z \ar[r] & 0 }
\]
with $Z\in \Gproj{B}$, $\Coker{s}\in \Gproj{A}$ and $f=(\iden_{M}\otimes \pi_X)\circ t$, $g=(\iden_N\otimes \pi_Y)\circ s$. We claim that the object $\G(X,Y,f,g)=(\Coker{g}, Y, h)$ is Gorenstein-projective. Thus, from \cite[Theorem $1.4$]{Zhang} we have to show that the map $h\colon M\otimes_A\Coker{g}\lxr Y$ is a monomorphism, $\Coker{g}\in \Gproj{A}$ and $\Coker{h}\in \Gproj{B}$. Since the map $\iden_N\otimes \pi_Y\colon N\otimes_B Y\lxr N\otimes_BZ$ is an epimorphism, it follows that $\Coker{g}=\Coker{s}$ and therefore $\Coker{g}$ lies in $\Gproj{A}$. Also,
we have $\pi=\pi_X$ since $g=(\iden_N\otimes \pi_Y)\circ s$. Then, from the relations
$f=(\iden_{M}\otimes \pi_X)\circ t$ and $f=(\iden_{M}\otimes \pi_X)\circ h$ it follows that $h=t$. 
This implies that $h$ is a monomorphism and $\Coker{h}=\Coker{t}=Z$ is Gorenstein-projective. We infer that $\G(X,Y,f,g)$ lies in $\Gproj{\bigl(\begin{smallmatrix}
 A & 0 \\
 _BM_A & B
\end{smallmatrix}\bigr)}$.
\end{proof}
\end{prop}

\section{Homological Embeddings and Gorenstein Artin Algebras}
\label{HomolEmbGorAlgebras}
Our purpose in this section is to provide a method for constructing Morita rings $\Lambda_{(0,0)}=\bigl(\begin{smallmatrix}
A & _AN_B \\
_BM_A & B
\end{smallmatrix}\bigr)$ which are Gorenstein Artin algebras. It turns out that our construction is strongly connected with the property of the functors $\mz_B\colon \Mod{B} \lxr \Mod{\Lambda_{(0,0)}}$ and $\mz_A\colon \Mod{A} \lxr \Mod{\Lambda_{(0,0)}}$ being homological embeddings. This section is divided into two subsections and the main result is stated in the second one.

\subsection{Homological Embeddings}
Let $\Lambda_{(\phi,\psi)}=\bigl(\begin{smallmatrix}
A & _AN_B \\
_BM_A & B
\end{smallmatrix}\bigr)$ be a Morita ring. Associated with the Morita ring $\Lambda_{(\phi,\psi)}$ are the following recollements of abelian categories (see Proposition~\ref{prop:propertiesrecollement})$\colon$
\[
\xymatrix@C=0.2cm{
\Mod{{B/\Image{\phi}}} \ar[rrr]^{\mathsf{I}_B} &&& \Mod{\Lambda_{(\phi,\psi)}} \ar[rrr]^{\mU_A } \ar @/_1.5pc/[lll]_{}  \ar
 @/^1.5pc/[lll]^{} &&& \Mod{A}
\ar @/_1.5pc/[lll]_{\mt_A} \ar
 @/^1.5pc/[lll]_{\mh_A}
 } \ \ \ \ \xymatrix@C=0.2cm{
\Mod{{A/\Image{\psi}}} \ar[rrr]^{\mathsf{I}_A} &&& \Mod{\Lambda_{(\phi,\psi)}} \ar[rrr]^{\mU_B } \ar @/_1.5pc/[lll]_{}  \ar
 @/^1.5pc/[lll]^{} &&& \Mod{B}
\ar @/_1.5pc/[lll]_{\mt_B} \ar
 @/^1.5pc/[lll]_{\mh_B}
 }
\]
Using the idempotent elements $e=\bigl(\begin{smallmatrix}
1_A & 0 \\
0 & 0
\end{smallmatrix}\bigr)$ and $f=\bigl(\begin{smallmatrix}
0 & 0 \\
0 & 1_B
\end{smallmatrix}\bigr)$ of $\Lambda_{(\phi,\psi)}$, we obtain easily that $\Mod{\Lambda/\Lambda e\Lambda}\simeq \Mod{{B/\Image{\phi}}}$,  $\Mod{e\Lambda e}\simeq \Mod{A}$,  $\Mod{\Lambda/\Lambda f\Lambda}\simeq \Mod{{A/\Image{\psi}}}$ and  $\Mod{f\Lambda f}\simeq \Mod{B}$. Note that for simplicity we denote the Morita ring $\Lambda_{(\phi,\psi)}$ by $\Lambda$.

In this subsection we investigate when the ideals $\langle e\rangle=\Lambda e\Lambda$ and $\langle f\rangle=\Lambda f\Lambda$ are stratifying. We recall first the notion of stratifying ideals due to Cline-Parshall-Scott \cite{CPS}.

Let $R$ be a ring and $e$ an idempotent element of $R$. Then we have the exact sequence
\[
\xymatrix{
 0 \ar[r] & \Ker{\mu_R} \ar[r]^{} & Re\otimes_{eRe}eR \ar[r]^{ \ \ \ \ \ \ \mu_R} & R \ar[r]^{} & R/ReR \ar[r] & 0 }
\]
where $\Image{\mu_R}=ReR$ and $\Ker{\mu_R}$ lies in $\Mod{R/ReR}$. The ideal $\langle e\rangle=ReR$ is called {\bf \textsf{stratifying}}, if the following two conditions hold$\colon$
\begin{enumerate}
\item The multiplication map $Re\otimes_{eRe}eR \lxr ReR$ is an isomorphism.

\item $\Tor_{eRe}^i(Re,eR)=0$, for all $i>0$.
\end{enumerate}
The surjective ring homomorphism $R\lxr R/ReR$ induces a fully faithful functor $\mathsf{I}_R\colon \Mod{R/ReR}\lxr \Mod{R}$. Then it is known from \cite{CPS} that the ideal $\langle e\rangle$ is stratifying if and only if the functor $\mathsf{I}_R$ is a {\bf \textsf{homological embedding}} \cite{Psaroud:homolrecol}, i.e. the exact functor $\mathsf{I}_R$ induces an isomorphism between the extension groups$\colon$
\[
\xymatrix{
 {\Ext}_{R/ReR}^n(X,Y) \ \ar[r]^{ \ \  \simeq } & \ {\Ext}_R^n(X,Y) }
\]
for all $X, Y$ in $\Mod{R/ReR}$ and $n\geq 0$. For more details on homological embeddings between abelian categories we refer to \cite{Psaroud:homolrecol}.

We now characterize when the ideals $\langle e\rangle$ and $\langle f\rangle$ are stratifying, or equivalently when the functors $\mathsf{I}_B\colon \Mod{{B/\Image{\phi}}} \lxr \Mod{\Lambda_{(\phi,\psi)}}$ and $\mathsf{I}_A\colon \Mod{{A/\Image{\psi}}} \lxr \Mod{\Lambda_{(\phi,\psi)}}$ are homological embeddings.

\begin{prop}
\label{propstratifyingideals}
Let $\Lambda_{(\phi,\psi)}$ be a Morita ring.
\begin{enumerate}
\item The ideal $\langle e\rangle$ is stratifying if and only if the map $\phi\colon M\otimes_AN\lxr B$ is a monomorphism and $\Tor_{i}^A(M,N)=0$ for all $i>0$.

\item The ideal $\langle f\rangle$ is stratifying if and only if the map $\psi\colon N\otimes_BM\lxr A$ is a monomorphism and $\Tor_{i}^B(N,M)=0$ for all $i>0$.

\end{enumerate}
\begin{proof}
We only prove (i) since part (ii) is dual. For simplicity we write $\Lambda$ for the Morita ring $\Lambda_{(\phi,\psi)}$. An easy computation shows that $f\Lambda e = M$ and $e\Lambda f=N$. Then, since $\Lambda e = e\Lambda e\oplus f\Lambda e$ and $e\Lambda = e\Lambda e\oplus e\Lambda f$, we have the following isomorphisms$\colon$
\[
{\Tor}_{e\Lambda e}^i (\Lambda e, e\Lambda) \iso {\Tor}_{e\Lambda e}^i (e\Lambda e\oplus f\Lambda e, e\Lambda e\oplus e\Lambda f) \iso {\Tor}_{e\Lambda e}^i (f\Lambda e, e\Lambda f) \iso {\Tor}_{A}^i (M, N)
\]
Also, the canonical map $\mu_{\Lambda}\colon \Lambda e\otimes_{e\Lambda e}e\Lambda\lxr \Lambda$ is a monomorphism if and only if the map $f\Lambda e\otimes_{e\Lambda e}e\Lambda \lxr f\Lambda$ is a monomorphism if and only if the map $f\Lambda e\otimes_{e\Lambda e}e\Lambda f\lxr f\Lambda f$ is a monomorphism, i.e. the map $\phi\colon M\otimes_AN\lxr B$ is a monomorphism. Hence, we infer that the ideal $\langle e\rangle$ is stratifying if and only if $\Tor_{i}^A(M,N)=0$ for all $i>0$ and the map $\phi\colon M\otimes_AN\lxr B$ is a monomorphism.
\end{proof}
\end{prop}

We provide examples of Morita rings where the conditions of Proposition~\ref{propstratifyingideals} are satisfied.

\begin{exam}
Let $\Lambda_{(\phi,\psi)}$ be a Morita ring. If $M=0$ we have the upper triangular matrix ring $\Lambda=\bigl(\begin{smallmatrix}
A & _AN_B \\
0 & B
\end{smallmatrix}\bigr)$ and the recollements $(\Mod{B},\Mod{\Lambda},\Mod{A})$ and $(\Mod{A},\Mod{\Lambda},\Mod{B})$, see \cite[Example~$2.12$]{Psaroud:homolrecol}. Then we obtain immediately from Proposition~\ref{propstratifyingideals}, that the functors $\mz_B\colon \Mod{B}\lxr \Mod{\Lambda}$ and $\mz_A\colon \Mod{A}\lxr \Mod{\Lambda}$ are homological embeddings. The same considerations hold when $N=0$. 
\end{exam}

\begin{exam}
Let $\Lambda_{(0,0)}$ be a Morita ring such that $_AN_B$ has an $A$-tight projective $\Lambda_{(0,0)}$-resolution and $_BM_A$ has a $B$-tight projective $\Lambda_{(0,0)}$-resolution, in the sense of \cite{GP}. This means that we have projective resolutions $\cdots \lxr {_AP_1} \lxr {_AP_0} \lxr {_AN} \lxr 0$ and $\cdots \lxr {_BQ_1} \lxr {_BQ_0} \lxr {_BM} \lxr 0$, such that $M\otimes_AP_i=0$ and $N\otimes_BQ_i=0$. Then, if we apply the functor $M\otimes_A-$ to the projective resolution of $N$ we obtain that $M\otimes_AN=0$ and $\Tor_i^A(M,N)=0$ for all $i>0$. Similarly, by applying the functor $N\otimes_B-$ to the projective resolution of $M$, it follows that $N\otimes_BM=0$ and $\Tor_i^B(N,M)=0$ for all $i>0$. Hence, from Propositions~\ref{propstratifyingideals} we infer that the functors $\mz_A\colon \Mod{A}\lxr \Mod{\Lambda_{(0,0)}}$ and $\mz_B\colon \Mod{B}\lxr \Mod{\Lambda_{(0,0)}}$ are homological embeddings. We refer to \cite{GP} for examples of Morita rings with tight resolutions.
\end{exam}

\begin{exam} 
Let $\Lambda$ be an Artin algebra with primitive idempotents $\{e_1, \ldots, e_n\}$. Let $\{S_1, \ldots, S_n\}$ be the corresponding simple $\Lambda$-modules. Assume that $S:=S_1$ is localizable, i.e. $\pd{_\Lambda S}\leq 1$ and $\Ext^{1}_{\Lambda}(S,S)=0$. If we consider the idempotent element $\alpha=e_2+\cdots+e_n$, then it is easy to see that $\alpha(S_1)=0$. This shows that $\add{S}$ is the kernel of the exact functor $\alpha\Lambda\otimes_{\Lambda}-\colon\smod\Lambda\lxr\smod{\alpha\Lambda \alpha}$, in particular the category $\smod{\Lambda/\Lambda\alpha\Lambda}$ is precisely the additive closure $\add{S}$ of $S$. From the short exact sequence $0\lxr \Lambda\alpha\Lambda\lxr \Lambda\lxr \Lambda/\Lambda\alpha\Lambda\lxr 0$ and since $\pd{_\Lambda\Lambda/\Lambda\alpha\Lambda}\leq 1$, it follows that $\Lambda\alpha\Lambda$ is a projective $\Lambda$-module. Then by \cite[Remark $3.2$]{KoenigNagase} we get that $\alpha\Lambda$ is a projective left $\alpha\Lambda\alpha$-module and the map $\Lambda\alpha\otimes_{\alpha\Lambda\alpha}\alpha\Lambda\lxr \Lambda\alpha\Lambda$ is an isomorphism. This implies that the map $e_1\Lambda\alpha\otimes_{\alpha\Lambda \alpha}\alpha\Lambda e_1\lxr e_1\Lambda e_1$ is a monomorphism and $\Tor_{\alpha\Lambda\alpha}^i(e_1\Lambda\alpha,\alpha\Lambda e_1)=0$ for all $i>0$. Note that we view $\Lambda$ as the Morita ring with $A=\alpha\Lambda\alpha$, $B=e_1\Lambda e_1$, $N=\alpha\Lambda e_1$ and $M=e_1\Lambda\alpha$, see Example~\ref{examplesMoritarings} (i). Hence, from Proposition~\ref{propstratifyingideals} we infer that the ideal $\Lambda\alpha\Lambda$ is stratifying. The above claim, that $\Lambda\alpha\Lambda$ being projective implies that $\Lambda\alpha\Lambda$ is a stratifying ideal, can be proved in a different way. We refer to \cite[Example $3.14$]{Psaroud:homolrecol} for more details.  
\end{exam}

We restrict now to the case where the bimodule homomorphisms $\phi$ and $\psi$ are zero, that is $\Lambda_{(0,0)}$ is the trivial extension $(A\times B)\ltimes M\oplus N$ and we have the recollements $(\Mod{A}, \Mod{\Lambda_{(0,0)}}, \Mod{B})$ and $(\Mod{B}, \Mod{\Lambda_{(0,0)}}, \Mod{A})$, see Proposition~\ref{prop:propertiesrecollement} and Example~\ref{examplesMoritarings} (iv). The following result, which is due to Beligiannis \cite[Corollary $4.4$]{Bel:RelativeCleft}, shows that under some conditions we can compute the extension groups induced by the
functors $\mz_{A}\colon \Mod{A}\lxr \Mod{\Lambda_{(0,0)}}$ and $\mz_{B}\colon \Mod{B}\lxr \Mod{\Lambda_{(0,0)}}$. 

\begin{lem}
\label{lemextisomBel}
Let $\Lambda_{(0,0)}$ be a Morita ring which is an Artin algebra. Assume that the right modules $M_A$ and $N_B$ are projective.
\begin{enumerate}
\item  For every $A$-modules $X, X'$ and $n\geq 0$ there are the following isomorphisms$\colon$

\smallskip

\begin{enumerate}
\item For $n=0, 1\colon$$\Ext^n_{\Lambda_{(0,0)}}(\mz_A(X),\mz_A(X'))\iso \Ext^n_{A}(X,X')$.

\smallskip

\item For $n=2k\colon$$\Ext^n_{\Lambda_{(0,0)}}(\mz_A(X),\mz_A(X'))\iso \Ext^{2k}_{A}(X,X')\oplus \Ext^{2(k-1)}_{A}(N\otimes_BM\otimes_AX,X')\oplus \Ext^{2(k-2)}_{A}((N\otimes_BM)^{\otimes^2}{\otimes_A}X,X')\oplus \cdots \oplus \Hom_{A}((N\otimes_BM)^{\otimes^k}{\otimes_A}X,X')$.

\smallskip

\item For $n=2k+1\colon$$\Ext^n_{\Lambda_{(0,0)}}(\mz_A(X),\mz_A(X'))\iso \Ext^{2k+1}_{A}(X,X')\oplus \Ext^{2k-1}_{A}(N\otimes_BM\otimes_AX,X')\oplus \Ext^{2k-3}_{A}((N\otimes_BM)^{\otimes^2}{\otimes_A}X,X')\oplus \cdots \oplus \Ext_{A}^1((N\otimes_BM)^{\otimes^k}{\otimes_A}X,X')$.
\end{enumerate}

\item For every $B$-modules $Y, Y'$ and $n\geq 0$ there are the following isomorphisms$\colon$

\smallskip

\begin{enumerate}
\item For $n=0, 1\colon$$\Ext^n_{\Lambda_{(0,0)}}(\mz_B(Y),\mz_B(Y'))\iso \Ext^n_{B}(Y,Y')$.

\smallskip

\item For $n=2k\colon$$\Ext^n_{\Lambda_{(0,0)}}(\mz_B(Y),\mz_B(Y'))\iso \Ext^{2k}_{B}(Y,Y')\oplus \Ext^{2(k-1)}_{B}(M\otimes_AN\otimes_BY,Y')\oplus \Ext^{2(k-2)}_{B}((M\otimes_AN)^{\otimes^2}{\otimes_B}Y,Y')\oplus \cdots \oplus \Hom_{B}((M\otimes_AN)^{\otimes^k}{\otimes_B}Y,Y')$.

\smallskip

\item For $n=2k+1\colon$$\Ext^n_{\Lambda_{(0,0)}}(\mz_B(Y),\mz_B(Y'))\iso \Ext^{2k+1}_{B}(Y,Y')\oplus \Ext^{2k-1}_{B}(M\otimes_AN\otimes_BY,Y')\oplus \Ext^{2k-3}_{B}((M\otimes_AN)^{\otimes^2}{\otimes_B}Y,Y')\oplus \cdots \oplus \Ext_{B}^1((M\otimes_AN)^{\otimes^k}{\otimes_B}Y,Y')$.
\end{enumerate}
\end{enumerate}
\begin{proof}
We only sketch the proof of (ii), statement (i) follows similarly. From Proposition~\ref{prop:propertiesrecollement} the functor $\mz_B\colon \Mod{B}\lxr \Mod{\Lambda_{(0,0)}}$ is fully faithful and from \cite[Remark $3.7$]{Psaroud:homolrecol} we always have the isomorphism $\Ext^1_{B}(Y,Y')\iso \Ext^1_{\Lambda_{(0,0)}}(\mz_B(Y),\mz_B(Y'))$ for all $B$-modules $Y$ and $Y'$. We explain now how we obtain the rest isomorphisms from \cite[Corollary $4.4$]{Bel:RelativeCleft}. First, from Example~\ref{examplesMoritarings} (iv), the Morita ring $\Lambda_{(0,0)}$ is isomorphic to the trivial extension ring $(A\times B)\ltimes M\oplus N$. Then, the module category $\Mod{\Lambda_{(0,0)}}$ is equivalent to the trivial extension of abelian categories $(\Mod{A}\times \Mod{B})\ltimes H$, where $H$ is the endofunctor $H\colon \Mod{A}\times \Mod{B}\lxr \Mod{A}\times \Mod{B}$, $H(X,Y)=(N\otimes_BY,M\otimes_AX)$. We refer to \cite{FGR} for more details on trivial extensions of abelian categories. We compute only $\Ext^2_{\Lambda_{(0,0)}}(\mz_B(Y),\mz_B(Y'))$. Using the description of $\Mod{\Lambda_{(0,0)}}$ as a trivial extension and \cite[Corollary $4.4$]{Bel:RelativeCleft}, it follows that $\Ext^2_{\Lambda_{(0,0)}}(\mz_B(Y),\mz_B(Y'))$ is isomorphic with the direct sum $\oplus_{i=0}^2\Ext^i_{A\times B}(H^{2-i}(0,Y),(0,Y'))$. The latter extension group is isomorphic with $\Ext^2_B(Y,Y')\oplus \Hom_B(M\otimes_AN\otimes_BY,Y')$, since $\Ext^1_{A\times B}(H(0,Y),(0,Y'))=\Ext^1_{A\times B}((N\otimes_BY,0),(0,Y'))=0$. Hence, $\Ext^2_{\Lambda_{(0,0)}}(\mz_B(Y),\mz_B(Y'))\iso \Ext^2_B(Y,Y')\oplus \Hom_B(M\otimes_AN\otimes_BY,Y')$. The rest isomorphisms follow in the same way, the details are left to the reader.
\end{proof}
\end{lem} 

As a consequence of Lemma~\ref{lemextisomBel} we have the next result. Note that it also follows from Proposition~\ref{propstratifyingideals}. 

\begin{cor}
\label{corhomemb}
Let $\Lambda_{(0,0)}$ be a Morita ring such that the modules $M_A$ and $N_B$ are projective modules.
\begin{enumerate}
\item The following are equivalent$\colon$
\begin{enumerate}
\item The functor $\mz_{A}\colon \Mod{A}\lxr \Mod{\Lambda_{(0,0)}}$ is a homological embedding.

\item $N\otimes_BM=0$.

\end{enumerate}

\item The following are equivalent$\colon$
\begin{enumerate}
\item The functor $\mz_{B}\colon \Mod{B}\lxr \Mod{\Lambda_{(0,0)}}$ is a homological embedding.

\item $M\otimes_AN=0$.

\end{enumerate}

\end{enumerate}
\begin{proof}
(i) (a) $\Longrightarrow$ (b) If the functor $\mz_A$ is a homological embedding, then from Lemma~\ref{lemextisomBel} (i) we get that $\Hom_{A}(N{\otimes_B}M{\otimes_A}X,X')=0$ for every $A$-module $X$ and $X'$. We infer that $N{\otimes_B}M=0$.

(b) $\Longrightarrow$ (a) If $N{\otimes_B}M=0$, then from Lemma~\ref{lemextisomBel} (i) it follows that $\Ext_{A}^n(X,X')\iso \Ext^{n}_{\Lambda_{(0,0)}}(\mz_A(X),\mz_A(X'))$ for every $A$-module $X$, $X'$ and $n\geq 0$. 

(ii) This follows as in (i) using Lemma~\ref{lemextisomBel} (ii).
\end{proof}
\end{cor}

The next result provides another reason
for our interest in stratifying ideals. It is a consequence of Proposition~\ref{propstratifyingideals} and the well known result of Cline-Parshal-Scott \cite{CPS} which relates stratifying ideals and recollements of derived module categories. For the notion of recollement of triangulated categories see \cite{BBD}, and for more details on deriving recollements of abelian categories we refer to \cite{PsaroudVitoria:2}.

\begin{cor}
Let $\Lambda_{(\phi,\psi)}$ be a Morita ring.
\begin{enumerate}
\item If the map $\phi\colon M\otimes_AN\lxr B$ is a monomorphism and $\Tor_{i}^A(M,N)=0$ for all $i>0$, then we have the following recollement of derived categories$\colon$
\[
\xymatrix@C=0.5cm{
\mathsf{D}(\Mod{{B/\Image{\phi}}}) \ar[rrr]^{\mathsf{D}(\mathsf{I}_A)} &&& \mathsf{D}(\Mod{\Lambda_{(\phi,\psi)}}) \ar[rrr]^{\mathsf{D}(\mU_A)} \ar @/_1.5pc/[lll]_{}  \ar
 @/^1.5pc/[lll]^{} &&& \mathsf{D}(\Mod{A})
\ar @/_1.5pc/[lll]_{} \ar
 @/^1.5pc/[lll]^{}
 }
\]

\item If the map $\psi\colon N\otimes_BM\lxr A$ is a monomorphism and $\Tor_{i}^B(N,M)=0$ for all $i>0$, then we have the following recollement of derived categories$\colon$
\[
\xymatrix@C=0.5cm{
\mathsf{D}(\Mod{{A/\Image{\psi}}}) \ar[rrr]^{\mathsf{D}(\mathsf{I}_B)} &&& \mathsf{D}(\Mod{\Lambda_{(\phi,\psi)}}) \ar[rrr]^{\mathsf{D}(\mU_B) } \ar @/_1.5pc/[lll]_{}  \ar
 @/^1.5pc/[lll]^{} &&& \mathsf{D}(\Mod{B})
\ar @/_1.5pc/[lll]_{} \ar
 @/^1.5pc/[lll]^{}
 }
\]
\end{enumerate}
\end{cor}

\subsection{Gorenstein Algebras}
Our aim in this subsection is to provide sufficient conditions for Morita rings with  zero bimodule homomorphisms to be Gorenstein Artin algebras. Recall from \cite{AR:cm, Ha} that an Artin algebra $\Lambda$ is called {\bf \textsf{Gorenstein}} if $\id{_\Lambda \Lambda}<\infty$ and $\id\Lambda_{\Lambda}<\infty$. Equivalently, $\Lambda$ is Gorenstein if and only if $\spli{\Lambda}=\sup\{\pd{_\Lambda I} \ | \ I\in \inj{\Lambda}\}<\infty$ and $\silp{\Lambda}=\sup\{\id{_\Lambda P} \ | \ P\in \proj{\Lambda}\}<\infty$, i.e. $\smod{\Lambda}$ is a Gorenstein abelian category in the sense of \cite{BR}.

We start with the next result which, under some conditions, gives isomorphisms between the extension groups induced from the adjoint pairs $(\mt_{A},\mU_{A})$ and $(\mt_{B},\mU_{B})$. It follows from \cite[Theorem~$3.10$]{Psaroud:homolrecol}, but for completeness we give a direct proof.

\begin{lem}
\label{lemextisom}
Let $\Lambda_{(\phi,\psi)}$ be a Morita ring. Let $X$ be an $A$-module and let $Y$ be a $B$-module.
\begin{enumerate}
\item Assume that the module $M_A$ is projective. Then for every $\Lambda_{(\phi,\psi)}$-module $(X',Y',f',g')$ and $n\geq 0$ we have an isomorphism$\colon$
\[
\xymatrix{
 {\Ext}_{\Lambda_{(\phi,\psi)}}^{n}(\mt_{A}(X),(X',Y',f',g')) \ \ar[r]^{ \ \ \ \ \ \ \ \ \ \ \ \iso} & \ {\Ext}_{A}^{n}(X,X') }
\]

\item Assume that the module $N_B$ is projective. Then for every $\Lambda_{(\phi,\psi)}$-module $(X',Y',f',g')$ and $n\geq 0$ we have an isomorphism$\colon$
\[
\xymatrix{
 {\Ext}_{\Lambda_{(\phi,\psi)}}^{n}(\mt_{B}(Y),(X',Y',f',g')) \ \ar[r]^{ \ \ \ \ \ \ \ \ \ \ \ \iso} & \ {\Ext}_{B}^{n}(Y,Y') }
\]
\end{enumerate}
\begin{proof}
(i) Let $X$ be an $A$-module and let $\cdots \lxr P_1 \lxr P_0 \lxr X \lxr 0$ be a projective resolution of $X$. Since the functor $M\otimes_A-$ is exact, it follows from Proposition~\ref {prop:projmod} and Remark~\ref{remMoritarings} (iv) that the sequence $\cdots \lxr \mt_{A}(P_{1}) \lxr \mt_{A}(P_0) \lxr \mt_{A}(X) \lxr 0$ is a projective resolution of $\mt_{A}(X)$. Let $(X',Y',f',g')$ be a $\Lambda_{(\phi,\psi)}$-module. Then, using the adjoint pair $(\mt_{A},\mU_{A})$ we have the following commutative diagram$\colon$
\[
\xymatrix{
  (\mt_A(X),(X',Y',f',g')) \ar[d]_{\iso} \ \ar@{>->}[r]^{} & (\mt_A(P_0),(X',Y',f',g')) \ar[d]_{\iso}
  \ar[r]^{} & (\mt_A(P_1),(X',Y',f',g')) \ar[d]_{\iso} \ar[r]^{} & \cdots  \\
  \Hom_{A}(X,X') \ \ar@{>->}[r]^{} & \Hom_{A}(P_0,X') \ar[r]^{} & \Hom_{A}(P_1,X') \ar[r]^{} & \cdots   }
\]
This implies that ${\Ext}_{\Lambda_(\phi,\psi)}^n(\mt_{A}(X),(X',Y',f',g'))\iso
{\Ext}_{A}^n(X,X')$ for every $n\geq 0$.

(ii) This follows similarly as in (i).
\end{proof}
\end{lem}

\begin{lem}
\label{leminjectivedim}
Let $\Lambda_{(\phi,\psi)}$ be a Morita ring which is an Artin algebra.
\begin{enumerate}
\item Assume that $M_A$ and $N_B$ are projective modules. If $\id{_{\Lambda_{(\phi,\psi)}}\Lambda_{(\phi,\psi)}}<\infty$, then$\colon$
\[
\left\{
  \begin{array}{lll}
   \id{_{A}A}<\infty, \ \ \id{_{B}B}<\infty.  & \hbox{} \\
           & \hbox{} \\
   \id{_{A}N}<\infty, \ \ \id{_{B}M}<\infty.  & \hbox{}
  \end{array}
\right.
\]
\item Assume that $_BM$ and $_AN$ are projective modules. If $\id{{\Lambda_{(\phi,\psi)}}_{\Lambda_{(\phi,\psi)}}} <\infty$, then$\colon$
\[
\left\{
  \begin{array}{lll}
   \id{A_{A}}<\infty, \ \ \id{B_{B}}<\infty.  & \hbox{} \\
           & \hbox{} \\
   \id{N_{B}}<\infty, \ \ \id{M_{A}}<\infty.  & \hbox{}
  \end{array}
\right.
\]
\end{enumerate}
\begin{proof}
(i) From Proposition~\ref{prop:projmod} we have $\id{_{\Lambda_{(\phi,\psi)}}\mt_{A}(A)}<\infty$ and $\id{_{\Lambda_{(\phi,\psi)}}\mt_{B}(B)}<\infty$. Then, from Lemma~\ref{lemextisom} (i) we have the following isomorphisms for every $A$-module $X$ and $n\geq 0\colon$
\[
{\Ext}^n_{\Lambda_{(\phi,\psi)}}(\mt_{A}(X),\mt_{A}(A))\iso {\Ext}^n_{A}(X,A) \ \ \ \ \text{and} \ \ \ \ {\Ext}^n_{\Lambda_{(\phi,\psi)}}(\mt_{A}(X),\mt_{B}(B))\iso {\Ext}^n_{A}(X,N)
\]
These isomorphisms imply that $\id{_{A}A}\leq \id{_{\Lambda_{(\phi,\psi)}}\mt_{A}(A)}<\infty$ and $\id{_{A}N}\leq \id{_{\Lambda_{(\phi,\psi)}}\mt_{B}(B)}<\infty$. Similarly, for every $B$-module $Y$ and $n\geq 0$ we have from Lemma~\ref{lemextisom} (ii) the following isomorphisms$\colon$
\[
{\Ext}^n_{\Lambda_{(\phi,\psi)}}(\mt_{B}(Y),\mt_{A}(A))\iso {\Ext}^n_{B}(Y,M) \ \ \ \ \text{and} \ \ \ \ {\Ext}^n_{\Lambda_{(\phi,\psi)}}(\mt_{B}(Y),\mt_{B}(B))\iso {\Ext}^n_{B}(Y,B)
\]
Hence, $\id{_{B}M}<\infty$ and $\id{_{B}B}<\infty$.

(ii) In this part we use right modules. If $X$ is a right $A$-module, then $\mt_{A}(X)=(X,X\otimes_AN,\iden_{X\otimes N},\Psi_X)$ and since $_AN$ is projective it follows that the functor $-\otimes_AN\colon \smod{A}\lxr \smod{B}$ is exact and therefore $\mt_{A}$ is exact. Similarly, since $_BM$ is projective we obtain that the functor $\mt_B$ is exact. Then using the isomorphisms of Lemma~\ref{lemextisom}, but for right modules now, the result follows as in case (i). 
\end{proof}
\end{lem}

It should be clear from the proof of the above result, that the assumption of $M_A$, resp. $N_B$, being projective, implies that $\id{_{A}A}<\infty$ and $\id{_{A}N}<\infty$, resp. $\id{_{B}B}<\infty$ and $\id{_{B}M}<\infty$. The same separation property also holds for part (ii). We continue with the next consequence of Lemma~\ref{leminjectivedim}.

\begin{cor}
\label{corsuffcondGoren}
Let $\Lambda_{(\phi,\psi)}$ be a Morita ring which is a Gorenstein Artin algebra.
\begin{enumerate}
\item If $M_A$ is a projective right $A$-module and $_AN$ is a projective left $A$-module, then the algebra $A$ is Gorenstein.

\item If $N_B$ is a projective right $B$-module and $_BM$ is a projective left $B$-module, then the algebra $B$ is Gorenstein.
\end{enumerate}
\end{cor}

Let $\Lambda$ be an Artin algebra and consider the Morita ring $\Delta_{(\phi,\phi)}=\bigl(\begin{smallmatrix}
\Lambda & \Lambda \\
\Lambda & \Lambda
\end{smallmatrix}\bigr)$. If $\Delta_{(\phi,\phi)}$ is Gorenstein, then by Corollary~\ref{corsuffcondGoren} the algebra $\Lambda$ is also Gorenstein. 
We mention that this was observed in Proposition~\ref{proppropertiesDelta} (i), where the the converse also holds in this case. Hence, Corollary~\ref{corsuffcondGoren} generalizes the one direction of Proposition~\ref{proppropertiesDelta} (i). We  give an example to show that the conditions in Corollary~\ref{corsuffcondGoren} are only sufficient.

\begin{exam} 
Let $\Lambda$ be a bimodule $d$-Calabi-Yau noetherian algebra over a field $k$, where $d\geq 2$ is an integer. Let $e$ be a non-trivial idempotent element of $\Lambda$ such that $\Lambda/\Lambda e\Lambda$ is a finite dimensional $k$-algebra. By ~\cite[Theorem $2.2$]{AIR} and it's proof, the algebra $e\Lambda e$ is Gorenstein and the $e\Lambda e$-module $e\Lambda$ is a non-projective Gorenstein-projectve. Note that from ~\cite[Proposition $2.4$]{AIR} the algebra $\Lambda$ has finite global dimension and therefore $\Lambda$ is Gorenstein. 
\end{exam}

In the rest of the subsection our aim is to consider the converse of Corollary~\ref{corsuffcondGoren}, that is how the Gorensteinness of $A$ and $B$ should be inherited to the whole Morita ring. We first need the following preliminary result. As usual we denote by $\du\colon \smod{\Lambda}\lxr \smod{\Lambda^{\op}}$ the duality for Artin algebras, see \cite{ARS}.

\begin{lem}\label{lempdid}
Let $\Lambda_{(0,0)}=\bigl(\begin{smallmatrix}
A & _AN_B \\
_BM_A & B
\end{smallmatrix}\bigr)$ be a Morita ring which is an Artin algebra.
\begin{enumerate}
\item Assume that $\pd{_AN}<\infty$. If $\id{_{\Lambda_{(0,0)}}\mz_A(A)}<\infty$ then $\id{_{\Lambda_{(0,0)}}\mz_A(N)}<\infty$.

\item Assume that $\pd{_BM}<\infty$. If $\id{_{\Lambda_{(0,0)}}\mz_B(B)}<\infty$ then $\id{_{\Lambda_{(0,0)}}\mz_B(M)}<\infty$.

\item Assume that $\pd{M_A}<\infty$. If $\pd{_{\Lambda_{(0,0)}}\mz_A(\du(A))}<\infty$ then $\pd{_{\Lambda_{(0,0)}}\mz_A\big(\Hom_{B}(M,\du(B))\big)}<\infty$.

\item Assume that $\pd{N_B}<\infty$. If $\pd{_{\Lambda_{(0,0)}}\mz_B(\du(B))}<\infty$ then $\pd{_{\Lambda_{(0,0)}}\mz_B\big(\Hom_{A}(N,\du(A))\big)}<\infty$.
\end{enumerate}
\begin{proof}
(i) Let $0\lxr P_n\lxr \cdots \lxr P_0\lxr {_AN}\lxr 0$ be a finite projective resolution of $N$. Then, since the functor $\mz_{A}\colon \smod{A}\lxr \smod{\Lambda_{(0,0)}}$ is exact (see Proposition~\ref{prop:propertiesrecollement}) and $\id{_{\Lambda_{(0,0)}}\mz_A(A)}<\infty$, it follows that $\id{_{\Lambda_{(0,0)}}\mz_A(N)}<\infty$. The proof of part (ii) is dual.

(iii) Since the projective dimension of $M_A$ is finite if and only if the injective dimension of $_A\du(M)$ is finite and we have an isomorphism $_A{\Hom_{B}(_BM_A,_B\du(B))}\iso {_A{\du(M)}}$, then  the result follows by applying the exact functor $\mz_{A}$ to a finite injective coresolution of $_A\du(M)$. The proof of part (iv) is dual.
\end{proof}
\end{lem}

The following is the main result of this section which provides sufficient conditions for Morita rings $\Lambda_{(0,0)}$ with zero bimodule homomorphisms such that $\silp{\Lambda_{(0,0)}}<\infty$ and $\spli{\Lambda_{(0,0)}}<\infty$. This result constitutes the second part of Theorem A presented in the Introduction.

\begin{thm}
\label{thmGorenstein}
Let $\Lambda_{(0,0)}=\bigl(\begin{smallmatrix}
A & _AN_B \\
_BM_A & B
\end{smallmatrix}\bigr)$ be a Morita ring which is an Artin algebra.
\begin{enumerate}
\item Assume the following conditions$\colon$
\begin{enumerate}
\item $M_A$ is projective and $\pd{_BM}<\infty$.

\item $N_B$ is projective and $\pd{_AN}<\infty$.

\item The functors $\mz_{A}$, $\mz_{B}$ are homological embeddings.
\end{enumerate}
If $\silp{A}<\infty$ and $\silp{B}<\infty$, then $\silp{\Lambda_{(0,0)}}<\infty$.

\item Assume the following conditions$\colon$
\begin{enumerate}
\item $_BM$ is projective and $\pd{M_A}<\infty$.

\item $_AN$ is projective and $\pd{N_B}<\infty$.

\item The functors $\mz_{A}$, $\mz_{B}$ are homological embeddings.
\end{enumerate}
If $\spli{A}<\infty$ and $\spli{B}<\infty$, then $\spli{\Lambda_{(0,0)}}<\infty$.
\end{enumerate}
\begin{proof}
(i) From Proposition~\ref{prop:projmod}, it is enough to consider the projective $\Lambda_{(0,0)}$-modules $\mt_{A}(A)$ and $\mt_{B}(B)$. Assume that $\pd{_BM}=\kappa<\infty$ and $\pd{_AN}=\lambda<\infty$. From Lemma~\ref{lem:exactseq}, we have the following exact sequences in $\smod{\Lambda_{(0,0)}}\colon$
\begin{equation}
\label{eqsesone}
\xymatrix{
  0 \ar[r]_{}^{} & \mz_{B}(M) \ar[r]^{} & \mt_{A}(A) \ar[r]^{} & \mz_{A}(A) \ar[r] & 0 }
\end{equation}
and
\begin{equation}
\label{eqsestwo}
\xymatrix{
  0 \ar[r]_{}^{} & \mz_{A}(N) \ar[r]^{} & \mt_{B}(B) \ar[r]^{} & \mz_{B}(B) \ar[r] & 0 }
\end{equation}
Thus, from Lemma~\ref{lempdid} we have to show that $\id{_{\Lambda_{(0,0)}}\mz_A(A)}<\infty$ and $\id{_{\Lambda_{(0,0)}}\mz_B(B)}<\infty$. We first prove that $\id{_{\Lambda_{(0,0)}}\mz_B(B)}<\infty$. Let $(X,Y,f,g)$ be a $\Lambda_{(0,0)}$-module. Then, from the morphism $(\iden_{X},f)\colon \mt_A(X)\lxr (X,Y,f,g)$, we derive the following exact sequences in $\smod{\Lambda_{(0,0)}}\colon$
\begin{equation}\label{eqfirst}
\xymatrix{
  0 \ar[r]_{}^{} & \mz_{B}(\Ker{f}) \ar[r]^{} & \mt_{A}(X) \ar[r]^{} & (X,\Image{f},k',0) \ar[r] & 0 }
\end{equation}
and
\begin{equation}\label{eqsecond}
\xymatrix{
  0 \ar[r]_{}^{} & (X,\Image{f},k',0) \ar[r]^{} & (X,Y,f,g) \ar[r]^{} & \mz_{B}(\Coker{f}) \ar[r] & 0 }
\end{equation}
Applying the functor $\Hom_{\Lambda_{(0,0)}}(-,\mz_B(B))$ to the exact sequence (\ref{eqfirst}), we obtain the following long exact $\Ext$-sequence$\colon$
\[
\xymatrix@C=0.3cm{
  \cdots \ar[r]_{}^{} & \Ext^n_{\Lambda_{(0,0)}}\big((X,\Image{f},k',0),\mz_{B}(B)\big) \ar[r]^{} & \Ext^n_{\Lambda_{(0,0)}}\big(\mt_{A}(X),\mz_{B}(B)\big) \ar[r]^{} & \Ext^n_{\Lambda_{(0,0)}}\big(\mz_{B}(\Ker{f}),\mz_{B}(B)\big) \ar[r] & \cdots }
\]
From Lemma~\ref{lemextisom} (i) it follows that $\Ext^n_{\Lambda_{(0,0)}}(\mt_{A}(X),\mz_{B}(B))=0$ for every $n\geq 0$. Let $\silp{B}=\mu<\infty$. Since the functor $\mz_{B}\colon \smod{B}\lxr \smod{\Lambda_{(0,0)}}$ is a homological embedding, we have $\Ext^n_{\Lambda_{(0,0)}}(\mz_{B}(\Ker{f}),\mz_{B}(B))=0$ for every $n\geq \mu+1$. We infer that $\Ext^n_{\Lambda_{(0,0)}}((X,\Image{f},k',0),\mz_{B}(B))=0$ for every $n\geq \mu+2$. Then from the following long exact sequence$\colon$
\[
\xymatrix@C=0.15cm{
  \cdots \ar[r]_{}^{} & \Ext^n_{\Lambda_{(0,0)}}\big(\mz_B(\Coker{f}),\mz_{B}(B)\big) \ar[r]^{} & \Ext^n_{\Lambda_{(0,0)}}\big((X,Y,f,g),\mz_{B}(B)\big) \ar[r]^{} & \Ext^n_{\Lambda_{(0,0)}}\big((X,\Image{f},k',0),\mz_{B}(B)\big) \ar[r] & \cdots }
\]
obtained from (\ref{eqsecond}), it follows that $\Ext^n_{\Lambda_{(0,0)}}((X,Y,f,g),\mz_{B}(B))=0$ for every $n\geq \mu+2$. Hence we have $\id{_{\Lambda_{(0,0)}}\mz_{B}(B)}\leq \mu+1$ and therefore from Lemma~\ref{lempdid} we infer that $\id{_{\Lambda_{(0,0)}}\mz_{B}(M)}\leq \kappa+\mu+1$. Next, for the injective dimension of $\mz_{A}(A)$, we consider the following exact sequence$\colon$
\begin{equation}
\label{eqthird}
\xymatrix{
  0 \ar[r]_{}^{} & \mz_{A}(\Ker{g}) \ar[r]^{} & \mt_{B}(Y) \ar[r]^{(g,\iden_{Y}) \ \ \ } & (X,Y,f,g) \ar[r] & \mz_{A}(\Coker{g}) \ar[r] & 0 }
\end{equation}
where $\Image{(g,\iden_{Y})}=(\Image{g},Y,0,l')$. Let $\silp{A}=\nu<\infty$. Then, applying the functor $\Hom_{\Lambda_{(0,0)}}(-,\mz_A(A))$ to the two short exact sequences obtained from $(\ref{eqthird})$, we derive as above that $\id{_{\Lambda_{(0,0)}}\mz_{A}(A)}\leq \nu+1$. Note that now we use Lemma~\ref{lemextisom} (ii) and that the functor $\mz_{A}\colon \smod{A}\lxr \smod{\Lambda_{(0,0)}}$ is a homological embedding. Since $\pd{_AN}=\lambda<\infty$, it follows from Lemma~\ref{lempdid} that $\id{_{\Lambda_{(0,0)}}\mz_{A}(N)}\leq \lambda+\nu+1$. Hence, from the exact sequences $(\ref{eqsesone})$ and $(\ref{eqsestwo})$ we have $\id{_{\Lambda_{(0,0)}}\mt_A(A)}\leq \max\{\kappa+\mu, \nu\}+1$ and $\id{_{\Lambda_{(0,0)}}\mt_B(B)}\leq \max\{\lambda+\nu, \mu\}+1$. We infer that $\silp{\Lambda_{(0,0)}}<\infty$.

(ii) This part follows by dual arguments but for completeness we sketch the proof. First, from Proposition~\ref{prop:projmod} it is enough to consider the injective $\Lambda_{(0,0)}$-modules $\mh_{A}(\du(A))$ and $\mh_{B}(\du(B))$. Then, from Lemma~\ref{lem:exactseq} we have the exact sequences in $\smod{\Lambda_{(0,0)}}\colon$ $0 \lxr \mz_{A}(\du(A)) \lxr \mh_{A}(\du(A)) \lxr \mz_{B}(\Hom_A(N,\du(A))) \lxr 0$
and $0 \lxr \mz_{B}(\du(B)) \lxr \mh_{B}(\du(B)) \lxr \mz_{A}(\Hom_B(M,\du(B))) \lxr 0$. Thus, from Lemma~\ref{lempdid} we have to show that $\pd{_{\Lambda_{(0,0)}}\mz_A(\du(A))}<\infty$ and $\pd{_{\Lambda_{(0,0)}}\mz_B(\du(B))}<\infty$. Also, for any $\Lambda_{(0,0)}$-module $(X,Y,f,g)$ we obtain, from the units of the adjoint pairs $(\mU_A,\mh_{A})$ and $(\mU_B,\mh_B)$, the exact sequences$\colon$
$0 \lxr \mz_{A}(\Ker{\pi(f)}) \lxr (X,Y,f,g) \lxr \mh_B(Y) \lxr \mz_{A}(\Coker{\pi(f)}) \lxr 0$
and $ 0 \lxr \mz_{B}(\Ker{\rho(g)}) \lxr (X,Y,f,g) \lxr \mh_A(X) \lxr \mz_{B}(\Coker{\rho(g)}) \lxr 0$.
Then, similarly with part (i) we show that $\spli{\Lambda_{(0,0)}}<\infty$. The details are left to the reader.
\end{proof}
\end{thm}

As a consequence we have the next result on the finiteness of the global dimension of $\Lambda_{(0,0)}$.

\begin{cor} 
\label{corgorensteinfiniteglobal}
Let $\Lambda_{(0,0)}=\bigl(\begin{smallmatrix}
A & _AN_B \\
_BM_A & B
\end{smallmatrix}\bigr)$ be a Morita ring which is an Artin algebra such that the modules $M_{A}$, $N_{B}$ are projective and the functors $\mz_{A}$, $\mz_{B}$ are homological embeddings. If $\gld{A}<\infty$ and $\gld{B}<\infty$, then $\gld{\Lambda_{(0,0)}}<\infty$.
\begin{proof} 
By the proof of Theorem~\ref{thmGorenstein} we have that $\id{_{\Lambda_{(0,0)}}\mz_{A}(N)}<\infty$ and $\id{_{\Lambda_{(0,0)}}\mz_{B}(M)}<\infty$. Since $\spli\Lambda_{(0,0)}<\infty$ it follows that $\pd{_{\Lambda_{(0,0)}}\mz_{A}(N)}<\infty$ and
$\pd{_{\Lambda_{(0,0)}}\mz_{B}(M)}<\infty$. Using that $\Lambda_{(0,0)}$ is the trivial extension ring $(A\times B)\ltimes M\oplus N$ (Example~\ref{examplesMoritarings} (iv)) and \cite[Proposition $5.19$]{GP}, we infer that the global dimension of $\Lambda_{(0,0)}$ is finite. 
\end{proof}
\end{cor}

We continue with the following result which gives us a class of Morita rings, in particular a class of trivial extension rings (Example~\ref{examplesMoritarings} (iv)), where Theorem~\ref{thmGorenstein} can be applied. 

\begin{cor}
\label{corGorensteinMorita}
Let $\Lambda_{(0,0)}=\bigl(\begin{smallmatrix}
A & _AN_A \\
_AN_A & A
\end{smallmatrix}\bigr)$ be a Morita ring which is an Artin algebra. Assume the following conditions$\colon$
\begin{enumerate}
\item $N_A$ and $_AN$ are projective.

\item $N\otimes_AN=0$.
\end{enumerate}
If $A$ is Gorenstein then the ring $\Lambda_{(0,0)}$ is Gorenstein. In particular, if $\gld{A}<\infty$ then $\gld{\Lambda_{(0,0)}}<\infty$.
\begin{proof}
The condition $N\otimes_AN=0$ implies that the functor $\mz_{A}\colon \smod{A}\lxr \smod{\Lambda_{(0,0)}}$ is a homological embedding, see Corollary~\ref{corhomemb}. The second functor we have to check that is a homological embedding is $\mz_{A}'\colon \smod{A}\lxr \smod{\Lambda_{(0,0)}}$ given by $\mz_{A}'(X)=(0,X,0,0)$. Note that this is the functor $\mz_2$ in the notation of subsection~\ref{subsectiondoulemorphism}. Exactly in the same way with Remark~\ref{remdoublemorphismcat} (iii), we show that the two recollements of $\smod{\Lambda_{(0,0)}}$ are equivalent. Both of them are of the form $(\smod{A}, \smod{\Lambda_{(0,0)}}, \smod{A})$, see Proposition~\ref{prop:propertiesrecollement}. Using this equivalence it follows that $\mz_{A}$ is a homological embedding if and only if $\mz'_{A}$ is a homological embedding. Then the result follows from Theorem~\ref{thmGorenstein} and Corollary~\ref{corgorensteinfiniteglobal}. 
\end{proof}
\end{cor}

The above method for constructing Gorenstein algebras is illustrated in the next example.

\begin{exam}
Let $A$ be a finite dimensional Gorenstein $k$-algebra, where $k$ is a field, and let $e$ and $f$ be two idempotents elements of $A$ such that $fAe=0$. Consider the $A$-$A$-bimodule $N:=Ae\otimes_kfA$. Then it follows easily that $N\otimes_AN=0$ and therefore from Corollary~\ref{corGorensteinMorita} we get the Gorenstein algebra$\colon$
\[
\Lambda_{(0,0)}=\begin{pmatrix}
A & _AN_A \\
_AN_A & A
\end{pmatrix}
\]
Note that $\Lambda_{(0,0)}$ is the trivial extension algebra $(A\times A)\ltimes N\oplus N$, see Example~\ref{examplesMoritarings} (iv).
\end{exam}

We close this section with an example of a Morita ring which is a Gorenstein algebra and the conditions of Theorem~\ref{thmGorenstein} are not satisfied. 

\begin{exam}Let $A$ be a ring and $M$ be a right $A$-module. Then from Example~\ref{examplesMoritarings} (ii) we have the Morita ring 
\[
\Lambda_{(\phi,\psi)} = \begin{pmatrix}
B & _BM_A \\
M^{*} & A
\end{pmatrix}
\] 
where $B=\End_{A}(M)$ and $M^{*}={_A{\Hom}_{A}(M,A)}_B$. Note that this Morita ring is  the  Auslander context, in the sense of Buchweitz \cite{Buchweitz}, defined by the pair $(A,M)$. If $M_{A}$ is a finitely generated projective right $A$-module, then from~\cite[Proposition $2.6$, Corollary $1.10$]{Buchweitz} it follows that the rings $A$ and $\Lambda_{(\phi,\psi)}$ are Morita equivalent, and therefore $A$ is Gorenstein if and only if $\Lambda_{(\phi,\psi)}$ is Gorenstein. Hence, if $A$ is a Gorenstein algebra and $M_A$ is a finitely generated projective module, then the Morita ring $\Lambda_{(\phi,\psi)}$ is Gorenstein. By Example~\ref{examplesMoritarings} (ii), the bimodule homomorphisms of this Morita ring are not zero, and also the rest assumptions of Theorem~\ref{thmGorenstein} are not satisfied in general.
\end{exam}

\section{Gorenstein Subcategories and Coherent Functors}
\label{SectionGorsubcatcoherentfun}

In this section we study the monomorphism category $\mono(\Lambda)$, see equation $(\ref{eqmono})$ in subsection~\ref{subMonocat}. In particular, we investigate the full subcategory $\C$ of $\mono(\Lambda)$ consisting of all monomorphisms $f\colon X\lxr Y$ such that the projective dimension of $X$ is finite. In the first subsection, we show that it is a Gorenstein subcategory of $\mono(\Lambda)$ when $\Lambda$ is a Gorenstein Artin algebra. In the second subsection, we prove that the category of coherent functors over the stable category of $\C$ is a Gorenstein abelian category.

\subsection{The Gorenstein subcategory of $\mono(\Lambda)$}
\label{subsectionGorensteinsubcat}
Let $\A$ be an abelian category with enough projective and injective objects and let $n$ be a non-negative integer. Recall from \cite[Theorem $2.2$, Chapter VII]{BR} that
$\A$ is $n$-Gorenstein if and only if every object has Gorenstein-projective dimension at most $n$.
For our purpose, we need the notion of a Gorenstein subcategory but now in the context of exact categories (see subsection~\ref{subMonocat}). Before that, we define Gorenstein-projective objects for exact categories.

\begin{defn}
Let $\A=(\A,\E)$ be an exact category with enough projective objects. An object $X$ in $\A$ is called {\bf \textsf{Gorenstein-projective}} if
there is an $\E$-acyclic complex of projective objects in $\A\colon$
\[
\mathsf{P}^{\bullet}\colon \
\xymatrix{
  \cdots \ar[r]_{}^{} & P^{-1} \ar[r]^{} & P^0 \ar[d]_{\kappa} \ar[r]^{d^0} & P^{1} \ar[r] & P^2 \ar[r] & \cdots \\
  & & X \ar[ur]_{\lambda} & & }
\]
such that $\Hom_{\A}(\textsf{P}^{\bullet},P)$ is acyclic for every object $P$ in $\Proj\A$ and $d^{0}=\lambda\circ \kappa$, where $\kappa\colon P^{0}\lxr X$ is a deflation and $\lambda\colon X\lxr P^{1}$ is an inflation. We denote by $\GProj{\A}$ the full subcategory of Gorenstein-projective objects of $\A$.
\end{defn}

For a complex being acyclic in an exact category we refer to\cite[Definition 10.1]{Buhler}. From now on, when we write $\A$ for an exact category we fix a class $\E$ of exact pairs.

\begin{defn}
\label{defnGorsubcat}
Let $\A$ be an exact category with enough projective objects. Then $\A$ is {\bf $n$\textsf{-Gorenstein}} for some non-negative integer $n$ if every object has Gorenstein-projective dimension at most $n$. Let $\B$ be an exact subcategory of $\A$. We call $\B$ an {\bf $n$\textsf{-Gorenstein subcategory}} of $\A$, if for all $X$ in $\B$ there exists an exact sequence $0\lxr G_{n}\lxr \cdots \lxr G_{0}\lxr X\lxr 0$ in $\B$ such that  $G_{j}\in \GProj\A$ for all $0\leq j\leq n$.
\end{defn}

Let $\Lambda$ be an Artin algebra and $\mathsf{T}_{2}(\Lambda)=\bigl(\begin{smallmatrix}
\Lambda & 0\\
\Lambda & \Lambda
\end{smallmatrix}\bigr)$ be the lower triangular matrix algebra. Let $\Delta_{(0, 0)}=\bigl(\begin{smallmatrix}
\Lambda & \Lambda\\
\Lambda & \Lambda
\end{smallmatrix}\bigr)$ be the Morita ring which is an Artin algebra, see Example~\ref{examplesMoritarings} (iv) and Example~\ref{examrecollementofDelta}. The following results provides a description of the Gorenstein-projective objects of $\mono(\Lambda)$. Recall that the singularity category $\mathsf{D}_{\mathsf{sg}}(\Lambda)$ of $\Lambda$ is defined to be the Verdier quotient $\mathsf{D}^{\mathsf{b}}(\smod\Lambda)/\mathsf{K}^{\mathsf{b}}(\proj\Lambda)$, see \cite{Buchweitz:unpublished, Or}.

\begin{prop}
\label{mono}
Let $\Lambda$ be an Artin algebra and set $\mathcal{X}:=\mono(\Lambda)\cap
\Gproj{\Delta_{(0,0)}}$. The following hold.
\begin{enumerate}

\item $\Gproj(\mono(\Lambda))=\{(X, Y, f, 0)\in \mono(\Lambda) \ | \ 0\lxr X\stackrel{f}{\lxr} Y \stackrel{p}{\lxr} \Coker{f} \lxr 0  \ \text{is an exact}$  $ \ \text{sequence sequence with terms in} \ \Gproj\Lambda \}$.

\item Assume that $\Lambda$ is Gorenstein. Then $\Gproj(\mono(\Lambda))$ is a full subcategory of $\mathcal{X}$ and there is a triangle equivalence $\mD_{\mathsf{sg}}(\mathsf{T}_2(\Lambda))\simeq\underline{\Gproj}(\mono(\Lambda))$.
\end{enumerate}
\begin{proof}
(i) Let $(X, Y, f, 0)$ be an object in $\Gproj(\mono(\Lambda))$. We show that $X$, $Y$ and $\Coker{f}$ are Gorenstein-projective $\Lambda$-modules. Since we have the short exact sequence $0\lxr X\lxr Y\lxr \Coker{f}\lxr 0$ and $\Gproj{\Lambda}$ is closed under extensions, it suffices to show that $X$ and $\Coker{f}$ are Gorenstein-projective. 

We first show that $\Coker{f}$ is Gorenstein-projective. From Lemma~\ref{lemprojinjMono} there is a totally acyclic complex 
\[
\mathsf{P}^{\bullet}\colon \
\xymatrix{
  \cdots \ar[r]_{}^{} & \mt_1(P^{-1})\oplus \mz_2(Q^{-1}) \ar[r]^{ \ d^{-1}} & \mt_1(P^{0})\oplus \mz_2(Q^{0}) \ar[r]^{} & \mt_1(P^{1})\oplus \mz_2(Q^{1}) \ar[r] & \cdots }
\]
where $(X,Y,f,0)=\Coker(d^{-1})$ and $\mt_1(P^i)$, $\mz_2(Q^i)$ belong to 
$\proj(\mono(\Lambda))$. By Lemma~\ref{adjointtriples} the cokernel functor $\Cok\colon \mono(\Lambda)\lxr \smod\Lambda$ is exact and preserves projectives, so $\Cok(\mathsf{P}^{\bullet})$ is an acyclic complex of projective $\Lambda$-modules and $\Cok(X,Y,f,0)=\Coker{f}$. Also, by Lemma~\ref{adjointtriples} we have the isomorphism $\Hom_{\Lambda}(\Cok(\mathsf{P}^{\bullet}), \Lambda)\iso \Hom_{\mono(\Lambda)}(\mathsf{P}^{\bullet}, \mz_2(\Lambda))$ and $\mz_2(\Lambda)$ is in $\proj(\mono(\Lambda))$ by Lemma~\ref{lemprojinjMono}. This implies that the complex $\Hom_{\Lambda}(\Cok(\mathsf{P}^{\bullet}), \Lambda)$ is acyclic and therefore the complex $\Cok(\mathsf{P}^{\bullet})$ is a totally acyclic complex of projective $\Lambda$-modules. Hence, the $\Lambda$-module $\Coker{f}$ is Gorenstein-projective. 

We claim that the functor $\mz_2\colon \smod{\Lambda}\lxr \mono(\Lambda)$ lifts Gorenstein-projectives to $\mono(\Lambda)$, i.e. $\mz_2(\Coker{f})$ lies in $\Gproj(\mono(\Lambda))$. Let $\mathsf{Q}^{\bullet}\colon \cdots \lxr Q^{-1}\lxr Q^0\lxr Q^1\lxr \cdots$ a totally acyclic complex of projective $\Lambda$-modules and $Y=\Coker(Q^{-1}\lxr Q^0)$. Clearly, from Lemma~\ref{lemprojinjMono} and Lemma~\ref{adjointtriples} it follows that $\mz_2(\mathsf{Q}^{\bullet})$ is an acyclic complex of projective objects in $\mono(\Lambda)$. Consider now the complexes $\Hom_{\mono(\Lambda)}(\mz_2(\mathsf{Q}^{\bullet}),\mt_1(\Lambda))$ and $\Hom_{\mono(\Lambda)}(\mz_2(\mathsf{Q}^{\bullet}),\mz_2(\Lambda))$. Since the complex $\Hom_{\Lambda}(\mathsf{Q}^{\bullet}, \Lambda)$ is acyclic,  using the adjoint pair $(\mU_2,\mt_1)$ and that the functor $\mz_2$ is fully faithful (Lemma~\ref{adjointtriples}), we obtain that the complexes $\Hom_{\mono(\Lambda)}(\mz_2(\mathsf{Q}^{\bullet}),\mt_1(\Lambda))$ and $\Hom_{\mono(\Lambda)}(\mz_2(\mathsf{Q}^{\bullet}),\mz_2(\Lambda))$ are acyclic. This shows that $\mz_2(\mathsf{Q}^{\bullet})$ is a totally acyclic complex. We infer that the object $\mz_2(Y)$ lies in $\Gproj(\mono(\Lambda))$. 

Consider the following exact sequence$\colon$
\begin{equation}
\label{exactsequence}
\xymatrix{
0 \ar[r]^{}& \mt_1(X) \ar[r]^{(\iden_{X},f) \ \ } & (X, Y, f, 0) \ar[r]^{(0,p) \ \ } & \mz_2(\Coker{f}) \ar[r] & 0}
\end{equation}
Since $\Gproj(\mono(\Lambda))$ is closed under kernels of epimorphisms, it follows that the object $\mt_1(X)$ lies in $\Gproj(\mono(\Lambda))$. We claim that $X$ is in $\Gproj\Lambda$. The object $\mt_1(X)$ is the cokernel of a totally acyclic complex $\mathsf{P}^{\bullet}$ as above. Applying the functor $\mU_2\colon \mono(\Lambda)\lxr \smod\Lambda$, we get from Lemma~\ref{corliftinggorproj} that $\mU_2(\mathsf{P}^{\bullet})$ is an acyclic complex of projective $\Lambda$-modules. Then using the adjunction $(\mU_2,\mt_1)$ it follows easily that $\mU_2(\mathsf{P}^{\bullet})$ is a totally acyclic complex, see also Proposition~\ref{proptotallyacyclic} and Corollary~\ref{corGorprojDelta}. Hence, the module $X$ is Gorenstein-projective. This completes the proof that $X$, $Y$ and $\Coker{f}$ lie in $\Gproj\Lambda$.

Conversely, let $(X, Y, f, 0)$ be an object in $\mono(\Lambda)$
such that $X, Y$ and $\Coker{f}$ lie in $\Gproj{\Lambda}$. Then, from Corollary~\ref{corGorprojDelta} we get that the object $\mt_1(X)$ lies in $\Gproj(\mono(\Lambda))$. Also, from the first part of the proof, the object $\mz_2(\Coker{f})$ belongs to $\Gproj(\mono(\Lambda))$. Since $\Gproj(\mono(\Lambda))$ is closed under extensions, we obtain from the exact sequence $(\ref{exactsequence})$ that $(X, Y, f, 0)$ lies in $\Gproj(\mono(\Lambda))$.

(ii) Since $\Lambda$ is Gorenstein, it follows from Proposition~\ref{proppropertiesDelta} that $\Delta_{(0,0)}$ is a Gorenstein algebra and that $\Gproj{\Delta_{(0,0)}}=\{(X, Y, f, g)$ $\in \smod{\Delta_{(0,0)}} \ | \ X, Y \in \Gproj{\Lambda}\}$.
 This means that $\mathcal{X}$ consists of all objects of the form
 $(X, Y, f,0)$, where $f\colon X\lxr Y$ is a monomorphism and $X, Y$ lie in $\Gproj{\Lambda}$. Thus $\mathcal{X}$ contains $\Gproj(\mono(\Lambda))$ by (i).

Set $\mathsf{S}_{2}(\Gproj\Lambda):=\{(X, Y, f)\in \smod\mathsf{T}_2(\Lambda) \ | \ 0\lxr X\lxr Y \lxr \Coker{f} \lxr 0  \ \text{exact in} \ \Gproj\Lambda\}$ and define the functor $F\colon \mathsf{S}_{2}(\Gproj\Lambda)\lxr \Gproj(\mono(\Lambda))$  by $F(X, Y, f)=(X, Y, f, 0)$ on objects $(X,Y,f)\in \mathsf{S}_{2}(\Gproj\Lambda)$. Then it is straightforward to check that $F$ is an equivalence of categories.
Since $\Lambda$ is Gorenstein, we have from \cite[Theorem $1.1$]{LiZhang} that $\Gproj\mathsf{T}_2(\Lambda)=\mathsf{S}_{2}(\Gproj\Lambda)$. This implies that the categories $\Gproj\mathsf{T}_2(\Lambda)$ and $\Gproj(\mono(\Lambda))$ are equivalent. Since the algebra $\mathsf{T}_2(\Lambda)$ is Gorenstein by \cite{Ha}, it follows from \cite{Buchweitz:unpublished} that there is a triangle equivalence between the singularity category $\mD_{\mathsf{sg}}(\mathsf{T}_2(\Lambda))$ and the stable category $\underline{\Gproj}\mathsf{T}_2(\Lambda)$. We infer that $\mD_{\mathsf{sg}}(\mathsf{T}_2(\Lambda))$ and $\underline{\Gproj}(\mono(\Lambda))$ are triangle equivalent.
\end{proof}
\end{prop}

As a consequence we obtain the next result due to Xiao-Wu Chen.

\begin{cor}\textnormal{\cite[Theorem 4.1]{Chen}}
Let $\Lambda$ be a selfinjective algebra. Then there is a triangle equivalence between the stable monomorphism category $\underline{\mono}(\Lambda)$ and the singularity category $\mD_{\mathsf{sg}}(\mathsf{T}_2(\Lambda))$.
\begin{proof}
From Proposition~\ref{mono} (i) and since $\Lambda$ is selfinjective, it follows that the category $\underline{\Gproj}(\mono(\Lambda))$ is precisely $\underline{\mono}(\Lambda)$. Then the result follows immediately from Proposition~\ref{mono} (ii).
\end{proof}
\end{cor}

Consider now the following subcategory of $\mono(\Lambda)\colon$
\begin{equation}
\label{subcofmonopdfinite}
\C:=\big\{(X, Y, f,0)\in \mono(\Lambda) \ | \ \pd{_{\Lambda}X}<\infty \big\}.
\end{equation}
We denote by $\mathscr{P}^{<\infty}(\Lambda)$ the full subcategory of $\smod\Lambda$ consisting of all $\Lambda$-modules of finite projective dimension.
Then $\mathscr{P}^{<\infty}(\Lambda)$ is an exact subcategory of $\smod{\Lambda}$, since it is extension closed, and this implies that $\C$ is also an exact subcategory of $\mono(\Lambda)$. The first main result on the structure of $\C$ is as follows. This result constitutes the first part of Theorem B presented in the Introduction.

\begin{thm}
\label{Gro}
 Let $\Lambda$ be an $n$-Gorenstein algebra for some non-negative integer $n$. Then $\C$ is an $n$-Gorenstein subcategory of $\mono(\Lambda)$.
\begin{proof}

Let $(X, Y, f, 0)$ be an object in $\C$ and consider the exact sequence $(\ref{exactsequence})$ in $\mono(\Lambda)$. Since $(\mU_2,\mt_1)$ is an adjoint pair of exact functors and both functors preserve projective objects (Lemma~\ref{adjointtriples}), we have the isomorpism $\Ext_{\mono(\Lambda)}^{i}((G_{1}, G_{2}, f, 0),\mt_{1}(X))\cong \Ext_{\Lambda}^{i}(\mU_2(G_{1}, G_{2}, f, 0), X)=\Ext_{\Lambda}^{i}(G_{2}, X)$ for all $i\geq 1$ and $(G_{1}, G_{2}, f, 0)$ in $\Gproj(\mono(\Lambda))$. By Proposition~\ref{mono} the $\Lambda$-module $G_2$ is Gorenstein-projective and since $\pd{_{\Lambda}X}<\infty$, it follows that $\Ext_{\Lambda}^{i}(G_{2}, X)=0$ for all $i\geq 1$ (recall that $(\Gproj{\Lambda}, \mathscr{P}^{<\infty}(\Lambda))$ is a cotorsion pair in $\smod{\Lambda}$, see \cite{BR}). Hence, $\Ext_{\mono(\Lambda)}^{i}((G_{1}, G_{2}, f, 0), \mt_{1}(X))=0$ for all $i\geq 1$ and $(G_{1}, G_{2}, f, 0)$ in $\Gproj(\mono(\Lambda))$. This implies that $(\ref{exactsequence})$ remains exact after applying $\Hom_{\mono(\Lambda)}((G_{1}, G_{2}, f, 0),-)$, for every $(G_{1}, G_{2}, f,$ $ 0)$ in $\Gproj(\mono(\Lambda))$. Since the algebra $\Lambda$ is $n$-Gorenstein, there exist the following two exact sequences of left $\Lambda$-modules$\colon$
\[
\xymatrix{
0 \ar[r]^{} & P_n \ar[r]^{a_n} & \cdots \ar[r] & P_1\ar[r]^{a_1} & P_0\ar[r]^{a_0} & X\ar[r] & 0 }
\]
and
\[
\xymatrix{
0 \ar[r]^{} & G_n \ar[r]^{b_n} & \cdots \ar[r] & G_1 \ar[r]^{b_1} & G_0\ar[r]^{b_0 \ \ \ } & \Coker{f}\ar[r] & 0 }
\]
where $P_{j}$ and $G_{j}$ are Gorenstein-projective $\Lambda$-modules for all $0\leq j\leq n$. Applying the exact functors $\mt_1$ and $\mz_2$, respectively, we get the exact sequences in $\C\colon$
\[
\xymatrix{
0 \ar[r]^{} & \mt_1(P_n) \ar[r]^{ \ \ \mt_1(a_n)} & \cdots \ar[r] &  \mt_1(P_1)\ar[r]^{\mt_1(a_1)} & \mt_1(P_0)\ar[r]^{\mt_1(a_0)} & \mt_1(X) \ar[r] & 0 }
\]
and
\[
\xymatrix{
0 \ar[r]^{} & \mz_2(G_n) \ar[r]^{ \ \ \mz_2(b_n)} & \cdots \ar[r] & \mz_2(G_1) \ar[r]^{\mz_2(b_1)} & \mz_2(G_0) \ar[r]^{\mz_2(b_0) \ \ \ } & \mz_2(\Coker{f})\ar[r] & 0 }
\]
where $\mt_{1}(P_j)$ and $\mz_{2}(G_j)$ belong to $\Gproj(\mono(\Lambda))$, for all $0\leq j\leq n$, by Proposition~\ref{mono}. Since the map $\Hom_{(\mono(\Lambda))}(\mz_2(G_0),(0,p))$ is surjective, we obtain from the Horseshoe Lemma the following exact commutative diagram$\colon$
\[
\xymatrix{
0 \ar[r]^{}& \mt_1(P_0)\ar[r]^{(1 \ 0)  \ \ \ \ \ \ \  \ } \ar@{->>}[d]^{\mt_1(a_0)} & \mt_1(P_0)\oplus \mz_2(G_0) \ar[d]^{\alpha_0} \ar[r]^{ \ \ \ \ {^t(0 \ 1)} } & \mz_2(G_0) \ar@{->>}[d]^{\mz_2(b_0)} \ar[r] & 0 \\
0 \ar[r]^{}& \mt_1(X)\ar[r]^{(\iden_{X},f) \ \ } & (X, Y, f, 0) \ar[r]^{(0,p)} & \mz_2(\Coker{f}) \ar[r] & 0}
\]
Now taking the exact sequence of the kernels and applying the functor $\Hom_{(\mono(\Lambda))}(\mz_2(G_1),-)$, we obtain that the map $\Hom_{(\mono(\Lambda))}(\mz_2(G_1),\Ker{a_0})\lxr \Hom_{\mono(\Lambda)}(\mz_2(G_1),\mz_2(\Ker{b_0}))$ is surjective. This follows since $\Ext^1_{\mono(\Lambda)}(\mz_2(G_1),\mt_1(\Ker{a_0}))\cong \Ext^1_{\Lambda}(G_1,\Ker{a_0})=0$.  
Then continuing in the same way we construct an exact sequence of $(X,Y,f,0)$ by objects in $\Gproj(\mono(\Lambda))$ of length at most $n$. According to Definition~\ref{defnGorsubcat}, we infer that $\C$ is an $n$-Gorenstein subcategory of $\mono(\Lambda)$.
\end{proof}
\end{thm}

\subsection{Categories of coherent functors and Gorensteinness}\label{subsectioncoherentgor}

In this subsection we continue our study on the subcategory $\C$ of $\mono(\Lambda)$, see $(\ref{subcofmonopdfinite})$.  More precisely, we consider the category of coherent functors $\smod{\underline{\C}}$ over the stable category of $\C$ and show that it is a Gorenstein abelian category. For this purpose, we generalise \cite[Theorem 3.11]{MatsuiTakahashi} from the abelian setting to exact categories. 

In subsection~\ref{subsectionGorensteinsubcat} we recalled the definition of the singularity category of an Artin algebra $\Lambda$. When we deal with an additive category $\A$, for instance an exact (sub)category, the notion of singularity category can be extended using the category of coherent functors over $\A$. This approach was recently investigated by Matsui and Takahashi in \cite{MatsuiTakahashi}. 
We now discuss this. Let $\A$ be an additive category with weak kernels, that is, for each morphism $f\colon X\lxr Y$ in $\A$ there exists a morphism
$g\colon Z\lxr X$ in $\A$  such that the sequence $\Hom_{\A}(-,Z) \lxr \Hom_{\A}(-,X) \lxr \Hom_{\A}(-,Y)$ is exact. We denote by $\smod{\A}$ the category of  coherent functors over $\A$, i.e. functors $F\colon \A^{\op}\lxr \Ab$ such that there is an exact sequence $\Hom_{\A}(-,X)\lxr \Hom_{\A}(-,Y)\lxr F\lxr 0$ with $X$ and $Y$ in $\A$. It is known that $\smod{\A}$ is an abelian category with enough projective objects. We refer to \cite{Auslander:coherent, Auslander:repdim} for more details on coherent functors. Then, Matsui and Takahashi \cite{MatsuiTakahashi} considered the Verdier quotient
\[
\mD_{\mathsf{sg}}(\smod{\A}):=\mathsf{D}^{\mathsf{b}}(\smod{\A})/\mathsf{K}^{\mathsf{b}}(\proj({\smod{\A}})) 
\]
and call it the {\bf \textsf{singularity category}} of $\smod{\A}$. We remark that this triangulated category is included in the general framework of the stabilization of an abelian or exact category studied by Beligiannis \cite{Bel:ABcontexts}.

In what follows, we show that the singularity category of $\mono(\Lambda)$ is trivial. We write $\smod\mono(\Lambda)$ for the category of coherent functors over the monomorphism category $\mono(\Lambda)$.

\begin{prop}
\label{proptrivialsingmono}
Let $\Lambda$ be an Artin algebra.
Then the following hold
\begin{enumerate}
\item The category $\smod{\mono(\Lambda)}$ is abelian.

\item We have$\colon$$\gld(\smod({\smod\Lambda}))\leq \gld(\smod\mono(\Lambda))\leq 2$.

\item The singularity category $\mD_{\mathsf{sg}}(\smod\mono(\Lambda))$ is trivial.
\end{enumerate}
\begin{proof} (i) Since $\mono(\Lambda)$ is closed under kernels by Lemma~\ref{adjointtriples}, it follows that $\mono(\Lambda)$ has weak kernels. Hence, the category of coherent functors $\smod(\mono(\Lambda))$ is abelian.

(ii) Let $F$ be a functor in $\smod\mono(\Lambda)$, that is, there is an exact sequence$\colon$
\[
\xymatrix@C=0.5cm{
\Hom_{\mono(\Lambda)}(-,(X_1,Y_1,f_1,0)) \ar[rr]^{(-,(a,b)) \ } && \Hom_{\mono(\Lambda)}(-,(X_0,Y_0,f_0,0)) \ar[r]^{} & F \ar[r] & 0 }
\]
where $(X_1,Y_1,f_1,0)$ and $(X_0,Y_0,f_0,0)$ are objects in $\mono(\Lambda)$.
Since we have the exact sequence
\[
\xymatrix@C=0.3cm{
0 \ar[rr] && (\Ker{a},\Ker{b},k,0) \ar[rr] && (X_1,Y_1,f_1,0) \ar[rr]^{(a,b) \ } && (X_0,Y_0,f_0,0), }
\]
and $\Ker{(a,b)}=(\Ker{a},\Ker{b},k,0)$ lies in $\mono(\Lambda)$, we obtain the following exact sequence$\colon$
\[
\xymatrix@C=0.5cm{
0\ar[r] & (-,\Ker{(a,b)}) \ar[r]^{} & (-,(X_1,Y_1,f_1,0)) \ar[r]^{} & (-,(X_0,Y_0,f_0,0)) \ar[r]^{} & F \ar[r] & 0 }
\]
This implies that $\gld(\smod\mono(\Lambda))\leq 2$. From Lemma~\ref{adjointtriples} we know that $(\mt_{1}, \mU_{1})$ is an adjoint pair between $\smod{\Lambda}$ and
$\mono(\Lambda)$ and the functor $\mt_{1}$ is fully faithful. Then \cite[Theorem 3.1]{Xi} yields that $\gld(\smod(\smod\Lambda))\leq \gld(\smod\mono(\Lambda))$. This completes the proof of (ii).

(iii) This statement follows immediately from (ii).
\end{proof}
\end{prop}

Although that the singularity category $\mD_{\mathsf{sg}}(\smod\mono(\Lambda))$ is trivial, we show in Corollary~\ref{corcoherentfunctorgoren} that if we restrict to the subcategory $\C$ of $\mono(\Lambda)$, then this singularity category is not at all trivial. 

Before we get there we need some more definitions.

\begin{defn}
An additive subcategory $\B$ of $\A$ is called {\bf \textsf{quasi-resolving}} if it contains $\Proj{\A}$ and  given a conflation $X\lxr Y\lxr Z$ with $Y$ and $Z$ in $\B$ then the object $X$ lies in $\B$.
A quasi-resolving subcategory $\B$ is called {\bf \textsf{resolving}} if it is closed under direct summands and extensions, i.e. given a conflation $X\lxr Y\lxr Z$ with $X$ and $Z$ in $\B$ then the object $Y$ lies in $\B$.
\end{defn}

Note that a resolving subcategory $\B$ of $\A$ is an exact subcategory of $\A$ since it is closed under extensions. Let $X$ be an object in $\A$. Since $\A$ has enough projective objects there exists a deflation $g\colon P\lxr X$ with $P\in \Proj{\A}$ (or a right $\Proj{\A}$-approximation). This means that there is an exact pair $K\lxr P\lxr X$, where the map $f\colon K\lxr P$ is an inflation. The object $K$ is called the first syzygy of $X$ and is denoted by $\Omega(X)$. The nth syzygy $\Omega^{n}(X)$ of $X$ is defined inductively as $\Omega(\Omega^{n-1}(X))$. We denote by $\Omega^{n}(\A)$ the subcategory of $\A$ consisting of all nth syzygies of objects in $\A$. Assume that there is a left $\Proj{\A}$-approximation $f\colon X\lxr P$, i.e. $f$ is an inflation with $P\in \Proj{\A}$ such that the map  $\Hom_{\A}(f,P')\colon \Hom_{\A}(P,P')\lxr \Hom_{\A}(X,P')$ is surjective for all $P'\in \Proj{\A}$. Then we have the exact pair $X\lxr P\lxr L$, where the map $g\colon P\lxr L$ is a deflation. The object $L$ is called the first cosyzygy
of $X$ and is denoted by $\Omega^{-1}(X)$. The nth cosyzygy $\Omega^{-n}(X)$ of $X$ is defined inductively as $\Omega^{-1}(\Omega^{-(n-1)}(X))$. We denote by $\Omega^{-n}(\A)$ the subcategory of $\A$ consisting of all nth cosyzygies of objects in $\A$.

We are now ready to prove the second main result of this section which generalizes \cite[Theorem 3.11]{MatsuiTakahashi} to the setting of exact categories. 

\begin{thm}
\label{mainthmexactcat}
Let $\A$ be an exact category with enough projective objects. Let $\B$ be a quasi-resolving subcategory of $\A$ such that $\Omega^{n}(\B)\subseteq \GProj\A$ for some non-negative integer $n$ and is closed under $\Omega^{-1}$. Then the following statements hold.
\begin{enumerate}
\item $\smod{\underline{\Omega^n(\B)}}$ is a Frobenius abelian category.

\item $\smod\underline{\B}$ is a $3n$-Gorenstein abelian category.
\end{enumerate}
Moreover, there are the following triangle equivalences$\colon$
\[
\xymatrix{
\mD_{\mathsf{sg}}(\smod\underline{\B}) \ar[r]^{\simeq \ \ } & \underline{\Gproj}(\smod{\underline{\B}}) } \ \ \text{and} \ \ \xymatrix{
\mD_{\mathsf{sg}}(\smod{\underline{\Omega^n(\B)}}) \ar[r]^{ \ \ \simeq } & \underline{\smod} \, \underline{\Omega^n(\B)} }
\]
\begin{proof} We divide the proof into four steps.

\textsf{Step} $1\colon$ We show that $\smod\underline{\B}$ is an abelian category with enough projective objects. It suffices to show that $\underline{\B}$ has weak kernels.  Let $m\colon M\lxr N$ be a morphism in $\B$. Since $\A$ has enough projective objects and $\Proj{\A}\subseteq \B$, there is a deflation $p\colon P\lxr N$ with $P\in \Proj{\A}$. Then, by the axiom $(\textsf{Ex} \, 2)$ of exact categories (see subsection~\ref{subMonocat}), we have the pullback diagram
\[
\xymatrix{
 L \ar[r]^{p'} \ar[d]_{m'} & M  \ar[d]^{m} \\
 P   \ar[r]^{p}    & N }
\]
such that the map $p'$ is also a deflation. From \cite[Proposition A.1]{Keller} and since $\B$ is a quasi-resolving subcategory of $\A$ we obtain the following conflation in $\B\colon$
\[
\xymatrix{ L \ar[r]^{f \ \ \ } & M\oplus P \ar[r]^{ \ \ g} & N }
\]
where $f=\bigl(\begin{smallmatrix}
p' \\
-u'
\end{smallmatrix}\bigr)$ and $g=\bigl(\begin{smallmatrix}
u & p\\
\end{smallmatrix}\bigr)$. Thus, for an object $X$ in $\B$ we have the exact sequence$\colon$
\[
\xymatrix{ \Hom_{\A}(X, L) \ar[r]^{ } & \Hom_{\A}(X, M\oplus P) \ar[r]^{ } & \Hom_{\A}(X, N) }
\]
Let $u\colon X\lxr M\oplus P$ be a morphism in $\A$ such that $\underline{u}$ is in the kernel of $\underline{\Hom}_{\A}(X, g)\colon\underline{\Hom}_{\A}(X, M\oplus P)\lxr \underline{\Hom}_{\A}(X, N)$. Then $u\circ g$ is the composition of some morphisms $a\colon X\lxr Q$ and $b\colon Q\lxr N$ in $\A$, where $Q\in \Proj\A$. There is a morphism $c\colon Q\lxr M\oplus P$ with $c\circ g=b$. So $(a\circ c-u)\circ g=0$. This implies that there is a morphism $d\colon X\lxr L$ such that $a\circ c - u=d\circ f$. We have $\underline{u}=\underline{d\circ f}$, which is in the image of $\underline{\Hom}_{\A}(X, f)\colon \underline{\Hom}_{\A}(X, L)\lxr \underline{\Hom}_{\A}(X, M)$. Thus the sequence
\[
\xymatrix{
\underline{\Hom}_{\A}(-, L)|_{\underline{\B}} \ar[r] & \underline{\Hom}_{\A}(-, M)|_{\underline{\B}} \ar[r] & \underline{\Hom}_{\A}(-, N)|_{\underline{\B}} }
\]
is exact in $\smod\underline{\B}$.

\textsf{Step} $2\colon$ We show that for any object $F$ in $\smod\underline{\B}$ there is a conflation
 $A\lxr B\lxr C$ in $\B$ which induces a projective resolution as follows$\colon$
\begin{equation}
\label{projresol}
\xymatrix@C=0.5cm{
   &  &  \cdots\cdots \ar[r] &  \underline{\Hom}_{\A}(-, \Omega^2(A))|_{\underline{\B}} \ar[r] &  \underline{\Hom}_{\A}(-, \Omega^2(B))|_{\underline{\B}} \ar[r] & \underline{\Hom}_{\A}(-, \Omega^2(C))|_{\underline{\B}}
                \ar@{->} `r/8pt[d] `/10pt[lll] `^dl[llll]|{ } `^r/3pt[dllll] [dlll] \\
     &  &   \underline{\Hom}_{\A}(-, \Omega(A))|_{\underline{\B}} \ar[r]   & \underline{\Hom}_{\A}(-, \Omega(B))|_{\underline{\B}}\ar[r] & \underline{\Hom}_{\A}(-, \Omega(C))|_{\underline{\B}} \ar[r] & \underline{\Hom}_{\A}(-, A)|_{\underline{\B}}
                \ar@{->} `r/8pt[d] `/10pt[lll] `^dl[llll]|{} `^r/3pt[dllll] [dlll] \\
     &  &  \underline{\Hom}_{\A}(-, B)|_{\underline{\B}} \ar[r]   & \underline{\Hom}_{\A}(-, C)|_{\underline{\B}} \ar[r] & F \ar[r] &  0
              }
\end{equation}
Let $F$ be a functor in $\smod\underline{\B}$. Then there is an exact sequence $\underline{\Hom}_{\A}(-, B)|_{\underline{\B}}\stackrel{\phi}{\lxr} \underline{\Hom}_{\A}(-, C)|_{\underline{\B}}\lxr F\lxr 0$ with $B, C\in \B$ and by Yoneda's Lemma the map $\phi$ is of the form $\underline{\Hom}_{\A}(-, u)|_{\underline{\B}}$ for some morphism $u\colon B\lxr C$. As in Step $1$, we obtain a conflation $A\lxr B\oplus Q\lxr C$ in $\B$ with $Q$ in $\Proj{\A}$. Since $\A$ has enough projective objects, there is a deflation $P\lxr C$ with $P$ projective in $\A$. Note that any deflation ending at the object $P$ splits. Then we can form the following pullback diagram
\[
\xymatrix{
 & \Omega(C) \ar[d] \ar@{=}[r] & \Omega(C) \ar[d]  \\
 A \ar[r]^{} \ar@{=}[d]_{} & A\oplus P  \ar[d]^{} \ar[r] & P\ar[d] \\
A \ar[r] & B\oplus Q   \ar[r]^{}  & C }
\]
where every row or column is a conflation. In particular, we get the conflation $\Omega(C)\lxr A\oplus P\lxr B\oplus Q$.  Consider now a deflation $Q'\lxr B$ with $Q'$ in $\Proj{\A}$. From the following commutative diagram
\[
\xymatrix{
 & \Omega(B)\oplus Q \ar[d] \ar@{=}[r] & \Omega(B)\oplus Q \ar[d]  \\
 \Omega(C) \ar[r]^{} \ar@{=}[d]_{} & \Omega(C)\oplus Q'\oplus Q  \ar[d]^{} \ar[r] & Q'\oplus Q\ar[d] \\
\Omega(C) \ar[r] & A\oplus P   \ar[r]^{}  & B\oplus Q }
\]
we obtain the conflation $\Omega(B)\oplus Q\lxr \Omega(C)\oplus P'\lxr A\oplus P$ where $P'=Q'\oplus Q$.  Iterating this procedure yields the conflations$\colon$$\Omega(B)\oplus Q\lxr \Omega(C)\oplus P^{\prime}\lxr A\oplus P$, $\Omega(A)\oplus P\lxr \Omega(B)\oplus P^{\prime\prime}\lxr \Omega(C)\oplus P^{\prime}$, $\Omega^{2}(C)\oplus P^{\prime}\lxr \Omega(A)\oplus P^{\prime\prime\prime}\lxr \Omega(B)\oplus P^{\prime\prime}$ and so on, where $P^{\prime}, P^{\prime\prime}, P^{\prime\prime\prime}$  belong to $\Proj{\A}$. Putting these conflations together and using Step $1$, the desired projective resolution of $F$ follows immediately.

\textsf{Step} $3\colon$ We show that for any object $F$ in $\smod\underline{\B}$ we have $\Ext_{\smod\underline{\B}}^{i}(F, \Proj(\smod\underline{\B}))=0$ for $i>3n$.
Firstly, given a conflation $A\lxr B\lxr C$ in $\B$ such that $\Ext^1_{\A}(C,\Proj{\A})=0$, it follows as in \cite[Lemma 2.2 (2)]{MatsuiTakahashi} that the sequence $\underline{\Hom}_{\A}(C, X)\lxr \underline{\Hom}_{\A}(B, X)\lxr \underline{\Hom}_{\A}(A, X)$ is exact for every $X\in \A$.
Since $\Omega^{j}(C)$ lies in $\Gproj\A$ for all $j\geq n$, we know that $\Ext^{1}_{\A}(\Omega^{j}(C), \Proj\A)=0$ for all $j\geq n$. Thus by the above fact
and Step 2  we obtain that $\underline{\Hom}_{\A}(\Omega^{j}C, Y)\lxr \underline{\Hom}_{\A}(\Omega^{j}B, Y)\lxr \underline{\Hom}_{\A}(\Omega^{j}A, Y)$ is exact for any $Y\in \B$ and all $j\geq n$. Then, applying the functor $(-,\underline{\Hom}_{\A}(-,Y)|_{\underline{\B}})$ to $(\ref{projresol})$ and using Yoneda's Lemma, we infer that  $\Ext_{\smod\underline{\B}}^{i}(F, \Proj(\smod\underline{\B}))=0$ for $i>3n$.

\textsf{Step} $4\colon$ Let $F$ be an object in $\smod\underline{\B}$. By Step $2$ there is a conflation $A\lxr B\lxr C$ in $\B$ which induces a projective resolution
of $F$ as indicated in diagram $(\ref{projresol})$. Set $G:=\Coker(\underline{\Hom}_{\A}(-, \Omega^{n}(B))|_{\underline{\B}}\lxr \underline{\Hom}_{\A}(-, \Omega^{n}(C))|_{\underline{\B}})$. We show that $G$ is a Gorenstein-projective object in $\smod\underline{\B}$. For simplicity, we write $L:=\Omega^{n}(A)$, $M:=\Omega^{n}(B)$ and $N:=\Omega^{n}(C)$.  Then from Step $2$ we get a conflation $L\lxr M\lxr N$. Since $\Omega^{n}(\B)$ is closed under $\Omega^{-1}$, there is a left $\Proj\A$-approximation $L\lxr Q$ with $Q\in \Proj\A$. We make the following pushout diagram$\colon$
\[
\xymatrix{
 L \ar[d] \ar[r] &  Q \ar[d]  \ar[r] & \Omega^{-1}(L) \ar@{=}[d]  \\
 M \ar[r]^{} \ar[d]_{} & N\oplus Q  \ar[d]^{} \ar[r] & \Omega^{-1}(L)  \\
N \ar@{=}[r] & N    &  }
\]
Note that the middle vertical conflation splits, i.e.  $\Ext^{1}_{\A}(N, Q)=0$ since $N\in \GProj\A$ and $Q\in \Proj\A$. Thus we obtain the conflation $M\lxr N\oplus Q\lxr \Omega^{-1}(L)$.
Iterating this procedure gives rise to the conflations$\colon$$N\oplus Q\lxr \Omega^{-1}(L)\oplus Q^{\prime}\lxr \Omega^{-1}M $, \ $\Omega^{-1}(L)\oplus Q^{\prime}\lxr \Omega^{-1}(M)\oplus Q^{\prime\prime}\lxr \Omega^{-1}(N)\oplus Q $, $\Omega^{-1}(M) \oplus Q^{\prime\prime}\lxr \Omega^{-1}(N)\oplus Q^{\prime\prime\prime}\lxr \Omega^{-2}(L)\oplus Q^{\prime}$ and so on, where $Q^{\prime}$, $Q^{\prime\prime}$, $Q^{\prime\prime\prime}$ are  in $\Proj{\A}$. Thus by Step $1$ we obtain an exact sequence$\colon \underline{\Hom}_{\A}(-, L)|_{\underline{\B}}\lxr \underline{\Hom}_{\A}(-, M)|_{\underline{\B}}\lxr\underline{\Hom}_{\A}(-, N)|_{\underline{\B}} \lxr \underline{\Hom}_{\A}(-, \Omega^{-1}(L))|_{\underline{\B}}\lxr \underline{\Hom}_{\A)}(-, \Omega^{-1}(M))|_{\underline{\B}}\lxr \underline{\Hom}_{\A}(-, \Omega^{-1}(N))|_{\underline{\B}}\lxr\cdots$. Combining this with $(\ref{projresol})$ we obtain an exact sequence of projective objects in $\smod{\underline{\B}}$ as follows$\colon$
\begin{equation}
\label{totallyacyclicresol}
\xymatrix@C=0.5cm{
   &  &  \cdots\cdots \ar[r] &  \underline{\Hom}_{\A}(-, \Omega(N))|_{\underline{\B}} \ar[r] &  \underline{\Hom}_{\A}(-, L)|_{\underline{\B}} \ar[r] & \underline{\Hom}_{\A}(-, M)|_{\underline{\B}}
                \ar@{->} `r/8pt[d] `/10pt[lll] `^dl[llll] `^r/3pt[dllll] [dlll] \\
     &  &   \underline{\Hom}_{\A}(-, N)|_{\underline{\B}} \ar[r]^{\alpha \ \ \ \ \ }   & \underline{\Hom}_{\A}(-, \Omega^{-1}(L))|_{\underline{\B}}\ar[r] & \underline{\Hom}_{\A}(-, \Omega^{-1}(M))|_{\underline{\B}} \ar[r] & \cdots
              }
\end{equation}
where $\Image{\alpha}=G$. By the construction of the above pushout diagrams using left $\Proj{\A}$-approximations and since  $\Omega^{n}(\B)\subseteq \GProj\A$, we get that $\Ext^{1}_{\A}(\Omega^{i}(\Omega^{n}(\B)), \Proj\A))=0$ for any $i\in \mathbb{Z}$. Let $Y$ be an object in $\B$. Then applying the functor  $(-,\underline{\Hom}_{\A}(-, Y)|_{\underline{\B}})$ to $(\ref{totallyacyclicresol})$ and using \cite[Lemma 2.2 (2)]{MatsuiTakahashi} as explained in \textsf{Step} $3$, we obtain an acyclic complex$\colon$ $\cdots\lxr \underline{\Hom}_{\A}(\Omega^{-1}(L), Y)\lxr \underline{\Hom}_{\A}(N, Y)\lxr \underline{\Hom}_{\A}(M, Y)\lxr\underline{\Hom}_{\A}(L, Y)\lxr\underline{\Hom}_{\A}(\Omega(N), Y)\lxr\cdots$. This implies that $(\ref{totallyacyclicresol})$ is a totally acyclic complex, equivalently, $G$ is a Gorenstein-projective object in $\smod\underline{\B}$.
Hence, we have shown that for every object $F$ in $\smod\underline{\B}$ there is a projective resolution as in $(\ref{projresol})$ such that the nth syzygy $G$ is Gorenstein-projective. We infer that $\smod\underline{\B}$ is a $3n$-Gorenstein abelian category. Moreover, from \cite[Corollary $4.13$]{Bel:ABcontexts} we obtain the desired triangle equivalence between $\mD_{\mathsf{sg}}(\smod\underline{\B})$ and $\underline{\GProj}(\smod\underline{\B})$.

It remains to show that $\smod{\underline{\Omega^n(\B)}}$ is a Frobenius abelian category, i.e $\smod\underline{\Omega^n(\B)}$ is of Gorenstein dimension at most zero. From \cite[Proposition 3.6]{MatsuiTakahashi}
it suffices to show that the stable category $\underline{\Omega^n(\B)}$ is triangulated. Recall that $\GProj\A$ is an exact Frobenius category and as in the abelian case it follows easily that $\GProj\A$ is extension closed. We claim  that $\Omega^n(\B)$ is an admissible subcategory of $\GProj\A$ (see \cite{Chen:threeresults}), that is, $\Omega^n(\B)$ is an extension closed subcategory of $\GProj\A$ such that for each object $B$ in $\Omega^n(\B)$ there are conflations $B\lxr P\lxr \Omega^{-1}(B)$ and $\Omega(B)\lxr Q\lxr B$ with $P$, $Q$ in $\Proj\A$. Note that $\Omega^n(\B)$ being admissible implies that it is an exact Frobenius category and therefore from\cite{Happelbook} it follows that $\underline{\Omega^n(\B)}$ is triangulated. Since we have $\Proj\A\subseteq \Omega^{n}(\B)\subseteq \GProj\A$ we only have to show that $\Omega^{n}(\B)$ is extension closed. For this, we first show that $\Omega^{n}(\B)=\B\cap \GProj\A$. 
Since $\B$ is quasi-resolving we have $\Omega^{n}(\B)\subseteq\B$. This implies that $\Omega^{n}(\B)\subseteq\B\cap \GProj\A$.
Let $X$ be an object in $\B\cap \GProj\A$. Then $X$ is Gorenstein-projective,  so $X\cong \Omega^{-n}(\Omega^{n}(X))$. Since $X\in \B$ and $\Omega^{n}(\B)$ is closed under $\Omega^{-1}$, we have that $\Omega^{n}(X)\in \Omega^{n}(\B)$ and $\Omega^{-n}(\Omega^{n}(X))\in \Omega^{n}(\B)$. This shows that $X\in \Omega^{n}(\B)$, i.e. $\B\cap \GProj\A\subseteq\Omega^{n}(\B)$. We now show that $\Omega^{n}(\B)=\B\cap \GProj\A$ is extension closed. Consider a conflation $X\lxr Y\lxr Z$ with $X$ and $Z$ in $\B\cap \GProj\A$. Then there is a left-$\Proj\A$ approximation $X\lxr P\lxr \Omega^{-1}(X)$. Taking the pushout diagram of these two conflations and since $\Ext^1(Z,P)=0$, we obtain the conflation $Y\lxr P\oplus Z\lxr \Omega^{-1}(X)$.
The object $P\oplus Z$ lies in $\B$ and the object $\Omega^{-1}(X)$ lies also in $\B$ since we assume that $\Omega^n(B)=\B\cap \GProj\A$ is closed under $\Omega^{-1}$. Since $\B$ is quasi-resolving, it follows that the object $Y$ lies in $\B$. Since $\GProj\A$ is closed under extensions, we conclude that the object $Y$ lies in $\B\cap \GProj(\A)$. This completes the proof that $\smod{\underline{\Omega^n(\B)}}$ is a Frobenius abelian category. Finally, from \cite[Corollary $4.13$]{Bel:ABcontexts} we get that a triangle equivalence between $\mD_{\mathsf{sg}}(\smod{\underline{\Omega^n(\B)}})$ and $\underline{\smod}\underline{\Omega^n(\B)}$.
\end{proof}
\end{thm}

\begin{rem}
\label{remrelationwithMT}
The first part of the proof of Theorem~\ref{mainthmexactcat} is devoted to show that the category of coherent functors  $\smod\underline{\B}$ is abelian. This is similar to \cite[Proposition 2.11(i)]{MatsuiTakahashi}. However, in the setting of exact categories we need to show how we obtain from the axioms the conflation which gives us the correct $\Hom$-exact sequence in order to conclude that $\underline{\B}$ has weak kernels. Part two of our proof is proved in the same way as \cite[Proposition 2.11(ii)]{MatsuiTakahashi}, but again we need to make clear that this construction works in our setting.
Similar comments hold for the rest of the proof. Moreover, as in \cite[Theorem 5.4]{MatsuiTakahashi}, we can deduce a triangle equivalence between $\mD_{\mathsf{sg}}(\smod\underline{\B})$ and $\mD_{\mathsf{sg}}(\smod{\underline{\Omega^n(\B)}})$. We suggest to the reader to work the details following the proof of  \cite[Theorem 5.4]{MatsuiTakahashi}. 
Clearly, the statement of Theorem~\ref{mainthmexactcat} as well as its proof are inspired by the work of Matsui and Takahashi \cite{MatsuiTakahashi}.
\end{rem}

We return to the case of the monomorphism category, but before we apply Theorem~\ref{mainthmexactcat} to the subcategory $\C=\{(X, Y, f,0)\in \mono(\Lambda) \ | \ \pd{_{\Lambda}X}<\infty\}$ of $\mono(\Lambda)$ we need the following result.

\begin{lem}
\label{resolving}
Let $\Lambda$ be an $n$-Gorenstein algebra. The following hold.
\begin{enumerate}
\item $\C$ is a resolving subcategory of $\mono(\Lambda)$.

\item The category $\Omega^{n}(\C)$ is a Frobenius subcategory of $\Gproj(\mono(\Lambda))$ and $\underline{\Omega^{n}(\C)}$ is a triangulated
subcategory of $\underline{\Gproj}(\mono(\Lambda))$.

\item $\Omega^{n}(\C)$ is closed under $\Omega^{-1}$.
\end{enumerate}
\begin{proof}
(i) From Lemma~\ref{lemprojinjMono} the category $\C$ contains the subcategory $\proj(\mono(\Lambda))$. As mentioned already (see before Theorem~\ref{Gro}), the subcategory  $\C$ is closed under extensions, kernels of epimorphisms and direct summands. This means that $\C$ is a resolving subcategory of $\mono(\Lambda)$.

(ii) By Theorem~\ref{Gro} we have the following description of the nth syzygies of $\C\colon$
\[
\Omega^{n}(\C)=\big\{(X,Y,f,0) \ | \ f\colon X\lxr Y \ \text{is a split monomorphism with} \ X\in \proj\Lambda \ \text{and} \ Y\in \Gproj\Lambda \big\}
\]
From Lemma~\ref{lemprojinjMono} and Lemma~\ref{mono}, it follows that $\proj(\mono(\Lambda))\subseteq\Omega^{n}(\C)\subseteq \Gproj(\mono(\Lambda))$. Thus $\Omega^{n}(\C)$ is a Frobenius subcategory of $\Gproj(\mono(\Lambda))$ and  $\underline{\Omega^{n}(\C)}$ is a triangulated subcategory of $\underline{\Gproj}(\mono(\Lambda))$.

(iii) From statement (ii), any object in $\Omega^{n}(\C)$ is of the form  $(P, P\oplus G, \bigl(\begin{smallmatrix}
{\rm Id}_{P}\\
0
\end{smallmatrix}\bigr), 0)$, where  $P\in \proj\Lambda$ and $G\in \Gproj\Lambda$. By defintion of a Gorenstein projective module, there is a monomorphism $u\colon G\lxr Q$ which is a left $\proj\Lambda$-approximation. Then we obtain the following map
\[
0\lxr (P, P\oplus G, \bigl(\begin{smallmatrix}
\iden_{P}\\
0
\end{smallmatrix}\bigr), 0)\xrightarrow{(\iden_{P}, \bigl(\begin{smallmatrix}
\iden_{P} & 0\\
0 & u
\end{smallmatrix}\bigr))} (P, P\oplus Q, \bigl(\begin{smallmatrix}
\iden_{P}\\
0
\end{smallmatrix}\bigr), 0)
\]
and it is easy to check that it is a left $\proj(\mono(\Lambda))$-approximation of $(P, P\oplus G, \bigl(\begin{smallmatrix}
\iden_{P}\\
0
\end{smallmatrix}\bigr), 0)$. Moreover, the cokernel $\Coker(\iden_{P}, \bigl(\begin{smallmatrix}
\iden_{P} & 0\\
0 & u
\end{smallmatrix}\bigr))=(0, \Coker u, 0, 0)$ lies in $\Omega^{n}(\C)$. This completes the proof.
\end{proof}
\end{lem}

\begin{rem}
We showed above that $\C$ is a resolving subcategory of $\mono(\Lambda)$. However, the category $\C$ is not resolving as a subcategory the double morphism category $\DMor(\smod\Lambda)$, i.e. of $\smod{\Delta_{(0,0)}}$. Indeed, by Proposition~\ref{prop:projmod} and for any $P\in \proj\Lambda$ the object $\mt_2(P)=(P,P,0,\iden_{P})$ is projective in $\smod{\Delta_{(0,0)}}$ but is not projective in $\C$. Thus, the category $\C$ doesn't contain all the projectives of $\smod{\Delta_{(0,0)}}$.
\end{rem}

As a consequence of Theorem~\ref{mainthmexactcat}  and Lemma~\ref{resolving} we obtain the following result, which is the second part of Theorem B presented in the Introduction.

\begin{cor}
\label{corcoherentfunctorgoren}
Let $\Lambda$ be an $n$-Gorenstein Artin algebra for some integer $n\geq 0$.
Then for the category of coherent functors over $\underline{\C}$ and $\underline{\Omega^n(\C)}$, respectively, the following hold$\colon$
 \begin{enumerate}

\item $\smod\underline{\C}$ is a $3n$-Gorenstein abelian category.

\item $\smod{\underline{\Omega^n(\C)}}$ is a Frobenius abelian category.
\end{enumerate}
Moreover, there are the following triangle equivalences$\colon$
\[
\xymatrix{
\mD_{\mathsf{sg}}(\smod\underline{\C}) \ar[r]^{\simeq \ \ } & \underline{\Gproj}(\smod{\underline{\C}}) } \ \ \text{and} \ \ \xymatrix{
\mD_{\mathsf{sg}}(\smod{\underline{\Omega^n(\C)}}) \ar[r]^{ \ \ \simeq } & \underline{\smod} \, \underline{\Omega^n(\C)} }
\]
\end{cor}

It should be noted that although the stable derived category $\mD_{\mathsf{sg}}(\smod{\mono(\Lambda)})$ is trivial (Proposition~\ref{proptrivialsingmono}), if we consider the Gorenstein subcategory $\C$ of $\mono(\Lambda)$ then $\mD_{\mathsf{sg}}(\smod\underline{\C})$ is not at all trivial.
This example shows that the singularity category of certain subcategories of exact or abelian categories provide us with interesting triangulated categories. We refer to \cite{MatsuiTakahashi} for more information and further applications of this approach. We close this section with the next remark where we clarify Corollary~\ref{corcoherentfunctorgoren}.

\begin{rem}
\label{finalremark}
Let $\Lambda$ be an Artin algebra and consider the monomorphism category $\mon(\Lambda)$ as a full subcategory of the morphism category $\Mor(\smod{\Lambda})$, i.e. $\smod{\mathsf{T}_2(\Lambda)}$, see subsection~\ref{subMonocat}. Note that the monomorphism category $\mono(\Lambda)$ that we have worked throughout the paper, i.e. viewing it as a full subcategory of the double morphism category $\DMor(\smod\Lambda)$, is equivalent to $\mon(\Lambda)$, see Lemma~\ref{lemequivmono}. Consider the  full subcategory $\D:=\{(X, Y, f) \ | \ f\colon X\lxr Y \ \text{is} \ \text{a} \ \Lambda\mbox{-}\text{monomorphism} \ \text{with} \ \pd{_{\Lambda}X}<\infty \}$ of $\mon(\Lambda)$. Note that it is equivalent with the subcategory $\C$, see $(\ref{subcofmonopdfinite})$, of $\mono(\Lambda)$.  Assume that $\Lambda$ is an $n$-Gorenstein algebra for some non-negative integer $n$. Then for $\D$ the following properties hold$\colon$

(i) $\D$ is a resolving subcategory of $\smod{\mathsf{T}_2(\Lambda)}$. Indeed, $\D$ contains $\proj\mathsf{T}_2(\Lambda)$, see \cite{ARS}. Furthermore, it follows easily that $\D$ is closed under extensions, kernels of epimorphisms and direct summands since  the subcategory $\mathscr{P}^{<\infty}(\Lambda)$  has these properties as well.

(ii) $\Omega^{n}(\D)$ is a Frobenius subcategory of $\Gproj\mathsf{T}_2(\Lambda)$ and is closed under $\Omega^{-1}$. This statement follows exactly as we argued in Lemma~\ref{resolving} for the subcategory $\C$.

From the above conditions, we can apply \cite[Theorem 3.11]{MatsuiTakahashi} to the pair $\D\subseteq \smod\mathsf{T}_2(\Lambda)$ and obtain that $\smod\underline{\D}$ is a $3n$-Gorenstein abelian category and there are equivalences $\mD_{\mathsf{sg}}(\smod\underline{\D})\simeq \underline{\Gproj}(\smod{\underline{\D}})$ and $\mD_{\mathsf{sg}}(\smod\underline{\Omega^n\D})\simeq \underline{\smod}\underline{\Omega^n(\D)}$. Since there is an equivalence between $\C$ and $\D$, Lemma~\ref{lemprojinjMono} implies also an equivalence between $\smod\underline{\C}$ and $\smod\underline{\D}$. Hence, in this way we recover the statements of Corollary~\ref{corcoherentfunctorgoren}. However, we believe that it is more useful to prove Corollary~\ref{corcoherentfunctorgoren} via Theorem~\ref{mainthmexactcat} as we did, since it only requires the pair of subcategories $\C\subseteq \mono(\Lambda)$  which is at the level of exact categories.
\end{rem}

\medskip

\medskip

\begin{ackn}
This work was developed during a stay of the authors to the University of Stuttgart in 2014. Both authors would like to express their gratitude to Steffen Koenig for the  warm hospitality, his support, as well as his comments and suggestions on this project. Moreover, the authors are grateful to Hongxing Chen for many helpful discussions. The first author is supported by the National Natural Science Foundation of China
(Grant No. 11101259) and the second author is supported by the Norwegian Research Council (NFR 221893) under the project \textit{Triangulated categories in Algebra}.
\end{ackn}

\end{document}